\documentclass[12pt,a4paper]{article}

\usepackage{amsmath,amssymb,amsthm,fancyhdr,supertabular,calc,eepic}
\usepackage{mathrsfs}
\usepackage[mathcal]{euscript}
\usepackage[a4paper,nofoot,headheight=14pt,body={15cm,23cm},nofoot,dvips]{geometry}
\usepackage{wasysym}

\usepackage[dvips]{graphicx}
\usepackage[dvips]{color}

\newtheorem*{result}{Result}
\newtheorem{theorem}{Theorem}
\newtheorem{lemma}[theorem]{Lemma}
\newtheorem{proposition}[theorem]{Proposition}

\newtheorem{corollary}[theorem]{Corollary}

\newtheorem*{theorema}{Theorem A}

\newtheorem*{frameworktheorem}{Theorem{\bf \ \ref{theorem:frameworkprojplane}}}

\newtheorem*{hypothesis}{Hypothesis}
\newtheorem*{remark}{Remark}

\newcommand{\Z}{\mathbb{Z}}

\newcommand{\Galph}{G_\alpha}
\newcommand{\Lalph}{L_\alpha}

\newcommand{\ordLalph}{|L_\alpha|}

\newcommand{\curlyc}[1]{\mathcal{C}_#1}
\newcommand{\spaceP}{\mathcal{P}}
\newcommand{\polysum}[1]{q^{#1}+\dots+q+1}
\newcommand{\polysumsq}[1]{q^{#1}+\dots+q^2+1}
\newcommand{\fracNR}{\frac{n_g}{r_g}}

\newcommand{\linel}{\mathfrak{L}}
\newcommand{\hatt}{\hat {\ } }

\begin{document}


\pagestyle{plain}

\pagenumbering{arabic}
\title{Transitive Projective Planes}
\author{Nick Gill\thanks{This paper contains results from the
    author's PhD thesis. I would like to thank my supervisor, Professor Jan Saxl. Professor Bill Kantor has also given much helpful advice for which I am very grateful.}}

\maketitle

\begin{abstract}
A long-standing conjecture is that any transitive finite projective plane is
Desarguesian. We make a contribution towards a proof of this conjecture by
showing that a group acting transitively on the the points of a non-Desarguesian
projective plane must not contain any components.
\end{abstract}

\section{Background definitions and main results}

We say that a projective plane is {\it transitive} (resp. {\it primitive}) if it
admits an automorphism group which is transitive (resp. primitive) on points.
Kantor\cite{kantor} has proved that
a projective plane $\spaceP$ of order $x$ admitting a point-primitive
automorphism group $G$ is Desarguesian and $G\geq PSL(3,x)$, or else
$x^2+x+1$ is a prime and $G$ is a regular or Frobenius group of order
dividing $(x^2+x+1)(x+1)$ or $(x^2+x+1)x$.

Kantor's result, which depends upon the Classification of Finite Simple Groups, represents the strongest success in the pursuit of a proof to the
conjecture mentioned in the abstract. A corollary of Kantor's result is that a
group acts primitively on the points of a projective plane $\spaceP$ if and only if it acts primitively on
the lines of $\spaceP$. We also know, by a combinatorial argument of Block, that
a group acts transitively on the points of a projective plane $\spaceP$ if and
only if it acts transitively on the lines of $\spaceP$\cite{block}.

Our primary result is the following:

\begin{theorema}
Suppose that $G$ acts transitively on a projective plane
$\spaceP$ of order $x$. Then one of the following cases holds:
\begin{itemize}
\item $\spaceP$ is Desarguesian, $G\geq PSL(3,x)$ and the action is
$2$-transitive on points;
\item $G$ does not contain a component. In particular all minimal normal subgroups of $G$ are elementary abelian.
\end{itemize}
\end{theorema}

Here a {\it component} $C$ of a group $G$ is defined to be a subnormal quasi-simple subgroup of $G$. We note that Theorem A implies that if an almost simple group (or almost quasi-simple group) $G$ acts transitively on the lines of a projective plane $\spaceP$ of order $x$ then $\spaceP$ is Desarguesian and $G$ has socle $PSL(3,x).$ Note that definitions for group theory terms used here are provided in Section \ref{section:basic}.

Theorem A also relates to two other results that already exist in the literature. The first is Kantor's result on primitive projective planes \cite{kantor} which has already been mentioned and which is used in the proof of Theorem A; Theorem A can be thought of as a generalization of Kantor's result. The second is Ho's result that a finite projective plane admitting more than one abelian Singer group is Desarguesian \cite[Theorem 1]{ho}; this result is implied by Theorem A and \cite[Lemma 4.3 and Theorem 2]{ho} - details are given in \cite{gill4}. In fact \cite{gill4} outlines a number of results about line-transitive projective planes that follow from Theorem A.

Finally we note that all groups and sets that we consider in this paper are finite.

\section{Overview of Proof}

To prove Theorem A we need to analyse many different possible transitive group
actions on finite projective planes. The framework for our analysis of the transitive projective planes will be given
by results in \cite{campraeg2} and \cite{camina}. The key theorem is the following:

\begin{theorem}\cite[Theorem 2]{camina}
Let $G$ act transitively on a projective plane $\spaceP$ and let $M$
be a minimal normal subgroup of $G$. Then $M$ is either abelian or simple.
\end{theorem}

In fact we are able to state our results more strongly by rewriting
this result in terms of {\it components}. Hence the
theorem which will provide the framework for our analysis is the
following:

\begin{theorem}\label{theorem:frameworkprojplane}
Suppose that $G$ acts transitively on a projective plane
$\spaceP$. Then $G$ contains at most one component.
\end{theorem}

The proof of this theorem, which involves rewriting proofs of similar
theorems from \cite{campraeg2} and \cite{camina}, is given in Section
\ref{section:projplaneframework}. In Section \ref{section:basic} we give the
basic lemmas and notation which will be used throughout the remainder
of the paper.

In the remaining sections we use Theorem \ref{theorem:frameworkprojplane} to
examine the possible unique components of a group $G$ acting transitively on a
projective plane. Existing results in the literature are generally limited to
the case where the component is simple and $G$ is almost simple.

\section{Framework results}\label{section:projplaneframework}

We prove Theorem \ref{theorem:frameworkprojplane}
which states that if a group $G$ acts transitively upon a projective
plane then $G$ contains at most one component.  Our proof of
Theorem \ref{theorem:frameworkprojplane} starts with some preliminary results.

Note first that if $C$ is a component of $G$ then
$C^\circ:= \ <C^g:g\in G> \ \cong C\circ C^{g_1}\circ\dots\circ C^{g_m}$ is
a normal subgroup of $G$ where $g_1,\dots,g_m\in G$; furthermore, if $C$ and
$D$ are components of $G$ with $C$ not $G$-conjugate to $D$ then
$[C,D]=1$ and so $[C^\circ,D^\circ]=1$.

We need some information about the fixed points of automorphisms of
a projective plane $\spaceP$ of order $x$: If an automorphism $g$ fixes
at least $x$ points then $g$ is called {\it quasicentral} and $g$ fixes $x+1$, $x+2$ or
$x+\sqrt{x}+1$ points\cite[4.1.7]{dembov}. In the first two cases $g$
fixes a {\it fan}, namely a line $\linel$ and a point $\alpha$ and all
the points on $\linel$ and all the lines incident with $\alpha$. The
distinction between the two cases depends on whether or not $\alpha$
lies on $\linel$. In the third case the set of fixed points and fixed lines of $g$ forms a subplane of $\spaceP$
of order $\sqrt{x}$.

In addition we have the following lemma:

\begin{lemma}\cite[3.1.2 and 4.1.6]{dembov}\label{lemma:maximalsubplane}
Let $\spaceP$ be a projective plane of order $x$. If $H$ is a group of
automorphisms of $\spaceP$ which does not fix (point-wise) a subplane of $\spaceP$
then the fixed set of $H$ lies inside a fan. If, on the other hand,
$H$ point-wise fixes a subplane of order $y$ then either $y^2=x$ or $y(y+1)\leq x-2$.
\end{lemma}

We are now ready to prove our first result which is very similar to
\cite[Theorem 3]{campraeg2}:

\begin{proposition}\label{proposition:commutingnormal}
Let $G$ be a transitive automorphism group of a projective plane
$\spaceP$ of order greater than $4$. Let $G$ have normal subgroups $M$ and
$N$ such that $M_\alpha\neq 1$ and $N_\alpha\neq1$ for some point
$\alpha$. Then $[N,M]\neq 1$.
\end{proposition}
\begin{proof}
Let $M$ and $N$ be two normal subgroups of $G$ such that there is a
point $\alpha$ so that $M_\alpha\neq 1$ and $N_\alpha\neq1$ and
$[M,N]=1$.

Consider the point $\beta\in\alpha N$ and let $n\in N$ be such that
$\beta =\alpha n$. If $m\in M_\alpha,$ then $\beta m = \alpha nm =
\alpha mn = \beta$. Thus $\alpha N$ is fixed point-wise by
$M_\alpha$. If $\beta\in\alpha N\backslash\{\alpha\}$ and $\linel$ is the
line through $\alpha$ and $\beta$, then $M_\alpha$ fixes $\linel$
set-wise. Thus there is a line $\linel$ through $\alpha$ which is fixed by
$M_\alpha$ and $M_\alpha$ fixes at least two points. A similar result
applies with $N$ replacing $M$.

Next we show that every line through $\alpha$ is fixed either by
$M_\alpha$ or $N_\alpha$. Assume that this is false and let $\linel$ be a
line through $\alpha$ which is fixed by neither. Since $G$ is
line-transitive, there is some point $\beta$ such that $M_\beta$ fixes
$\linel$. Now, since $[M,N]=1$, $N_\alpha$ acts on the set of fixed lines of
$M_\beta$. Thus each image of $\linel$ under the action of $N_\alpha$ is a
line through $\alpha$ fixed by $M_\beta$. Since $N_\alpha$ does not
fix $\linel$, it follows that $M_\beta$ fixes $\alpha$. However, this means
that $M_\beta=M_\alpha$ and hence $M_\alpha$ fixes $\linel$ which is a
contradiction to our assumption.

Thus, for one of $M_\alpha$ and $N_\alpha$, the number of lines through
$\alpha$ which are fixed must be at least $k/2$. Without loss of
generality, this is true
for $N_\alpha$. We now show that the set of fixed points of $N_\alpha$ forms a
subplane of $\spaceP$. By the lemma above it is sufficient to prove
that $N_G(N_\alpha)$ acts transitively on the set of lines fixed by
$N_\alpha$; to show this we demonstrate that $N_\linel=N_\alpha$ for
any line $\linel$ fixed by $N_\alpha$.

Let $\linel$ be any line through $\alpha$ which is fixed by
$N_\alpha$. Let $m\in M$ such that $\linel m\neq \linel$. Then, since
$[M,N]=1$, it follows that $\linel mN_\linel=\linel N_\linel m=\linel
m$, that is $N_\linel$ fixes $\linel m$ and so $N_\linel$ fixes
$\linel m\cap \linel=\{\beta\}$, say. Then $N_\alpha\subseteq
N_\linel\subseteq N_\beta$, and since $N_\alpha$ is conjugate to
$N_\beta$, we obtain $N_\alpha=N_\linel$.

Since $N$ is normal in $G$, $N_G(N_\linel)$ is transitive on the lines
fixed by $N_\linel=N_\alpha$. Thus the fixed set of $N_\alpha$ is a
subplane of $\spaceP$ with line size at least $k/2$. This is a
contradiction of the lemma above.
\end{proof}

\begin{corollary}\label{corollary:allconjugate}
Suppose that $G$ acts transitively on a projective plane
$\spaceP$. Then all components of $G$ are conjugate in $G$.
\end{corollary}
\begin{proof}
If $\spaceP$ is Desarguesian then $G$ contains at most one component
and the statement holds.

By \cite[3.2.15]{dembov} a non-Desarguesian projective plane has order
at least $9$. Thus by the previous theorem any two normal subgroups
$M$ and $N$ of $G$ with $M_\alpha\neq 1$ and $N_\alpha\neq1$ for some point
$\alpha$ satisfy $[N,M]\neq 1$.

Now suppose that $C$ and $D$ are components of $G$ which are not
conjugate in $G$. Then $C^\circ$ and $D^\circ$ are distinct normal subgroups of
$G$. Note that any component contains an
involution and, since the number of points in $\spaceP$ is odd,
each involution must fix a point. The theorem implies that $[C^\circ,
D^\circ]\neq 1$. This is a contradiction.
\end{proof}

We can now prove Theorem \ref{theorem:frameworkprojplane}. Our method
of proof is very similar to that of Camina \cite[Theorem 1]{camina}. First we
state some preliminary results:

\begin{lemma}\cite[Theorem 1]{campraeg2}\label{lemma:faithful}
Let $\spaceP$ be a finite linear space and let $G$ be a
line-transitive automorphism group of $\spaceP$. Let $N$ be a normal
subgroup of $G$. Then $N$ acts faithfully on each of its point orbits.
\end{lemma}

\begin{lemma}\cite[XIII.13.1]{hughespiper}\label{lemma:abelian}
Let $A$ be an abelian automorphism group of a projective plane of
order $x$ then $|A|\leq x^2+x+1.$
\end{lemma}

\begin{frameworktheorem}
Suppose that $G$ acts transitively on a projective plane
$\spaceP$. Then $G$ contains at most one component.
\end{frameworktheorem}
\begin{proof}
By Corollary \ref{corollary:allconjugate}, $\spaceP$ is non-Desarguesian of order $x$ and all
components are conjugate in $G$. Let $C$ be a component of $G$ and
let $C^\circ$ be the normal closure of $C$ in $G$. Write $C^\circ=C_1\circ\dots\circ C_m$
with each $C_i$ isomorphic to $C$ and suppose that $m\geq2$.

Let $D$ be a Sylow $2$-subgroup of $C^\circ$. Since $\spaceP$ has an odd
number of points there is a point $\alpha$ so that $D$ fixes
$\alpha$. Thus $(C_i)_\alpha\neq 1$ for $1\leq i\leq m$. Since $G$
acts transitively on $\spaceP$ this is true for all points
$\alpha$. Choose $\alpha$ so that $(C_1)_\alpha$ has maximal
order. Observe that $[C_2, (C_1)_\alpha]=1$ so $\alpha C_2$ consists
of points fixed by $(C_1)_\alpha$.

Now $C^\circ$ is faithful on all its point orbits by Lemma
\ref{lemma:faithful}. This implies that $\alpha C_2$ contains at
least 5 points as $C_2$ is quasisimple and normal in $C^\circ$. The
fixed set  of $(C_1)_\alpha$ is either a subplane or lies inside a
fan. But, since $C_2$ does not fix any point, we conclude that
$(C_1)_\alpha$ fixes a subplane whose order is at most $\sqrt{x}$.

We now show that for any line $\linel$ incident with $\alpha$ there is
a $j$ so that $(C_j)_\alpha$ fixes $\linel.$ Choose a line $\linel$
incident with $\alpha$. If $(C_1)_\alpha$ fixes $\linel$ there is
nothing to prove. We know that there exists a line, $\linel_1,$ which
is incident with $\alpha$ and is fixed by $(C_1)_\alpha$. But $G$ is
transitive on lines so there is $g\in G$ with $\linel_1g=\linel$. Then
$\beta=\alpha g$ is incident with $\linel$ and $((C_1)_\alpha)^g$
fixes $\linel$. But there exists $j$ so that
$((C_1)_\alpha)^g=(C_j)_\beta$ since $g$ permutes the factors
$C_i$. Let $i\neq j$. Then $(C_i)_\alpha$ commutes with $(C_j)_\beta$
and so acts on the set of lines fixed by $(C_j)_\beta$. If
$(C_i)_\alpha$ fixes $\linel$ then we have proved our claim. If not we
see that $(C_j)_\beta$  fixes at least two lines through $\alpha$ and
so fixes $\alpha$. However $((C_1)_\alpha)^g=(C_j)_\beta$ so by the
maximality of $(C_1)_\alpha$ we have $(C_j)_\alpha=(C_j)_\beta$ and
the claim is proved.

Let $y$ be the order of the subplane fixed by $(C_i)_\alpha.$ Then $m(y+1)\geq
x+1$. If $y=\sqrt{x}$ then this implies that $m\geq \sqrt{x}$. If $y\neq
\sqrt{x}$ then Lemma \ref{lemma:maximalsubplane} implies that $y(y+1)\leq x-2$.
Thus $m\geq y+1$ and so $m\geq \sqrt{x+1}>\sqrt{x}$.

Since $C^\circ$ has an abelian subgroup of order at least $5^m$ it follows
from Lemma \ref{lemma:abelian} that $x^2+x+1\geq 5^m\geq
5^{\sqrt{x}}$. This has no solutions.
\end{proof}

\section{Basic Results and Notation}\label{section:basic}

The notation outlined in this section will hold throughout the rest of the
paper. We also state here a number of basic results which will be used
repeatedly throughout the paper.

\subsection{Projective Plane Results}

Consider a projective plane $\spaceP$ of order $x$ with $v=x^2+x+1$ points and
lines.

\begin{lemma}\cite[p.33]{kantor}\label{lemma:projsyl}
Let $G$ act transitively on a projective plane with $\Galph$ a point-stabilizer.
Then if $p_1$ is a prime $\equiv 2 (3)$ then $\Galph$ contains some Sylow
$p_1$-subgroup of $G$. Moreover, $\Galph$ contains a subgroup of index at most
$3$ in a Sylow $3$-subgroup of $G$.
\end{lemma}

For $g$ an element of $G$ we write $n_g$ for the size of the
$G$-conjugacy class of $g$ in $G$ and $r_g$ for the number of these
conjugates lying in a point-stabilizer $G_\alpha$, for some fixed
point $\alpha$ in $\spaceP$. Furthermore, $d_g$ is the number of fixed
points of $g$. We
will sometimes also write $r_g(B)$ for the number of $G$-conjugates of
$g$ lying in a subgroup $B$ of $G$, so $r_g=r_g(\Galph)$.

We know already that if an automorphism $g$ fixes at least $x$ points
then $g$ is called {\it quasicentral} and $g$ fixes $x+1$, $x+2$ or
$x+\sqrt{x}+1$ points\cite[4.1.7]{dembov}. Furthermore, if an automorphism has $x+1$ or $x+2$ fixed points then it is known as a
{\it perspectivity} and Wagner has proved that if $G$ contains a
nontrivial perspectivity and $G$ acts transitively on $\spaceP$ then
$\spaceP$ is Desarguesian and $G\geq PSL(3,x)$\cite{wagner}.

Now any involution is quasicentral (\cite[3.1.6]{dembov}) and so all
the groups $G$ that we consider contain quasicentral automorphisms. By Wagner's
result we will be interested in the situation when $x$ is a square, say $x=u^2$,
and all quasicentral automorphisms, in particular all involutions, have
$u^2+u+1$ fixed points.

We will be particularly interested in properties of integers of the form $u^2+u+1$ where $u$ is an integer.

\begin{lemma}\label{lemma:easy}
If $x=u^2$ then $x^2+x+1=(u^2+u+1)(u^2-u+1)$, where $(u^2+u+1,u^2-u+1)=1$.
\end{lemma}

\begin{lemma}\cite[p.11]{ljundggren}\label{lemma:primefixed}
If $u^2+u+1=p_1^a$ where $p_1$ is a prime, then either $p_1^a=p_1$ or $p_1^a=7^3$.
\end{lemma}

\begin{lemma}\cite[p.33]{kantor}\label{lemma:kantorinequality}
If $x=u^2$ and $x^2+x+1=p^am$ for a prime $p$ with $a>1$, then either $m>8p^a$ or $p^a=u^2\pm u+1=7^3.$
\end{lemma}


\begin{lemma}\label{lemma:counting}
Let $x=u^2$ and let $g$ be an involution acting on projective plane $\spaceP$
with $u^2+u+1$ fixed points. Then
\begin{itemize}
\item   $\frac{n_g}{r_g}=u^2-u+1$;
\item $d_g=u^2+u+1$;
\item $v=\frac{n_g}{r_g}d_g$ and $(\frac{n_g}{r_g},d_g)=1.$
\end{itemize}
\end{lemma}
\begin{proof}
Count pairs of the form $(\alpha,g)$, where $\alpha$ is a point and $g$ is an involution fixing $\alpha$, in two different ways. Then
$$|\{(\alpha,g):\alpha g=\alpha\}|=vr_g= n_gd_g$$
We know already that $d_g=u^2+u+1$ thus we must have $\frac{n_g}{r_g}=u^2-u+1$ and the result follows.
\end{proof}

\begin{lemma}\label{lemma:largestprime}
Suppose that $g$ is an involution acting on projective plane $\spaceP$ with
$u^2+u+1$ fixed points. If $n_g=2^cp^am$ where $(m,2p)=1$ then the
largest power of $p$ in $v$ is less than or equal to $max(p^a,m+2\sqrt{m}+2)$.
\end{lemma}
\begin{proof}
If $p|\frac{n_g}{r_g}$ then clearly the highest power of $p$ dividing $v$ divides $p^a$. If not, then $u^2-u+1=\frac{n_g}{r_g}$ divides $m$. Then the highest power of $p$ dividing $v$ divides $d_g=u^2+u+1<(u^2-u+1)+2\sqrt{u^2-u+1}+2$.
\end{proof}

It is in our exploitation of the last two results that our treatment
will differ substantially from that of Kantor in the primitive
case. We will make use of the equalities outlined in Lemma
\ref{lemma:counting}, taking $g$ to be a member of a small conjugacy
class of involutions.

\subsection{Group Theory Results and Notation}

We begin with a general lemma which will be useful throughout the chapter.

\begin{lemma}\label{lemma:directproduct}
Let $C<A\times B$. Suppose $|A|<|B:N|$ where $N$ is the largest proper normal subgroup of $B$. Then either:
\begin{itemize}
\item   $C\leq A\times B_1$ for $B_1<B$; or
\item $C=A_1\times B$ for $A_1\leq A$.
\end{itemize}
\end{lemma}
\begin{proof}
Suppose $C\not\leq A\times B_1$ for $B_1<B$. Then define $B_1=\{(1,b):(a,b)\in C\}\cong B$ and observe that the projection $C\to A, (a,b)\mapsto a$ has kernel $K=\{(1,b)\in C\}\lhd B_1$. But $|B_1:K|\leq |A|<|B:N|$ where $N$ is the largest proper normal subgroup of $B$. Thus $K=B_1$ and $C=A_1\times B$ for some $A_1\leq A$ as required.
\end{proof}

Now we want to show that a group $G$ with unique component $L$ cannot
act transitively on a projective plane $\spaceP$ unless it contains a
non-trivial perspectivity.

Recall that $L$ is a component of $G$ provided $L$ is a subnormal quasi-simple subgroup of $G$; a {\it quasi-simple} group $C$ is one such that $C=C'$ ($C$ is equal to its commutator subgroup) and $C/Z(C)$ is simple. We also define an {\it almost simple} group to be a group $G$ such that $N\unlhd G\leq Aut(N)$ where $N$ is a non-abelian simple group; an almost simple group can also be thought of as a group with {\it non-abelian simple socle}, the {\it socle} of a group $G$ being the product of the minimal normal subgroups of $G$. For a fuller discussion see \cite{aschbacher3}.

We write $H.G$ for an extension of a group $H$ by a group $G$ and $H:G$ for a split
extension. An integer $n$ denotes a cyclic group of order $n$, while
$[n]$ (resp. $[q^n]$) denotes an arbitrary soluble group of order $n$ (resp. $q^n$) and $p^n$
denotes an elementary abelian group of order $p^n$ where $p$ is a
prime. We write $|H|_p$ for the highest divisor of $|H|$ which is a
power of a prime $p$.

Put $L_\alpha = G_\alpha\cap L$, the
stabilizer of a point $\alpha$ in the action of $L$ on $\spaceP$. In
general, we will set $M$ to be a maximal subgroup of the component $L$
which contains $\Lalph$. Define $L^\dag:=L/Z(L)$ and $M^\dag:=M/(Z(L)\cap M)$.

Write $G=(L\circ C_G(L)).N$ where $N$ is a subgroup of $Out L$. Then
$G/C_G(L)$ is an almost simple group and we use results about
the maximal subgroups of such groups:

When $L^\dag$ is a classical simple group we use the results of
Aschbacher\cite{aschbacher2} as described in Kleidman and Liebeck
\cite{kl}. These results give information about the maximal subgroups
of a group $L^\dag.N$ with simple socle $L^\dag$ a classical group.

We will sometimes precede the structure of a subgroup of a projective
group with $\hatt$ which means that we are giving the structure of the
pre-image in the corresponding universal group (we call this {\bf hat notation}). For a given element
$g\in L$ we will often write $g^*$ for an element in the corresponding
universal group which projects onto $g$. The symbol $^*$ will also be
used in a different way, with groups, e.g. $P_1^*$, to signal that a
group is a subgroup of a section of $L$ or $L^\dag$. Write $GF(q)$ for the finite field of size $q$.

We now prove a small result which will be very useful:
\begin{lemma}\label{lemma:semilinear}
Suppose that $G$ has a unique component $L$ and $G$ acts transitively on the set of points of a projective plane $\spaceP$. Then, except
when $L=P\Omega ^+(8,q)$, there exists $L\leq H\leq G$ such that $H/C_H(L)\leq \Gamma L$ and $H$ acts transitively on the set of points of $\spaceP$.
Here $\Gamma L$ is the full semilinear classical group associated with $L$.
\end{lemma}
\begin{proof}
The result is trivial except when $L^\dag=PSL(n,q)$ while
$G/C_G(L)$ contains an inverse-transpose automorphism of $L$ and when
$L=Sp(4,2^f)$ while $G/C_G(L)$ contains a graph automorphism of $L$. In both
cases $G$ contains a normal subgroup $H$ of index 2 such that
$H/C_H(L)\leq \Gamma L$. Since we are acting on a set
of odd order, any transitive action of $G$ induces a transitive action
of $H$ as required.
\end{proof}

Lemma \ref{lemma:semilinear} implies that, to prove Theorem A, it is enough to show that the subgroup $H$ cannot act transitively upon a non-Desarguesian projective plane as this implies that the same must hold for $G$. Thus, except when $L^\dag=P\Omega^+(8,q)$, we assume that
$G/C_G(L)\leq \Gamma L$.

We will write $M\in \curlyc{i}$ to mean that $M^\dag$ is in the $i$-th family of
 natural maximal subgroups of $L^\dag$ given by Kleidman and
 Liebeck\cite{kl}. When $M$ is parabolic we will write $M=P_m$ to mean
 that $M$ is a maximal parabolic subgroup fixing a totally singular
 subspace $W$ of dimension $m$ inside the natural classical geometry
 $V$ of dimension $n$.

When $L^\dag$ is an exceptional simple group we use different sources to find information about
maximal subgroups $M$ of $L$. When $M$ is parabolic we refer to
\cite{carter,gorenstein,gl}. In some other cases, the maximal subgroups are
completely enumerated; in particular for $L^\dag={^2B_2(q)}$\cite{suzuki2},
for $L^\dag={^2G_2(q)}$\cite{kleidman2, ward}, for
$L^\dag ={G_2(q)}$\cite{kleidman2, cooperstein}, for
$L^\dag ={^2F_4'(q)}$\cite{malle,atlas} and for $L^\dag ={^3D_4(q)}$\cite{kleidman3}.

In both classical and exceptional cases, we appeal to a result of
Liebeck and Saxl \cite{ls} and
Kantor\cite{kantor} which gives the maximal subgroups of odd index in
an almost simple group. In particular, when the socle is a finite
simple classical group acting on a classical geometry $V$, such a
maximal subgroup either lies in $\curlyc{1}$ (stabilizers of totally
singular or non-singular subspaces) for characteristic $2$ or, when
the characteristic is odd, lies in $\curlyc{1}$, $\curlyc{2}$
(stabilizers of decompositions into subspaces of fixed dimension, $V=\oplus_{i=1}^tV_i$) or $\curlyc{5}$ (stabilizers of subfields) or is
in a small set of listed exceptions.

Finally, when $L^\dag $ is a sporadic simple group we refer to
\cite{aschbacher} which, amongst many other things, lists the maximal
subgroups of odd index.

Our analysis becomes slightly simpler by using the following
result of Camina and Praeger which is a corollary of Lemma \ref{lemma:faithful}:

\begin{lemma}\cite[Corollary 1]{campraeg2}\label{lemma:regularcentre}
Let $N$ be an abelian normal subgroup of a group $G.$ Suppose that $G$
acts line-transitively on a finite linear space $\spaceP$. Then $N$
acts semiregularly on the points of $\spaceP$.
\end{lemma}

In the case where $\spaceP$ is a projective plane we can apply Lemma
\ref{lemma:projsyl}. Thus if $L$ is a unique component of $G$ then $Z(L)$ is
normal in $G$ and must have order only divisible by primes congruent to
$1(3)$ or by $3$ to the first power. In the case where $L$ is a group of
Lie type, for instance, this implies that $L$ is simple unless it is
isomorphic to $E_6(q)$, ${^2E_6(q)}$, $U(n,q)$ or $PSL(n,q)$ for
certain $n$.

\subsection{Hypothesis}\label{section:hypothesis}

Finally we state our hypothesis for the rest of the paper:

\begin{hypothesis}
\begin{enumerate}
\item Suppose that $G$ is a group with a unique component $L$;
\item Suppose that $G$ acts transitively on a set of points of order $v=x^2+x+1$ where $x=u^2, u\in\Z, u\geq 2$;
\item Suppose that all involutions fix $u^2+u+1$ points;
\item Suppose that $\Lalph\leq M$ where $M$ is a maximal subgroup of $L$ of odd
index and that $v>|L:M|$;
\item Except when $L^\dag=P\Omega^+(8,q)$, suppose that $G/C_G(L)\leq \Gamma L$;
\item  Finally suppose that $Z(L)$ has order only divisible by primes congruent to
$1(3)$ or by $3$ to the first power.
\end{enumerate}
\end{hypothesis}

Throughout the rest of the paper we will set $L^\dag$ to be in a particular
family of simple groups and will prove the following result (which, in turn, implies Theorem A):

\begin{result}
If $L\neq PSL(2,q),$ then our hypothesis leads to a contradiction. If $L=PSL(2,q),$ then our hypothesis along with two extra suppositions (described in Section \ref{section:psl23}) leads to a contradiction.
\end{result}

This result is entirely group theoretic and makes no reference to the geometry of projective planes. Note also that Lemmas \ref{lemma:projsyl} to \ref{lemma:largestprime} all apply under our hypothesis since they depend only on the number of points $x^2+x+1$.

\section{$L^\dag $ is alternating or sporadic}

In this section we prove that, if $L^\dag $ is alternating or sporadic, then the hypothesis in Section \ref{section:hypothesis} leads to a contradiction. This implies the following proposition:

\begin{proposition}\label{proposition:alternatingsporadic}
Suppose $G$ has a unique component $L$ such that $L^\dag $ is
isomorphic to an alternating group, $A_n$ with $n\geq 5$, or a
sporadic simple group. Then $G$ does not act transitively on a
projective plane.
\end{proposition}

When $L^\dag $ is a sporadic simple group, the maximal subgroups of $L^\dag $
of odd index are given by Aschbacher\cite{aschbacher}. Aschbacher's
list implies that any maximal subgroup $M$ of odd index in $L$ has
index divisible by $9$ or by a prime congruent to $2(3)$. Since
$\Lalph$ must lie in such a maximal subgroup this contradicts Lemma
\ref{lemma:projsyl}.

Suppose that $L^\dag \cong A_n$, the alternating group on $n$
letters. If $n\neq 6,7$ then $Z(L)\leq 2$ \cite{schur}; thus, by Lemma
\ref{lemma:regularcentre}, $L=L^\dag =A_n$. If $n=6,7$ then $Z(L)\leq 6$
and so, by Lemma \ref{lemma:regularcentre}, $L=A_n$ or $L=3.A_n$.

Assume for the moment that $n>7$ and so $L=A_n$. Let $g\in L=A_n$ be a
double transposition. Then $n_g=\frac{n(n-1)(n-2)(n-3)}{8}$. Now $A_n$
contains an abelian subgroup, $H,$ of size
$2^{\lfloor\frac{n}{2}\rfloor-1}$ which contains at least
$\lfloor\frac{n}{2}\rfloor(\lfloor\frac{n}{2}\rfloor-1)$
$L$-conjugates of $g$.

Since $H$ lies inside a Sylow $2$-subgroup of $L,$ we know that $H$
lies in $\Lalph$ for some point $\alpha$. We conclude that
$$\fracNR\leq \frac{n(n-1)(n-2)(n-3)}{8\lfloor\frac{n}{2}\rfloor(\lfloor\frac{n}{2}\rfloor-1)}.$$
Next we refer to Lemma \ref{lemma:abelian} and observe that $|H|\leq
  v$. Furthermore, for $u>2$, $v<2(\fracNR)^2$. Hence
\begin{eqnarray*}
&&2^{\lfloor\frac{n}{2}\rfloor-1}\leq
  2\frac{n^2(n-1)^2(n-2)^2(n-3)^2}{2^6\lfloor\frac{n}{2}\rfloor^2(\lfloor\frac{n}{2}\rfloor-1)^2}\\
&\implies&2^{\lfloor\frac{n}{2}\rfloor}<n^4\\
&\implies& n\leq43.
\end{eqnarray*}
If $u=2$ then $v=21$ and again we can conclude that $n\leq 43$. Now to
examine the cases where $7<n\leq 43$ we use a method similar
to that in \cite[Section 5]{cnp}.

Consider the usual permutation action of $L=A_n$ as $Alt(\Omega)$,
acting on a set $\Omega$ of size $n$. Then $\Lalph$ contains a Sylow
$p$-subgroup of $L$ for every prime $p\equiv 2(3)$ and a subgroup of
index 3 in a Sylow $3$-subgroup of $L$.

Let $\Gamma$ be the longest orbit of $\Lalph$ in
$\Omega$. If $8\leq n\leq 10$ then, since $\Lalph$ contains a Sylow
$2$-group and a Sylow $5$-group of $L$, $\Lalph^\Gamma$ must be
primitive; if $11\leq n\leq 21$ then the same conclusion comes from
the primes $2$ and $11$; if $22\leq n\leq 33$ then the same conclusion
comes from the primes $2$ and $17$; and if $34\leq n\leq 43$ then the
same conclusion comes from the primes $2$ and $29$. Now $\Lalph^\Gamma$ has odd index in
$Alt(\Gamma)$ and $5$ does not divide the index. By \cite{ls} this
means that $\Lalph^\Gamma$ contains $Alt(\Gamma)$.

For $n\geq 11, n\neq 39$,  we claim that $|\Gamma|\geq n-2$. This is proved
using Lemma \ref{lemma:projsyl} for each individual value of $n$. We
do not reproduce this here but consider, for instance, when $n=16$:
Then $\Lalph$ contains elements with cycle type $(11)$ and $(8,8)$ and
so $|\Gamma|=16\geq n-2$.

Let us examine this case, where $n\geq 11, n\neq 39$. Consider again, $g,$ a double transposition with
$n_g=\frac{n(n-1)(n-2)(n-3)}{8}.$ Then $r_g\geq
\frac{(n-2)(n-3)(n-4)(n-5)}{8}$ and so $\fracNR\leq
\frac{n(n-1)}{(n-4)(n-5)}< 3$ for $n\geq 11$. This is impossible.

For $n=39$ it turns out, using Lemma \ref{lemma:projsyl}, that
$|\Gamma|\geq 34$. Then $\fracNR<3$ and this case is excluded.

For $n=8$ or $10$, the same argument gives $|\Gamma|=n$ and no action
exists. For $n=9$, $|\Gamma|\geq 5$ and, referring to \cite{ls},
$\Lalph$ lies in an intransitive subgroup of $L$ and this contradicts
Lemma \ref{lemma:projsyl}.

Now suppose $n\leq 7$. If $n=5$ or $6$ then Lemma \ref{lemma:projsyl}
implies that $|L:\Lalph|\leq 3$. This is impossible since no subgroup of
such small index exists in $L$. We are left with $n=7$.

When $n=7$ we know that $\Lalph$ contains an element of order
$5$. Examining \cite{atlas} this means that $M^\dag =S_5$ or $A_6$. In
fact we must have $\Lalph=S_5$ or $A_6$. In both cases $\fracNR$ is
not an integer. Thus all cases are excluded.

\begin{remark}
It is worth noting that we could prove Proposition \ref{proposition:alternatingsporadic} directly by appealing to \cite[Theorem
  1]{goho} and then dealing with the cases where $n<21$.
\end{remark}

\section{$L^\dag =PSL(n,q)$, $n>3$}

In this section we assume that $n>3$ and prove that, if $L^\dag =PSL(n,q)$, then the hypothesis in Section \ref{section:hypothesis} leads to a contradiction. This implies the following proposition:

\begin{proposition}\label{proposition:linearbig}
Suppose $G$ has a unique component such that $L^\dag $ is isomorphic to
$PSL(n,q)$ with $n>3$. Then $G$ does not act transitively on a
projective plane.
\end{proposition}

Consider $SL(n,q)$ acting naturally on a vector space $V$. Then recall that a {\it transvection}, $g^*$ say, in $SL(n,q)$ is an automorphism of $V$ such that $g^*-I$ has rank 1 and square 0. We now state the following preliminary result:

\begin{lemma}\label{lemma:nonfusing}
Let $C$ be a conjugacy classes of involutions in $L$ corresponding to either,
\begin{itemize}
\item   diagonalizable involutions in the natural modular
representation of $SL(n,q)$ with $q$ odd; or to
\item   the projective image of transvections in $SL(n,q)$, where
  $q=2^a$ for some integer $a$.
\end{itemize}
Then $C$ is invariant under $\Gamma L$.
\end{lemma}
\begin{proof}
Consider the diagonalizable case first. We need to consider the actions by conjugation of automorphisms of $SL(n,q)$ on a diagonal matrix,
$$g^*=\left(\begin{array}{cccccc} -1 & &&&& \\  & \ddots &&&& \\ &&-1&&& \\ &&&1&& \\ &&&&\ddots& \\ &&&&&1 \end{array}\right).$$
Clearly a field automorphism will preserve $g^*$. Similarly an automorphism lying in $GL(n,q)$ of form,
$$\left(\begin{array}{cccc} 1 &&& \\  & \ddots && \\ &&1& \\ &&&a \end{array}\right)$$
where $a\in GF(q)^*$, also preserves $g^*$. These generate the full outer automorphism group of $SL(n,q)$ in $\Gamma L(n,q)$ and we are done.
In the case where we have a transvection then we consider the actions by conjugation of automorphisms of $SL(n,q)$ on a matrix,
$$g^*=\left(\begin{array}{ccccc} 1&1 & 0&\dots& 0 \\ &1&&\ddots&\vdots \\&& &\ddots & 0 \\ &&&&1 \end{array}\right).$$
Clearly both field automorphisms and the automorphism in $GL(n,q)$ exhibited above preserve $g^*$ and we are done.
\end{proof}

Much of the ensuing treatment will involve counting involutions
$g$. We will take care to ensure that $g$ is always of one of the two
types in this lemma thus ensuring that $n_g=r_g(L)=|L:C_L(g)|$ and
$r_g=r_g(\Lalph)$. Also, observe that we may exclude $PSL(4,2)\cong
A_8$. We begin by restricting the
family within which $M$, a maximal subgroup of $L$ containing
$\Lalph$, may lie:

\subsection{$\Lalph$ must lie in a parabolic subgroup}

By Liebeck and Saxl \cite{ls}, we know that $L_\alpha$ lies inside a maximal subgroup $M$ where
\begin{itemize}
\item       for $q$ odd, $M\in\curlyc{1},\curlyc{2}$ or $\curlyc{5}$; or $n=4$;
\item       for $q$ even, $M\in\curlyc{1}$.
\end{itemize}

\begin{lemma}
$L_\alpha$ cannot lie inside a maximal subgroup from families
$\curlyc{i}, i>1$.
\end{lemma}
\begin{proof}
We may assume that $q$ is odd. In $SL(n,q)$, define
$$g^*=\left(\begin{array}{ccccc} -1 & &&& \\  & -1 &&& \\ &&1&& \\ &&&\ddots& \\ &&&&1 \end{array}\right).$$
Then $g^*$ is centralized in $SL(n,q)$ by $(SL(2,q)\times
SL(n-2,q)).(q-1)$ Then the projective image,$g$ , of $g^*$ is an
involution in $L$ and $n_g$ divides
$$\frac{q^{2(n-2)}(q^{n-1}+\dots+q+1)(q^{n-2}+\dots+q+1)}{q+1}.$$

Examining the order of subgroups $M$ in $\curlyc{2}$ of $\curlyc{5}$
we find that $|M|_p\leq q^{\frac{1}{4}(n-1)n}$ and hence
$|L:M|_p\geq q^{\frac{1}{4}(n-1)n}$. Since $n>3$, we know that $q^2$
divides the index of any maximal subgroup in $\curlyc{2}$ or
$\curlyc{5}$. In the case where $n=4$, the only maximal subgroups of
odd index which do not lie in families $\curlyc{1}$, $\curlyc{2}$ or
$\curlyc{5}$ also have index divisible by $q^2$. Hence $p\geq7$ by
Lemma \ref{lemma:projsyl}. Then, by Lemma \ref{lemma:largestprime},
the largest power of $p$ in $v$ is $q^{2(n-2)}$.

Thus, for $n>4$, $q^{\frac{1}{2}n(n-1)-2(n-2)}=q^{\frac{1}{2}(n^2-5n+8)}$ divides the order of $\Lalph$.\label{TWO} We therefore need to have $\frac{1}{2}(n^2-5n+8)\leq \frac{1}{4}(n-1)n$ and so $n< 7$.

If $n$ is $5$ or $6$ then the only possibility that fits this
inequality is when $M=N_L(L(n,q_0))$ for $q=q_0^2$. But then $|L:M|$
is even and so this case can be excluded. This possibility can also be
excluded when $n=4$. However when $n=4$ we also need to consider the
following further possibilities (note that when $n=4$ we can assume
that $L=PSL(4,q)$):

\begin{itemize}\label{THREE}
\item   $M=\hatt (SL(2,q)\times SL(2,q)).(q-1).2$. (Recall that we use hat notation ($\hatt$) to indicate that we are giving the structure of the pre-image of $M$ in $SL(4,q)$.) In this case
$|L:M|=n_g=\frac{1}{2}q^4(q^2+1)(q^2+q+1)$. Then we know that the
maximum power of $p$ in $v$ is $q^4$ hence $\Lalph$ contains Sylow
$p$-subgroups of $M$. However the index of a parabolic subgroup in
$SL(2,q)$ is even, hence we must have $\hatt (SL(2,q)\times
SL(2,q)).2<\Lalph$. Then we know that for some $\alpha$, $\Lalph>
\hatt \left(\begin{array}{cc} SL(2,q)& \\  & SL(2,q)
\end{array}\right).$ Since $\Lalph$ also contains a Sylow $2$-subgroup
of $PSL(4,q)$, this implies that $\Lalph$ must contain the projective image of $\left(\begin{array}{cccc} 1&&& \\ &-1&& \\ &&1& \\ &&&-1 \end{array}\right)$ which is $L$-conjugate to $g$ and so $r_g\geq q^2(q+1)^2.$ Thus $\fracNR\leq\frac{1}{2}q^2(q^2+1)$ and $v\leq q^4(q^2+1)(q^2+q+1)$ and so $v=\frac{1}{2}q^4(q^2+1)(q^2+q+1)$ contradicting Lemma \ref{lemma:kantorinequality}.

\item   $M=L(4,q_0).[\frac{c}{(q-1,4)}(q_0-1,4)]$ where $c=(q-1)/(q_0-1,\frac{q-1}{(q-1,4)}))$ and $q=q_0^3$. \label{FOUR}Then $|L:M|=(q_0^{12}(q_0^8+q_0^4+1)(q_0^6+q_0^3+1)(q_0^4+q_0^2+1))/(\frac{c}{(q-1,4)}(q_0-1,4))$. Now we know that $p\equiv1(3)$ and so the highest power of $3$ in $c$ is 3. Then we have $9\big||L:M|$ which is impossible.

\item $M$ is of odd index but does not lie in families $\curlyc{1}, \curlyc{2}$ or $\curlyc{5}$. Examining \cite{kl, ls} we find that there are two possibilities: Either $M\in\curlyc{6}$ and $M\cong 2^4.A_6$ or $M\in \curlyc{8}$ and $M\cong PGSp(4,q)$. In the former case, $q^6$ divides $|L:M|$ which is a contradiction. In the latter case, since $p\equiv 1(3)$, we find that $9$ divides $|L:M|$ which, again, is a contradiction.\label{FORTY-NINE}

\end{itemize}
\end{proof}

Thus we assume from here on that $\Lalph$ lies inside
$M\in\curlyc{1}$. This means that $\Lalph$ must always lie inside a parabolic
subgroup, $P_m$, which stabilizes a subspace $W$ of dimension $m$ in
the natural vector space for $G$. We now seek to bound $m$.

\subsection{$\Lalph$ lies in $P_m$, $m$ small}

We begin by noting some preliminary facts which we will use to
establish which parabolic groups $P_m$ are possible candidates to contain
$\Lalph$. In particular we will show that $m$ is small.

\begin{lemma}\label{lemma:choosing}
Suppose $\Lalph$ lies inside $P_m$. For $r|{n\choose m}$, $r$ prime, there exists an integer $a$ such that $(1+q^a+\dots+q^{a(r-1)})$ divides $|L:P_m|$ which, in turn, divides $v.$\label{fifty}
\end{lemma}
\begin{corollary}\label{corollary:divisibles}
Suppose $\Lalph$ lies inside $P_m$.
\begin{itemize}
\item       If $p\equiv 1(3)$ then for all primes $r$ dividing ${n\choose m}$, we must have $r\equiv 1(3)$ or $r=3$ and $9\not|{n\choose m}$.
\item       If $p$ is odd then ${n\choose m}$ must be odd and so either
        \begin{itemize}
        \item       $n$ is odd; or
        \item       $n$ is even and $m$ is even.
        \end{itemize}
\item       If $p=2$ then ${n\choose m}\not\equiv0(4)$.
\end{itemize}
\end{corollary}
\begin{proof}
We need only prove the final statement. Suppose $4|{n\choose m}$. Then
either $(q^2+1)|v$ or $(q+1)^2|v$. This means that either $v$ is
divisible by a prime congruent to $2(3)$ or that $9\big|v$. Both of
these are impossible.
\end{proof}

Note that, since $(n,q)\neq (4,2)$, the smallest index of a parabolic subgroup in $PSL(n,q)$, $n\geq 4$ is 31 (\cite[table 5.2A]{kl}). Since $x$ is a square we know that $v\geq 91$ and so $d_g<2\fracNR$.

\subsubsection{Case: $n$ odd, $p$ odd}\label{section:noddpodd}

In this case $L$ contains the projective image, $g$, of
$$g^*=\left(\begin{array}{cccc} -1 & && \\  & \ddots && \\ &&-1& \\ &&&1 \end{array}\right).$$\label{FIVE, THIRTYEIGHT}
Then $n_g=q^{n-1}(\polysum{n-1})$. Furthermore, since $n\geq 4$, $g$
is conjugate in $G$ to the projective image, $h$, of at least one
other diagonal matrix. Then $g$ and $h$ commute and lie in an
elementary abelian 2-group. Since $\Lalph$ contains a Sylow 2-subgroup
of $L$, we must have $r_g\geq 2$.

Thus $\fracNR\leq\frac{1}{2} q^{n-1}(\polysum{n-1})$, $d_g\leq q^{n-1}(\polysum{n-1})$ and $v\leq \frac{1}{2} q^{2n-2}(\polysum{n-1})^2$. Now observe that,
\begin{eqnarray*}
\frac{1}{2}(\polysum{n-1})^2\geq q^{2n-1} &\implies & (q^n-1)^2\geq 2q^{2n-1}(q-1)^2 \\
&\implies & q^{2n} \geq  2q^{2n-1}(q-1)^2\\
&\implies & q\geq 2(q-1)^2 \\
&\implies & q<3.
\end{eqnarray*}
We know that $q\geq 3$ hence $\frac{1}{2}(\polysum{n-1})^2<
q^{2n-1}$ and $v<q^{4n-3}$. But $|L:P_m|>q^{m(n-m)}$ hence, for $n\geq23$, we have $m\leq4$. We use Corollary \ref{corollary:divisibles} to narrow down the possibilities:

\begin{enumerate}
\item       For $p\equiv 1(3)$ we find, by explicit calculation using
Corollary \ref{corollary:divisibles}, that $m\leq 4$ for all $n$. In fact, checking small $n$ we find that if $m=1,2$ then $n\geq 7$; if $m=3$ then $n\geq 39$; if $m=4$ then $n>70$.\label{SIX}

\item       For $p\not\equiv 1(3)$ then $\frac{n_g}{r_g}|3(\polysum{n-1})$. Hence $d_g<3.q^n$ and so $v<9q^{2n}$. For $n\geq11$ this implies that $m\leq 2$.\label{SEVEN}

Checking the cases where $n<11$ we find that $m\leq 2$ or
$(n,m)=(7,3)$. This final case will be dealt with along with other exceptional
cases at the end of Section \ref{section:exceptional}.
\end{enumerate}

\subsubsection{Case: $n$ even, $p$ odd}\label{section:nevenpodd}

Note that in this case we must have $m$ even and $L$ contains the projective image, $g$, of
$$g^*=\left(\begin{array}{ccccc} -1 &&&& \\  & \ddots &&& \\ &&-1&& \\ &&&1& \\ &&&&1 \end{array}\right).$$
Now $n_g=q^{2(n-2)}(\polysumsq{n-2})(\polysum{n-2}).$ Again $r_g\geq 2$ and so $\fracNR\leq\frac{1}{2}q^{2(n-2)}(\polysumsq{n-2})(\polysum{n-2})$. This gives $d_g\leq q^{2(n-2)}(\polysumsq{n-2})(\polysum{n-2})$ and so $v\leq \frac{1}{2}q^{4(n-2)}(\polysumsq{n-2})^2(\polysum{n-2})^2$.

In a similar fashion to before we know that, for $q\geq 3$ and $n\geq 4$,
$$\frac{1}{2}(\polysumsq{n-2})^2(\polysum{n-2})^2<q^{4n-7}$$
and so $v<q^{8n-15}$. But $|PSL(n,q):P_m|>q^{m(n-m)}$ hence, for $n\geq70$, we have $m\leq8$. Once again we use Corollary \ref{corollary:divisibles} to narrow down the possibilities: \label{TENa, THIRTYEIGHT}

\begin{enumerate}
\item For $p\equiv1(3)$, we find that $n<70$ implies that $m=2$. In fact $(n,m)=(14,2), (38,2) \mbox{ or } (62,2)$.\label{ELEVEN}

\item For $p\not\equiv1(3)$,
$\fracNR|3(\polysumsq{n-2})(\polysum{n-2})< 3q^{2n-3}$. Thus
$v<9q^{4n-5}$. But $|G:P_m|>q^{m(n-m)}$. Thus for $n\geq 18$ we must
have $m\leq4$. For $n<18$, $m\leq 4$ or $(n,m)=(14,6)$. This final
case will be dealt with along with other exceptional cases in Section
\ref{section:exceptional}.\label{TWELVE}
\end{enumerate}

\subsubsection{Case: $p=2$}

In this case $G$ contains the projective image, $g$, of
$$g^*=\left(\begin{array}{ccccc} 1 &0 & \cdots &0 & 1 \\  & 1 &&& 0 \\ && \ddots && \vdots \\ &&&1& 0 \\ &&&&1 \end{array}\right).$$

Here $g^*$ is a transvection and $n_g=(q^{n-1}-1)(\polysum{n-1})$. Examining a Sylow-2 subgroup of $PSL(n,q)$ we see that it contains at least $2(q^{n-1}-1)$ $L$-conjugates of $g$. Since $\Lalph$ must contain one such Sylow 2-subgroup, we conclude that $r_g\geq 2(q^{n-1}-1)$ and so\label{FOURTEEN,FIFTEEN} $\fracNR<\frac{1}{2}(\polysum{n-1})$. Since $d_g<2\fracNR$,  $v<\frac{1}{2}(\polysum{n-1})^2$. Also, since $\Lalph<P_m$ and $|PSL(n,q):P_m|>q^{m(n-m)}$, we conclude that, for $n\geq 10$, $m\leq 2$.

For $n<10$, the fact that $4\not| {n\choose m}$ implies that $(n,m)=(7,3), (8,4)$ or $(9,4)$ if $m>2$. We rule these three possibilities out in turn:

\begin{itemize}
\item   $(9,4)$: This gives $q^{4(9-4)}>q^{2n}$ which is a contradiction.
\item   $(8,4)$: In this case, $(q^4+1)\big| |G:P_4|$ which is impossible.\label{SIXTEEN}
\item   $(7,3)$: In this case, $|G:P_3|=(q^2-q+1)(\polysum{4})(\polysum{6})>\frac{1}{2}(\polysum{6})^2>v$ which is a contradiction.
\end{itemize}

Note that if $m=2$ and $n\equiv 0,1(4)$ then $(q^2+1)\big|v$ which is
impossible. Hence when $m=2$ we assume that $n\equiv 2,3(4)$.

\subsubsection{Cases to be examined}

We now state those values of $m$ for which $\Lalph<P_m$ gives a potential transitive action of $G$:

\begin{enumerate}
\item   $p=2$: $m=1$ ($n\geq5$) or $2$ ($n\geq 6$);
\item $p\not\equiv1(3)$, $p$ odd:
    \begin{itemize}
    \item   $n$ odd: $m=1$ ($n\geq 5)$, $m=2$ ($n\geq 7$) or $(n,m)=(7,3)$;
    \item   $n$ even: $m=2$ ($n\geq 6)$, $m=4$($n\geq 12$) or $(n,m)=(14,6)$;
    \end{itemize}
\item   $p\equiv 1(3)$:
    \begin{itemize}
    \item   $n$ even: $m=2$ ($n =14$ or $n\geq 38$), $m=4, 6, 8$ ($n>70$);
    \item   $n$ odd: $m=1,2$ ($n\geq 7)$, $m=3$ ($n\geq 39$), $m=4$ ($n>70$).
    \end{itemize}
\end{enumerate}

\begin{remark}
Note that $n=4$ is now done. We will assume that $n\geq 5$ from now
on.
\end{remark}

All that remains is to go through the listed cases one at a time
assuming that $\Lalph$ lies inside the given $P_m$ and so $|L:P_m|$
divides $v$. We seek a contradiction. We begin with a preliminary
lemma and corollary which will be useful for counting the number of involutions in $\Lalph$:

\begin{lemma}\label{lemma:conjugatesplit}
Suppose that $q$ is an odd prime power. Assume that the following two
matrices are involutions in $SL(n,q)$, then they are conjugate in $SL(n,q)$:
$$\left(\begin{array}{cc} V & X_1 \\ 0 & W \\ \end{array}\right) \ , \
\left(\begin{array}{cc} V & 0 \\ 0 & W \\ \end{array}\right) $$
where $V\in GL(m,q)$, $W\in GL(n-m,q)$ and $X_1\in M(m\times (n-m),q)$, the set of $m$ by $n-m$ matrices over the field of $q$ elements.
\end{lemma}
\begin{proof}
Since these matrices are involutions we must have
$$VX_1+X_1W=0.$$
Take $X$ such that $2X=-X_1W$. Then $AX=X_1+XW$ and we find that:\label{ONE}
$$\left(\begin{array}{cc} I & X \\ 0 & I \\ \end{array}\right) \left(\begin{array}{cc} V & X_1 \\ 0 & W \\ \end{array}\right) =  \left(\begin{array}{cc} V & 0 \\ 0 & W \\ \end{array}\right)
\left(\begin{array}{cc} I & X \\ 0 & I \\ \end{array}\right). $$
\end{proof}

\begin{corollary}\label{corollary:projectinvolution}
Let $q$ be odd and suppose that $\Lalph$ lies inside a parabolic
subgroup, $P_m$, of $L$ where $P_m=\hatt A:(B:C)$ with $C=q-1$ and
$$
A=\left(\begin{array}{cc} I & M(m\times (n-m),q)  \\  & I
\end{array}\right), \ \
B= \left(\begin{array}{cc} SL(m,q) &  \\  & SL(n-m,q)
\end{array}\right).$$

Define $\pi(\Lalph)$ to be equal to the following set:
$$\left\{\left(\begin{array}{cc} Y_1 &   \\  & Y_2
\end{array}\right)\big|\left(\begin{array}{cc} Y_1 & Z  \\  & Y_2
\end{array}\right)\in A:(B:C), \ {\rm for \ some \ } Z\in M(m\times (n-m),q)\right\},$$
the projection of $P_m$ onto the Levi quotient restricted to $\Lalph$.
Now assume that $\Lalph$ contains an involution $g$
which is the projective image of an involution in $SL(n,q)$,
$g^*=\left(\begin{array}{cc} X_1 & Y \\  & X_2
\end{array}\right)$.

Then $r_g$ is greater than or equal to the number of
$\pi(\Lalph)$-conjugates of the block diagonal matrix
$\left(\begin{array}{cc} X_1 & \\  & X_2 \end{array}\right)$ in
$\pi(\Lalph).$

\end{corollary}

Recall that, in our statement of the corollary, we use hat notation ($\hatt$) to indicate that we are giving the structure of the pre-image of $P_m$ in $SL(n,q)$. Note that in what follows we will assume that $\Lalph$ lies in a
parabolic subgroup which is $L$-conjugate to one of the above form. In
fact, in $PSL(n,q)$ where $n\geq 3$, there are two conjugacy classes
of parabolic subgroups. However, since these two classes are fused by a graph
automorphism, our method extends trivially to cover the other class.

\subsection{Remaining Cases}
\subsubsection{Case: $p=2, m=1$}

Take $g^*$ a transvection as before, with $n_g=(q^{n-1}-1)(\polysum{n-1})$. Recall that $r_g\geq 2(q^{n-1}-1)$ and so $\fracNR\leq \frac{1}{2}(\polysum{n-1})$ and so $v<\frac{1}{2}(\polysum{n-1})^2.$

Then we suppose that $\Lalph=\hatt A.B.C\leq
P_1=\hatt[q^{n-1}]:(SL(n-1,q).(q-1))$. Since $\Lalph$ contains a Sylow
$2$-subgroup of $L$, $A=[q^{n-1}]$ with $B\leq SL(n-1,q)$,
$C\leq(q-1)$. Now $|L:P_1|=\polysum{n-1}$ and thus
$|SL(n-1,q):B|<\frac{1}{2}(\polysum{n-1})$. We know that $B$ contains
a Sylow $2$-subgroup of $SL(n-1,q)$ and so we are in one of the
following situations:

\begin{itemize}
\item       $B\leq P_{m_1}^*$, a parabolic subgroup of
  $SL(n-1,q)$. For $n\geq 5$ and $m_1\geq 2$ observe that
  $|SL(n-1,q):P_{m_1}^*|>q^{2(n-3)}> \frac{1}{2}(\polysum{n-1})$ which
  is impossible. Thus $m_1=1$ and $B< [q^{n-2}]:GL(n-2,q)$. In this
  case $(\polysum{n-1})(\polysum{n-2})\big|v$ and $B=[q^{n-2}]:B_1^*$
  where $|GL(n-2,q):B_1^*|<q$. Thus $B>B_1^*>SL(n-2,q)$.
\item       $B=SL(n-1,q)$.
\end{itemize}

Consider the second situation first. We know that, for some $\alpha$,
$\pi(\Lalph)$ contains $\left(\begin{array}{cc} 1 &  \\  & SL(n-1,q) \end{array}\right).$
We also know that projective images of the following matrices are conjugate in $L$:
$$g^*=\left(\begin{array}{ccccc} 1 &0 & \cdots &0 & 1 \\  & 1 &&& 0 \\
&& \ddots && \vdots \\ &&&1& 0 \\ &&&&1 \end{array}\right), \ \
h^*=\left(\begin{array}{cccccc} 1&&&&& \\ &1 &0 & \cdots &0 & 1 \\  && 1 &&& 0 \\ &&& \ddots && \vdots \\ &&&&1& 0 \\ &&&&&1 \end{array}\right). $$
Thus, by Corollary \ref{corollary:projectinvolution}, $r_g\geq r_g(\hatt SL(n-1,q))\geq(q^{n-2}-1)(\polysum{n-2})$. This implies that $\fracNR<q(q+1)$ and $v\leq q^4+q^2+1$. This is a contradiction for $n\geq 5$.\label{SEVENTEEN}

Thus we assume that we are in the first situation. The same argument though implies that $r_g\geq r_g(\hatt SL(n-2,q))\geq (q^{n-3}-1)(\polysum{n-3})$. \label{EIGHTEEN}This implies that $\fracNR<(q^2+1)^2$ and so  $\fracNR\leq q^4+q^2+1$. This means that $v\leq q^8+4q^6+7q^4+6q^2+3$. We know that $(\polysum{n-1})(\polysum{n-2})|v$ which gives a contradiction for $n\geq 6$.

For $n=5$ we find that $(q^3+q^2+q+1)|v$ hence $(q^2+1)|v$ which implies that a prime $p_1\equiv2(3)$ divides $v$ which is a contradiction.

\subsubsection{Case: $p=2, m=2$}

We assume here that $n\geq 6$ and $\Lalph\leq P_2\cong \hatt
[q^{2(n-2)}]:(SL(2,q)\times SL(n-2,q)).(q-1)$. Now $P_2$ has index $(\polysum{n-1})(\polysum{n-2})/(q+1)$. \label{NINETEEN}We know, as before, that $v<\frac{1}{2}(\polysum{n-1})^2$ hence $|P_2:\Lalph|<q(q+1)$. Now observe that $SL(n-2,q)$ does not have a subgroup of index less than $q(q+1)$ hence $\Lalph> SL(n-2,q)$. As for $m=1$, this implies that $v\leq q^8+4q^6+7q^4+6q^2+3$. This must be greater than the index of $P_2$ and so we must have $n=6$.\label{TWENTY}

In fact when we examine $n=6$ we find that, to satisfy the bound, we must have $q=2$. Explicit calculation of $n_g$, $r_g$ and $|L:P_2|$ excludes this possibility.

\begin{remark}
From here on we assume that $p$ is odd and $n\geq 5$.
\end{remark}

\subsubsection{Case: $p$ odd, $p\not\equiv 1(3)$, $n$ odd, $m$=1}

For the next two cases take $g$ as before for $p$ odd and $n$ odd with $n_g=q^{n-1}(\polysum{n-1})$. We suppose that $\Lalph=\hatt A.B.C<P_1=\hatt[q^{n-1}]:(SL(n-1,q).(q-1))$. Here $A\leq [q^{n-1}]$, $B \leq SL(n-1,q)$ and $C\leq q-1$. Note that $|L:P_1|=\polysum{n-1}$.

Suppose first that $p\neq 3$. Then $\fracNR|\polysum{n-1}$ and
so $v<2(\polysum{n-1})^2$. Then $|P_1:\Lalph|<2(\polysum{n-1})$. Now
$\Lalph$ contains a Sylow-$p$ subgroup of $L$ since $p\equiv 2(3)$. Hence $B$ either lies in a parabolic subgroup, $P_{m_1}^*$, of $SL(n-1,q)$ or $B=SL(n-1,q)$.

Observe that if $m_1$ is odd then $|SL(n-1,q):P_{m_1}^*|$ is
even. Thus we must assume that $m_1$ is even, in which case
$|SL(n-1,q):P_{m_1}^*|>q^{2(n-3)}>2(\polysum{n-1})$ for $n\geq
6$. This is a contradiction. For $n=5$, $P_2^*$ also has even index in
$SL(4,q)$ so can be excluded. Hence we assume that $B=SL(n-1,q)$ and
$|C|$ is even. We know that, for some $\alpha$,
$\pi(\Lalph)$ contains $\left(\begin{array}{cc} \pm1 &  \\  & SL(n-1,q).2
\end{array}\right).$ Thus, appealing to Corollary
\ref{corollary:projectinvolution}, we conclude that $r_g\geq r_g(\hatt
SL(n-1,q).2)\geq q^{n-2}(\polysum{n-2})$ and so $\fracNR<q(q+1)$. This
means that $v\leq q^4+q^2+1$ which is a contradiction for $n\geq
5$.\label{TWENTYONE}

We are left with the case where $p=3$. Now $\Lalph$ contains a group of index $3$ in a Sylow-$3$ subgroup of $L$ and $|L:\Lalph|$ is odd. Hence $B$ either lies in a parabolic subgroup, $P_{m_1}^*$ of $SL(n-1,q)$ or $B=SL(n-1,q)$. The case where $B=SL(n-1,q)$ is ruled out exactly as for $p\neq 3$.

Consider $B\leq P_{m_1}^*<SL(n-1,q)$  and suppose that $n\geq 8$. Then $v>\polysum{7}>1333$ and $\fracNR>31$.
This, combined with the fact that $\fracNR\leq 3(\polysum{n-1})$, means that $v<12(\polysum{n-1})^2$.

Now $B$ lies in $P_{m_1}^*$ and so $m_1$ must be even. Then
$|SL(n-1,q):P_{m_1}^*|>q^{2(n-3)}>12(\polysum{n-1})$ for $n\geq 8$
which is a contradiction. We are left with $n=5$ or $7$. If $n=5$ then
we exclude it as for $p\neq 3$. \label{TWENTYTWO}

For $n=7$, we know that $d_g<2\fracNR\leq 6(\polysum{6})$ and so
$v<18(\polysum{6})^2$. Thus we require that
$q^{2(7-3)}<|SL(n-1,q):P_{m_1}^*|<18(\polysum{6})$. This is
impossible for $q\geq 9$.

When $q=3$ we find that $\fracNR|3(\polysum{6})=3279$.  Now $\fracNR=u^2-u+1$ for some integer $u$ and so $\fracNR\leq \polysum{6}$ and we refer to the case where $p\neq 3$.

\begin{remark}
Note that we have now covered all possible cases where $n=5$ and we assume that $n\geq 6$ from here on.
\end{remark}

\subsubsection{Case: $p$ odd, $p\not\equiv 1(3)$, $n$ odd, $m=2$}

In this case $\Lalph=\hatt A.B.C\leq P_2\cong \hatt
[q^{2(n-2)}]:(SL(2,q)\times SL(n-2,q)).(q-1)$ where $A\leq [q^{n-1}]$,
$B\leq SL(2,q)\times SL(n-2,q)$ and $C\leq q-1$. Now $|L:P_2|=(\polysumsq{n-3})(\polysum{n-1})$.\label{TWENTYTHREE}

Now we know that $\fracNR|3(\polysum{n-1})$. Thus
$v<12(\polysum{n-1})^2$ and hence $|P_2:\Lalph|<12(q+1)^2$. If
$(n,q)\neq (7,3)$ then no subgroup of $SL(n-2,q)$ has index less than
$12(q+1)^2$ unless $(n,q)=(7,3)$. If $(n,q)=(7,3)$ then the only
subgroups of $SL(5,q)$ with indices less than $12(3+1)^2$ are the
parabolic subgroups. These have indices in $SL(5,q)$ divisible by $11$
and so can be excluded. This implies \label{TWENTYFOUR,
  TWENTYFIVE,TWENTYSIX}that in all cases $B=B^*\times SL(n-2,q)$ for
$B^*$ some subgroup of $SL(2,q)$.

Now $B=B^*\times SL(n-2,q)$ implies that $\pi(\Lalph)\geq SL(n-2,q).2$
and so, by Corollary \ref{corollary:projectinvolution}, $r_g>r_g(\hatt SL(n-2,q))>q^{n-3}(\polysum{n-3})$ and $\fracNR<q^2(q^2+1)$ and so $v<q^8+q^4+1$. This gives a contradiction for $n\geq 7$.

\subsubsection{Case: $p$ odd, $p\not\equiv 1(3)$, $n$ even, $m=2$}

For the next two cases, take $g$ as earlier for $p$ odd and $n$ even. Then $n_g=q^{2(n-2)}(\polysum{n-2})(\polysumsq{n-2})$.
As in the previous case, $\Lalph=\hatt A.B.C\leq P_2\cong \hatt
[q^{2(n-2)}]:(SL(2,q)\times SL(n-2,q)).(q-1)$ where $A\leq
[q^{2(n-2)}]$, $B\leq (SL(2,q)\times SL(n-2,q))$, $C\leq q-1$ and
$\pi(\Lalph)=\hatt B.C$. Now $P_2$ has index in $L$,
$(\polysumsq{n-2})(\polysum{n-2})$.

We know, by Lemma \ref{lemma:directproduct}, that one of the following must hold:
\begin{itemize}\label{TWENTYSEVEN}
\item   $B\leq(SL(2,q)\times B_1)$ for some $B_1<SL(n-2,q)$;

\item $B=(B_2\times SL(n-2,q))$ for some $B_2\leq SL(2,q)$.

\end{itemize}

Consider the second possibility. As previously Corollary \ref{corollary:projectinvolution} implies that $r_g\geq r_g(\hatt SL(n-2,q))\geq q^{2(n-4)}(\polysum{n-4})(\polysumsq{n-4})$. Then $\fracNR\leq q^4(q^2+1)^2$ and $v\leq q^{18}$ which is a contradiction for $n>11$. We will need to consider $n=6,8,10$.\label{TWENTYEIGHT, TWENTYNINE, TWO}

We turn to the first possibility above. We know that
$\fracNR|3(\polysum{n-2})(\polysumsq{n-2})$. This implies that
$v<9(\polysum{n-2})^3(\polysumsq{n-2})$ and so
$|P_2:\Lalph|<9(\polysum{n-2})^2$. Thus we must have $B_1$ lying
inside a parabolic subgroup, $P_{m_1}^*$, in $SL(n-2,q)$ with
$|SL(n-2,q):P_{m_1}^*|<9(\polysum{n-2})^2$. We know that $m_1$ must be
even. If $m_1\geq4$ then we know that
$|SL(n-2,q):P_{m_1}^*|>q^{4(n-2-4)}$ which is a contradiction for
$n\geq 12$. Thus $n-2\leq8$ in which case ${m_1}=4$ is not allowed and
so this can also be excluded. Thus we must have ${m_1}=2$. However we
know that ${n\choose2}$ is odd and so $n\equiv 2(4)$, hence
$n-2\equiv0(4)$, hence ${n-2\choose2}$ is even and $|SL(n-2,q):P_2^*|$
is even by Lemma \ref{lemma:choosing}. We may exclude this possibility.\label{THIRTY}

We are left with the possibility that $n=6,8$ or $10$ and
$B=B_2\times SL(n-2,q)$ for some $B_2\leq SL(2,q)$.

Observe first that $A.B.C/A$ acts on the non-identity elements of $A$ by
conjugation. Since $B=B_2\times SL(n-2,q)$, this action has orbits of
size divisible by $q^{n-2}-1$. When $p=3$, $q^{n-2}-1$ does not divide $\frac{q^{2(n-2)}}{3}-1$ hence in all cases we may assume that $A=[q^{2(n-2)}]$.


Then, for some $\alpha$, $A:B$ (or its transpose) has the following form and contains the following conjugate of $g^*$:
$$h^*=\left(\begin{array}{ccccc} I_{2\times 2} &&&& \\  & -I_{2\times 2} &&& \\ &&1&& \\ &&&\ddots & \\ &&&& 1
\end{array}\right)\in
\left(\begin{array}{cc} B_2 & A \\  & SL(n-2,q)   \end{array}\right).$$
\label{THIRTYONE}
Observe that $|A:C_A(h^*)|=q^4$. Thus $r_g\geq q^4 r_g(\hatt
 SL(n-2,q))\geq q^{2n-4}(\polysum{n-4})(\polysumsq{n-4})$. Thus
 $\fracNR\leq (q^2+1)^2$. In fact we may assume that $\fracNR\leq
 q^4+q^2+1$ and so $d_g\leq q^4+3q^2+3$ and $v\leq (q^4+q^2+1)(q^4+3q^2+3).$

Now $|L:P_2|=(\polysumsq{n-2})(\polysum{n-2})>(q^4+q^2+1)(q^4+3q^2+3)$
for $n\geq 6, q\geq 3$. This is a contradiction.

\begin{remark}
Observe that we have now completed the case where $n=6$. We assume that $n\geq 7$ from now on.
\end{remark}

\subsubsection{Case: $p$ odd, $p\not\equiv 1(3)$, $n$ even, $m=4$}

We assume, for this case, that $n\geq 12$. Similarly to the previous
case, $\Lalph=\hatt A.B.C\leq P_4\cong \hatt
[q^{4(n-4)}]:(SL(4,q)\times SL(n-4,q)).(q-1)$ where $A\leq
[q^{4(n-4)}]$, $B\leq (SL(4,q)\times SL(n-4,q))$, $C\leq q-1$ and
$\pi(\Lalph)=\hatt B.C$.

As before, $n_g=q^{2(n-2)}(\polysum{n-2})(\polysumsq{n-2})$ and so $\fracNR|3(\polysum{n-2})(\polysumsq{n-2})$. This implies that $v<9(\polysum{n-2})^3(\polysumsq{n-2})$. Then we have
$$|L:P_4| |P_4:\Lalph|< 9(\polysum{n-2})^3(\polysumsq{n-2})$$
Since $ 9(\polysum{n-2})^3(\polysumsq{n-2})<q^{4n-4}$ we must have $|P_4:\Lalph|<q^{12}$.
We know, by Lemma \ref{lemma:directproduct}, that one of the following must hold:\label{THIRTYFOUR, TWO}

\begin{itemize}
\item   $B\leq (SL(2,q)\times B_1)$ for some $B_1<SL(n-4,q)$. In this case $|SL(n-4,q):B_1|<q^{12}$. For $n\geq 12$ this implies that $B_1$ lies in the parabolic subgroup $P_1^*$ of $SL(n-4,q)$. But this has even index and so can be excluded.

\item $B=(B_2\times SL(n-4,q))$ for some $B_2\leq SL(4,q)$.
\end{itemize}

Thus the second possibility must hold. As before Corollary \ref{corollary:projectinvolution} implies that $r_g\geq r_g(\hatt SL(n-4,q))\geq q^{2(n-6)}(\polysum{n-6})(\polysumsq{n-6})$. Then $\fracNR<q^8(q^4+1)^2$ and
$$d_g<\fracNR+2\sqrt{\fracNR}+2<(q^8+q^4+3)q^4(q^4+1)$$
giving $v\leq q^{12}(q^4+1)^3(q^8+q^4+3)$ which is a contradiction for $n\geq12$.\label{THIRTYFIVE}


\subsubsection{Case: $p$ odd, $p\equiv 1(3)$, $n$ even, $m=2,4,6$ or $8$}

We will take $g$ to be the projective image of,
$$g^*=\left(\begin{array}{ccccc} -1 &&&& \\  & \ddots &&& \\ &&-1&& \\ &&&1& \\ &&&&1 \end{array}\right).$$
Then $n_g=q^{2(n-2)}(\polysumsq{n-2})(\polysum{n-2})$  and we know that $v<q^{8n-15}$. Recall that when $m=2$ we may assume that $n=14$ or $n\geq 38$, otherwise $n>70$.

Let $\Lalph=\hatt A.B.C\leq P_m\cong \hatt [q^{2(n-m)}]:(SL(m,q)\times
SL(n-m,q)).(q-1)$ where $A\leq [q^{m(n-m)}]$, $B\leq(SL(m,q)\times
SL(n-m,q))$, $C\leq q-1$ and $\pi(\Lalph)=\hatt B.C$. Note that
$|L:P_m|>q^{m(n-m)}$ and so $|P_m:\Lalph|<q^{8n-15-mn+m^2}$.

There are two possibilities for $B$, by Lemma \ref{lemma:directproduct}:
\begin{itemize}
\item       $B=(B_2\times SL(n-m,q))$ for some $B_2\leq
SL(m,q)$. Then Corollary \ref{corollary:projectinvolution} implies that $r_g\geq r_g(\hatt SL(n-m,q))\geq q^{2(n-m-2)}(\polysum{n-m-2})(\polysumsq{n-m-2})$. Then $\fracNR\leq q^{2m}(q^m+1)^2$ and $v\leq q^{8m+3}$
Thus we need $m(n-m)<8m+3$ which implies that $m>\frac{n-8}{2}$ which is a contradiction.

\item       $B\leq (SL(m,q)\times B_1)$ for some $B_1<SL(n-m,q)$. By Liebeck and Saxl \cite{ls}, the projective image of $B_1$ in $PSL(n-m,q)$ must lie in families $\curlyc{1}, \curlyc{2}$ or $\curlyc{5}$. The latter two possibilities imply that,
\begin{eqnarray*}
&&\frac{1}{4}n(n-1)<8n-15-mn+m^2 \\
&\implies & n^2-(33-m)n+(60-m^2)<0\\
&\implies & n<33-m\\
&\implies & n=14, m=2.
\end{eqnarray*}
We examine the remaining situation with $n=14, m=2$. Then one subgroup in $\curlyc{2}$ has index less than $q^{8n-15-mn+m^2}=q^{6n-11}$, namely the projective image of $Q_2\cong (SL(6,q)\times SL(6,q)).(q-1).2$ which has even index in $PSL(12,q)$. Similarly the only subgroup in $\curlyc{5}$ with index less than $q^{6n-11}$ is $N_{PSL(12,q)}(PSL(12,q_0))$ where $q=q_0^2$. This also has even index in $PSL(12,q)$ and so can be excluded.

Thus $B_1$ lies in a parabolic subgroup $P_{m_1}^*$ of
$SL(n-m,q)$. Since $n-m$ is even, we must have ${m_1}$ even to have
$i:=|SL(n-m,q):P_{m_1}^*|$ odd. Observe that $q^{m_1(n-m-m_1)}<i<
q^{8n-15-mn+m^2}$. Suppose first that $m+m_1\geq10$. The upper and
lower bounds for $i$ imply that
\begin{eqnarray*}
&&(10-m)(n-10)<8n-15-mn+m^2 \\
&\implies & 2n<m^2-10m+85\\
&\implies & n<35, m=2.
\end{eqnarray*}
We examine the remaining situation with $n<35, m=2$. Referring to
Corollary \ref{corollary:divisibles} the only value of $n$ less than
35 for which $P_2$ has admissible index is $n=14$. But in this case
$m_1=8$ is too large to define a parabolic group in $SL(12,q)$. This
case is excluded. Thus we assume that $m+m_1\leq 8$ and $m\leq 6$. We
split into cases:
\begin{itemize}
\item   Suppose that $m=6$ and so $m_1=2$. Then $|L:P_6|$ odd implies that ${n\choose 6}$ is odd and hence $n\equiv 2(4)$. However this implies that ${{n-6}\choose 2}$ is even and so $i$ is even which is impossible.

\item Suppose that $m=4$ and so $m_1\leq 4$. Recall that, by Corollary \ref{corollary:divisibles}, $5$ does not divide ${n\choose 4}$ hence $n\equiv 4(5)$. However this implies that $5$ divides ${{n-4}\choose m_1}$ which implies, by Lemma \ref{lemma:choosing}, that $i$ is divisible by a prime $p_1\equiv2(3)$ which is impossible.

\item Suppose that $m=2$ and so $m_1\leq 6$. We exclude $m_1=2$ or $6$ in the same way as we excluded $m_1=2$ for $m=6$. We exclude $m_1=4$ in the same way as we excluded $m_1=4$ for $m=4$. Hence we are done.
\end{itemize}

\end{itemize}

\subsubsection{Case: $p$ odd, $p\equiv 1(3)$, $n$ odd, $m=1,2,3$ or $4$}

We will take $g$ to be the projective image of,
$$g^*=\left(\begin{array}{cccc} -1 & && \\  & \ddots && \\ &&-1& \\ &&&1 \end{array}\right).$$
Then $n_g=q^{n-1}(\polysum{n-1})$ and we know that
$v<q^{4n-3}$. Furthermore, by Lemma \ref{lemma:largestprime}, we know
that $|v|_p\leq q^{n-1}$. Recall that, for $m=1$ or $2$, we have $n=7$ or $n\geq 13$, for $m=3$ we have $n\geq 39$ and for $m=4$ we have $n>70$.

Then, in this case, $\Lalph=\hatt A.B.C\leq
P_m=\hatt[q^{n-m}]:(SL(n-m,q).(q-1))$ where $A\leq [q^{n-m}]$, $B\leq
SL(n-m,q)$, $C\leq q-1$ and $\pi(\Lalph)=\hatt B.C$. Note that
$|L:P_m|>q^{m(n-m)}$ and so $|SL(n-m,q):B|<q^{4n-3-mn+m^2}$.

There are two possibilities for $B$, by Lemma \ref{lemma:directproduct}:
\begin{itemize}
\item       $B=(B_2\times SL(n-m,q))$ for some $B_2\leq
SL(m,q)$. We know that $2\leq C$ and so, by Corollary \ref{corollary:projectinvolution}, $r_g\geq r_g(\hatt SL(n-m,q).2)\geq q^{n-m-1}(\polysum{n-m-1})$. Hence $\fracNR<q^m(q^m+1)$ and $v\leq q^{4m}+q^{2m}+1$. Thus we must have
\begin{eqnarray*}
&&m(n-m)<4m+1 \\
&\implies & m^2+(4-n)m+1>0\\
&\implies & m>n-5.
\end{eqnarray*}
This is a contradiction.

\item $B\leq (SL(m,q)\times B_1)$ for some $B_1<SL(n-m,q)$. By Liebeck and Saxl \cite{ls}, the projective image of $B_1$ in $PSL(n-m,q)$ must lie in a subgroup $M$ of $PSL(m,q)$ from  families $\curlyc{1}, \curlyc{2}$ or $\curlyc{5}$. The latter two possibilities imply that,
\begin{eqnarray*}
&&\frac{1}{4}n(n-1)<4n-3-mn+m^2 \\
&\implies & n^2-(17-4m)n+(12-4m^2)<0\\
&\implies & n<17-2m.
\end{eqnarray*}
This implies that either $m=2$ and $n=7$ or $m=1$ and $n=7,13$. In
fact, when $m=1$ and $n=13$ the initial inequality is not satisfied
and this possibility can be excluded. When $m=2$ and $n=7$, the only
possibility is if $B_1\leq M=N_{L_5(q)}(L_5(q_0))$ where $q=q_0^2$. But
$|SL(n-2,q):M|$ is even here and can be excluded. When $m=1$ and $n=7$
we must have $M$ a subgroup of $SL(6,q)$ in $\curlyc{2}$ or
$\curlyc{5}$ and $|SL(6,q):M|<q^{19}$. The only such subgroups are
$M=\hatt (SL(3,q))^2.(q-1).2$ and $M=N_{L(6,q)}(L(6,q_0))$ where
$q=q_0^2$. Both of these subgroups have even index in $SL(6,q)$ and
hence $B_1$ does not lie inside such an $M$.

Thus $B_1$ lies in a parabolic subgroup, $P_{m_1}^*$ of
$SL(n-m,q)$. Write $i:=|SL(n-m,q):P_{m_1}^*|$ and observe that $q^{m_1(n-m-m_1)}<i< q^{4n-3-mn+m^2}$. Suppose first that $m+m_1\geq5$. The upper and lower bounds for $i$ imply that
\begin{eqnarray*}
&&(5-m)(n-5)<4n-3-mn+m^2 \\
&\implies & n<m^2-5m+28.
\end{eqnarray*}
This implies that $n<24$ and either $m=1$ or $m=2$. These cases imply
that $m_1\geq 3$. Now for $i$ to be divisible only by primes
congruent to $1(3)$ or by $3$ but not $9$, we must have
${{n-m}\choose m_1}$ divisible only by primes congruent to $1(3)$ or
by $3$ but not $9$ and hence $n-m\geq 39$ which is a contradiction.

Thus $m+m_1\leq 4$ and $m\leq 3$. Note that if $m$ is odd then $m_1$ must be even since $i$ is odd implies that ${{n-m}\choose m_1}$ is odd. This excludes $m=3$ and ensures that, for $m=1$, $m_1=2$.

Observe some facts about the remaining cases:
\begin{itemize}
\item Suppose that $m=1$ and $m_1=2$. We must have $n\geq 39$ to
ensure that $n$ and ${{n-1}\choose 2}$ are divisible only
by primes congruent to $1(3)$ or by $3$ but not $9$. Then we have
$B_1\leq P_2^*\cong [q^{2(n-3)}]:(SL(2,q)\times SL(n-3,q)).(q-1)$ and,
since $|SL(n-1,q): P_2^*|>q^{2(n-3)}$, then $|P_2^*:B_1|<q^{n+4}$.

\item Suppose that $m=2$. If $n=7$ then $B_1$ lies inside a parabolic
subgroup of $SL(5,q)$. But $5$ divides ${5\choose j}$ for $j=1,2$
which is not allowed. Thus $n\geq 39$ as this is the next smallest
number with allowable divisors of ${n\choose 2}$. Consider
$m_1=2$. Since ${n\choose 2}$ is odd we must have $n\equiv 3(4)$ and
so ${{n-2}\choose 2}$ is even which is a contradiction. Hence
${m_1}=1$ and $B_1\leq P_1^*\cong [q^{n-3}]:SL(n-3,q).(q-1)$. Now
$|SL(n-2,q):P_1^*|\geq q^{n-3}$ and so $|P_1^*:B_1|<q^{n+4}$.
\end{itemize}

Now the only subgroup of $SL(n-3,q)$ in $\curlyc{1}, \curlyc{2}$ or
$\curlyc{5}$ with index less than $q^{n+4}$ is a parabolic subgroup
$P_1^*$ which has even index. Thus, for $m=1$ and $m=2$,
$B_1\geq SL(n-3,q).2$ and so, by Corollary \ref{corollary:projectinvolution},
$r_g\geq r_g(\hatt SL(n-3,q).2)\geq q^{n-4}(\polysum{n-4})$. Hence
$\fracNR<q^3(q^3+1)$ and $v\leq q^{12}+q^6+1$ which is a
contradiction.

\end{itemize}

\subsubsection{Exceptional cases}\label{section:exceptional}

We have deferred two cases in the process of our proof. Firstly we
need to consider the possibility that $n=7, p\not\equiv 1(3)$ is odd
and $\Lalph\leq P_3$, a parabolic subgroup stabilizing a $3$-dimensional
subspace in the vector space for $G$. We exclude this possibility as follows:

Refer to Section \ref{section:noddpodd} when $np$ is odd and suppose that
$\Lalph<P_3$. In this case $\fracNR|3(\polysum{6})$ and
$|L:P_3|=(\polysum{6})(q^6+q^4+q^3+q^2+1)$. Thus $v>q^{12}$ and
$\fracNR>q^5\geq 243$. \label{EIGHT}

Suppose first that $\fracNR<\polysum{6}$. Then $u^2-u+1 =
\fracNR\leq\frac{3}{5}(\polysum{6})$ and $u^2+u+1 =
d_g<q^6+q^4+q^3+q^2+1$ since $\fracNR>243$. Thus $v<|L:P_3|$ which is
a contradiction.\label{NINE}

Then consider the case where $\fracNR\geq\polysum{6}$. We must have
\label{TEN}$v\geq 3(\polysum{6})(q^6+q^4+q^3+q^2+1).$ Suppose that
$\fracNR=\polysum{6}$. Then our lower bound on $v$ implies that
$d_g\geq 3(q^6+q^4+q^3+q^2+1)>2\fracNR$ which is impossible. The only
other possibility is that $\fracNR=3(\polysum{6}) = u^2-u+1$. But then
$u^2+u+1 = d_g<7(q^6+q^4+q^3+q^2+1)$ which again is impossible for
$q\geq 7$. For $q=3,5$ we find that $3(\polysum{6})\neq u^2-u+1$ for
integer $u$ and so these cases can be excluded.

The second possibility that we need to consider is when $n=14,
p\not\equiv 1(3)$ is odd and $\Lalph\leq P_6$, a parabolic subgroup
stabilizing a $6$-dimensional subspace in the vector space for $G$. We
exclude this possibility as follows:

Refer to Section \ref{section:nevenpodd} when $n$ is even and $p$ is
odd and observe that $v<9q^{51}$ and $n_g<q^{49}$. Furthermore
$$\Lalph\leq P_6=\hatt [q^{48}]:(SL(6,q)\times SL(8,q)).(q-1)$$
which has index greater than $q^{48}$. Thus $|P_6:\Lalph|<9q^3$. Now
$SL(6,q)$ and $SL(8,q)$ do not have any subgroups with index this
small, hence $\Lalph>\hatt A.(SL(6,q)\times
SL(8,q))$\label{THIRTEEN} where $A=[q^{48}]\cap \Lalph$. Observe that
$|[q^{48}]:A|\leq 3$. In fact, $A.(SL(6,q)\times SL(8,q))/A$ acts
by conjugation on the non-identity elements of $A$ with orbits of size
divisible by $q^5+\dots+q+1$, hence $A=[q^{48}]$. Then, for some $\alpha$,
$A:(SL(6,q)\times SL(8,q))$ (or its transpose) has the following form
and contains the following conjugate of $g^*$:
$$h^*=\left(\begin{array}{cccc} -1 &&& \\  & I_{5\times 5} && \\ &&-1&
\\ &&&I_{7\times 7}  \end{array}\right)\in
\left(\begin{array}{cc} SL(6,q) & A \\  & SL(8,q)   \end{array}\right).$$
Let $h$ be the projective image of $h^*$. Then $r_g>r_h(\hatt
(SL(6,q)\times SL(8,q)))>q^{10}.q^{14}=q^{24}$. Then $h$ is certainly
centralized by a subgroup of $A$ of size no more than $q^{36}$. Hence
$r_g>q^{36}$. This implies that $\fracNR<q^{13}$ and
$v<q^{27}$ which is a contradiction.

\section{$L=PSL(2,q)$ or $L^\dag =PSL(3,q)$}\label{section:psl23}

In this section we prove firstly that if $L^\dag =PSL(3,q)$ then the hypothesis in Section \ref{section:hypothesis} leads to a contradiction. In the case where $L=PSL(2,q)$ we add two extra suppositions to the hypothesis. For $g\in G$ let $Fix g$ be the set of fixed points of $g$; then our extra suppositions are as follows:
\begin{itemize}
\item Let $g,h\in G$ with $g$ an involution, $h^2=g$. Then $Fix h = Fix g$ or else $|Fix h|=u+1, u+2$ or $u+\sqrt{u}+1$.
\item Let $g,h\in G$ with $g$ an involution, $[g,h]=1$. Then $Fix h = Fix g$ or else $|Fix h\cap Fix g|\leq u+\sqrt{u}+1.$
\end{itemize}
We prove that, with the addition of these suppositions, if $L=PSL(2,q)$, then the hypothesis in Section \ref{section:hypothesis} leads to a contradiction.

To understand the implications of this, suppose for a moment that $G$ is acting on a projective plane of order $x$. Recall that then $g$ fixes a Baer subplane and so $h$, as described in our extra suppositions, either fixes this Baer subplane or else acts as an automorphism of this subplane. Then Lemma \ref{lemma:maximalsubplane} implies that these suppositions must hold. Hence in proving a contradiction we prove the following proposition:

\begin{proposition}\label{proposition:linearsmall}
Suppose that $G$ contains a minimal normal subgroup $L$ isomorphic to
$PSL(2,q)$ with $q\geq 4$ or that $G$ has a unique component $L$ such
that $L^\dag $ is isomorphic to $PSL(3,q)$ with $q\geq 2$. If $G$ acts
transitively on a projective plane $\spaceP$ of order $x$ then
$\spaceP$ is Desarguesian and $G\geq PSL(3,x)$.
\end{proposition}




\subsection{Preliminary facts}

We will need some preliminary facts about $PSL(2,q)$ and
$PSL(3,q)$. As before we assume that $(G/C_G(L))/Z(L)\leq P\Gamma
L(n,q)$ since $|Aut(L): P\Gamma L(n,q)|\leq 2, \ n=2,3$. Observe that
both $PSL(2,q)$ and $PSL(3,q)$ have a single conjugacy class of
involutions of size, in odd characteristic, $\frac{1}{2}q(q\pm1)$ and
$q^2(q^2+q+1)$ respectively and, in even characteristic, $q^2-1$ and
$(q^2-1)(q^2+q+1)$ respectively. Both also have the property that a
Sylow 2-subgroup contains at least 2 such involutions. Since a
point-stabilizer must contain such a Sylow 2-subgroup we conclude that
$r_g\geq 2$. Note also that $PSL(3,q)$ has a single conjugacy class of
transvections and this class does not fuse with any other in $P\Gamma L(3,q)$. \label{lsone}

Liebeck and Saxl\cite{ls} assert that, for $PSL(3,q)$, the maximal subgroups of odd degree lie, as before, in families $\curlyc{1}, \curlyc{2}$ and $\curlyc{5}$ for $q>2$. Note that $PSL(3,2)\cong PSL(2,7)$ and so we will deal with this group in the $PSL(2,q)$ case. We state a result of \cite{moore, wiman} (outlined in \cite{dickson}) which gives the structure of all the subgroups of $PSL(2,q)$:

\begin{theorem}\label{theorem:pslsubgroup}
Let $q$ be a power of the prime $p$. Let $d=(q-1,2)$. Then a subgroup of $PSL(2,q)$ is isomorphic to one of the following groups.
\begin{enumerate}
\item   The dihedral groups of order $2(q\pm 1)/d$ and their subgroups.
\item   A parabolic group $P_1$ of order $q(q-1)/d$ and its subgroups. A Sylow $p$-subgroup $P$ of $P_1$ is elementary abelian, $P\lhd P_1$ and the factor group $P_1/P$ is a cyclic group of order $(q-1)/d$.
\item   $PSL(2,r)$ or $PGL(2,r)$, where $r$ is a power of $p$ such that $r^m =q$.
\item   $A_4, S_4$ or $A_5$.
\end{enumerate}
\end{theorem}

Note that when $p=2$, the above list is complete without the final
entry. Dickson also outlines the conjugacy classes of subgroups of $PSL(2,q)$; in particular it is easy to see that there
are unique $PSL(2,q)$ conjugacy classes of the maximal dihedral
subgroups of size $2(q\pm 1)/d$ as well as a unique $PSL(2,q)$
conjugacy class of parabolic subgroups $P_1$.

The result of Liebeck and Saxl\cite{ls} asserts that all of the families of maximal subgroups can, for some $q$, contain a subgroup of odd index in $PSL(2,q)$ thus, when $L=PSL(2,q)$, we will simply go through the possibilities given in Theorem \ref{theorem:pslsubgroup}.

In the $PSL(3,q)$ case we will also need to know the subgroups of $GL(2,q)$
which can be easily obtained from the subgroups of $PSL(2,q)$.

\begin{theorem}\label{theorem:gl2subgroups}
$H$, a subgroup of $GL(2,q)$, $q=p^a$, is amongst the following up to
conjugacy in $GL(2,q)$. Note that the last two cases may be omitted
when $p=2$.
\begin{enumerate}
\item   $H$ is cyclic;
\item   $H =AD$ where
    $$A\leq\left\{\left(\begin{array}{cc} 1 & 0 \\ \lambda & 1 \\ \end{array}\right):\lambda\in GF(q)\right\}$$
    and $D\leq N(A)$, is a subgroup of the group of diagonal matrices;
\item   $H=\langle c, S\rangle$ where $c|q^2-1$, $S^2$ is a scalar 2-element in $c$;
\item   $H=\langle D,S\rangle$ where $D$ is a subgroup of the group of diagonal matrices, $S$ is an anti-diagonal 2-element and $|H:D|=2$;
\item   $H=\langle SL(2,p^b),V\rangle$ or contains $\langle SL(2,p^b),V\rangle$ as a subgroup of index 2 and here $b|a$, $V$ is a scalar matrix. In the second case, $p^b>3$;
\item   $H/\langle-I\rangle$ is isomorphic to $S_4\times C$, $A_4\times C$, or (with $p\neq 5$) $A_5\times C$, where $C$ is a scalar subgroup of $GL(2,q)/\langle-I\rangle$;
\item   $H/\langle-I\rangle$ contains $A_4\times C$ as a subgroup of index 2 and $A_4$ as a subgroup with cyclic quotient group, $C$ is a scalar subgroup of $GL(2,q)/\langle-I\rangle$.
\end{enumerate}
\end{theorem}
\begin{proof}
In this proof and subsequently, we will refer to subgroups of
$GL(2,q)$ as being {\it of type y}, where $y$ is a number between 1 and 7
corresponding to the list above.

When the characteristic is odd, the proof of this result is given in \cite[Theorem
3.4]{bloom}. When the characteristic is even we know that
$GL(2,q)\cong PSL(2,q)\times (q-1)$. Then, for $H<GL(2,q)$ either
$H\geq SL(2,q)$ and we are in type 5 above, or we have $H\leq
H_1\times (q-1)$ where $H_1$ is maximal in $PSL(2,q)$.

If $H_1=D_{2(q-1)}$ then $H$ is clearly of type 1 or 4. Similarly if $H_1=D_{2(q+1)}$ then $H$ is of type 1 or 3; if $H_1=P_1$ then $H$ is of type 2 in $GL(2,q)$.

Now consider $H\leq PSL(2,q_0)\times (q-1)$. Any maximal subgroup
of $PSL(2,q_0)$ must be an intersection with $D_{2(q\pm1)}$ or $P_1$ (and so is already accounted for) or else
equals $PSL(2,q_1)$ where $q=q_1^b$.

Thus we must consider $H\leq PSL(2,q_1)\times (q-1)$ and $H\not\leq
B\times (q-1)$ for $B<PSL(2,q_1)$. Provided $q_1>2$ this implies
that $H$ is a subgroup of $GL(2,q)$ of type 5. If $q_1=2$ then $PSL(2,q_1)\leq D_{2(q\pm1)}$ and the case is already
accounted for.
\end{proof}

Note that a subgroup of type 1 in $GL(2,q)$ is never maximal in $GL(2,q)$.
 Furthermore type 5 includes $GL(2,q)$ itself. We now proceed with our analysis.

\subsection{$L=PSL(2,q)$}

Assume that $L=PSL(2,q), q\geq 4$. Suppose first that
$G/C_G(L)$ contains $PGL(2,q)$. Then $G$ has a normal subgroup
$N$ of index $2$, $N/C_N(L)$ contains only field automorphisms
and $N$ acts transitively on our set of size $x^2+x+1$. Proving a contradiction for $N$ will give a contradiction for $G$, hence it is enough to assume in general that $G/C_G(L)$ contains only field automorphisms and $|G/C_G(L)|\leq |PSL(2,q)|.\log_pq$.\label{lstwo}

For $q=4,5$ or $9$, $L$ is isomorphic to an alternating group. This
case has already been examined and so these values of $q$
can be excluded. Observe that $P_1$, a parabolic subgroup of
$PSL(2,q)$, has odd index if and only if $p=2$. Furthermore if $p=2$
then $\Lalph\leq P_1$ since $\Lalph$ must contain a Sylow $2$-subgroup of
$PSL(2,q)$. This implies that $n_g=q^2-1$, $r_g=q-1$ and
$u^2-u+1=\fracNR=q+1$. But then $u^2-u=q$ which is impossible. Hence
we assume $\Lalph$ does not lie in a parabolic subgroup of
$PSL(2,q)$ and that $p$ is odd.

Now the only maximal subgroups of $PSL(2,q)$ which contain a Sylow
$p$-subgroup of $PSL(2,q)$ are the parabolic subgroups. Also, for
$q=3^a$ with $a\geq 3$, the only maximal subgroups containing a
subgroup of index $p$ in a Sylow $p$-subgroup of $PSL(2,q)$ are the
parabolic subgroups. Thus Lemma \ref{lemma:projsyl} implies that
$p\equiv 1(3)$ and we assume this from here on.\label{lsthree} Note that,
for an involution $g\in PSL(2,q)$, $n_g= \frac{1}{2}q(q\pm 1)$.

We examine the non-parabolic subgroups of $L$ as candidates to be
$\Lalph$, using Theorem \ref{theorem:pslsubgroup}.

If $\Lalph=A_4$ then $r_g=3$ and, since $r_g\big| n_g$ and $p\equiv
1(3)$, we must have
$n_g=\frac{1}{2}q(q-1)$ and $q\equiv 3(4)$. Similarly if $\Lalph=A_5$
then $r_g=15$ and $q\equiv 3(4)$. But then $\frac{q+1}{4}$
divides $|L:\Lalph|$. Since $\frac{q+1}{4}\equiv 2(3)$ this
contradicts Lemma \ref{lemma:projsyl}.

If $\Lalph=S_4$ then $r_g=9$ and once more $q\equiv 3(4)$. In fact
$\fracNR=\frac{q(q-1)}{18}$. Then in
$PSL(2,q)$ there is a unique conjugacy class of elements of order
4. Let $h$ be such an element and observe that $r_h=6$. Now the fixed
set of $h$ lies inside the fixed set of $g=h^2$ and
$d_h=\frac{1}{3}d_g=\frac{1}{3}(u^2+u+1)$. Referring to our first extra supposition this implies that $|Fix h|=u+1, u+2$ or $u+\sqrt{u}+1$. Since $|Fix h|$ divides $|Fix g|$ we have $\frac{1}{3}(u^2+u+1)= u+\sqrt{u}+1$ and $u=4$. But then $\frac{q(q-1)}{18}=\fracNR=13$ which is impossible.

Now suppose that $\Lalph\leq D_{q\pm 1}$ so $q\pm 1\equiv
0(4)$. Then $\fracNR=\frac{\frac{1}{2}q(q\mp 1)}{\frac{1}{2}|\Lalph|+1}$. Now
$|\fracNR|_p\neq 1$ and so $|\fracNR|_p=|v|_p=q$. Thus $|\Lalph|+2$
divides $q\mp1$.

Define $m:= \frac{q\pm1}{|\Lalph|}$ and assume first that
$m>1$. Observe that $v=q\frac{q\pm1}{\ordLalph}\frac{q\mp1}{2}a$ for some integer
$a$ and $d_g=\frac{|\Lalph|+2}{2}\frac{q\pm1}{\ordLalph}a$. If
$\ordLalph=4$ then $\fracNR=\frac{q(q\mp1)}{6}$ and, in fact, since
$q\equiv 1(3)$, $\fracNR=\frac{q(q-1)}{6}$. But then
$d_g=\frac{3(q+1)}{4}$ and, since $\frac{q+1}{4}\equiv 2(3),$ this is a
contradiction. Thus $\ordLalph>4$.\label{psl2A}

Now observe that $m(|\Lalph|+2)>q\mp1$; furthermore if $(m-1)(|\Lalph|+2)=q\mp1$
then $q\pm1-|\Lalph|+2m-2=q\mp1$. Reducing modulo $4$, this equation
gives $2m\equiv 0(4)$ which is a contradiction since $m\big| v$. Thus
$(m-2)(|\Lalph|+2)\geq q\mp1$. This implies that $m\geq \ordLalph+1$
and so $\ordLalph^2+\ordLalph\leq q\pm 1$.\label{psl2B}

Since $\fracNR<d_g$ we have
\begin{eqnarray*}
&&\frac{q(q\mp1)}{|\Lalph|+2}<\frac{|\Lalph|+2}{2}\frac{q\pm1}{\ordLalph}a \\
&\implies& 2|\Lalph|q(q\mp1)<(\ordLalph^2+4\ordLalph+4)(q\pm1)a \\
&\implies& \ordLalph<\frac{q+1}{q}a.
\end{eqnarray*}
The final inequality follows by using the fact that $\ordLalph>4$
and $\ordLalph^2+\ordLalph\leq q\pm 1$. It then implies that $a>3$.\label{psl2C}

Take $h$ of maximal order in $\Lalph$. Since $\ordLalph>4$ we know
that $h$ is not an involution and $n_h=q(q\mp 1)$ and so
$\frac{n_h}{r_h}=\frac{q(q\mp1)}{2}.$ Thus
$d_h=\frac{q\pm1}{\ordLalph}a$ which means that $d_h<d_g.$ Now $[h,g]=1$ and so, referring to our second extra supposition, $d_h^2<3d_g$ and so $\frac{(q\pm1)^2}{\ordLalph^2}a^2<3\frac{|\Lalph|+2}{2}\frac{q\pm1}{\ordLalph}a$.
This implies that $q\pm 1<\frac{1}{2}\ordLalph^2+\ordLalph$ which is a
contradiction.\label{psl2D}





Hence $m=1$ and $|\Lalph|=q\pm1$. We have two situations. If $q\equiv
3(4)$ then $n_g=\frac{1}{2}q(q-1)$ and $r_g=\frac{1}{2}(q+1)+1$. This
means that $\fracNR$ is a not an integer, which is impossible. If
$q\equiv 1(4)$ then
$\fracNR=\frac{\frac{1}{2}q(q+1)}{\frac{1}{2}(q-1)+1}=q$. Since
$|L:\Lalph|=\frac{1}{2}q(q+1)$ we must have $d_g$ a multiple of
$\frac{q+1}{2}$. The only possibility is that $d_g=\frac{3(q+1)}{2}$
which means that $q=13$ and $v=273$.

In this case $|Fix g|= 21$. But a Sylow $2$-subgroup of $PSL(2,q)$ which centralizes $g$
fixes 9 points; this contradicts our second extra supposition.

Now suppose that $\Lalph=PGL(2,r)$ and $q=r^a$ where $a\equiv
2(4)$. Thus $q\equiv 1(4)$ and
$\fracNR=\frac{\frac{1}{2}q(q+1)}{r^2}$. Now
$\frac{q}{r^2}=|\fracNR|_p\neq|v|_p\geq\frac{q}{r}$ and so
$|\fracNR|_p=1$ and $r=\sqrt{q}$. Then
$u^2-u+1=\fracNR=\frac{1}{2}(q+1)$. Then $u=\frac{c+1}{2}$ where
$c=\sqrt{2q-1}$. This implies that $u^2+u+1=\frac{q+3+2c}{2}$. Now
$|L:\Lalph|=\frac{1}{2}(q+1)\sqrt{q}$ and so $\sqrt{q}$ divides $u^2+u+1$. Now observe that
$\sqrt{q}(\frac{\sqrt{q}+5}{2})>\frac{q+3+2c}{2}$. Furthermore
$\sqrt{q}(\frac{\sqrt{q}-1}{2})<\fracNR$. Thus
$d_g=\sqrt{q}(\frac{\sqrt{q}+e}{2})$ where $e=1$ or $3$.\label{psl2E}

Now $2u=d_g-\fracNR=\frac{e\sqrt{q}-1}{2}$. We also know that
$u=\frac{c+1}{2}$ and so we must have
$e\sqrt{q}-3=2\sqrt{2q-1}$. Since $e=1$ or $3$ we must have
$e=3$. Then
\begin{eqnarray*}
2\sqrt{2q-1}=3\sqrt{q}-3&\implies&2\sqrt{2q}>3\sqrt{q}-3 \\
&\implies& q<(\frac{3}{3-2\sqrt{2}})^2<18^2.
\end{eqnarray*}
This implies that $q=7^2$ or $13^2$. But neither of these satisfy the
equality $2\sqrt{2q-1}=3\sqrt{q}-3$ and so can be excluded.

Now suppose that $\Lalph = PSL(2,r)$ and $q=r^a$ where $a$ is
odd. Then $\fracNR=\frac{\frac{1}{2}q(q\pm 1)}{\frac{1}{2}r(r\pm 1)}$
where $q\mp 1\equiv 0(4)$. Now let $h$ be an element of order
$\frac{r\pm1}{2}$. Then $\frac{n_h}{r_h}=\frac{q(q\mp 1)}{r(r\mp 1)}.$
If $r\equiv 3(4)$ then
$$\fracNR=r^{a-1}(r^{a-1}+r^{a-2}+\dots+r+1)>r^{a-1}(r^{a-1}-r^{a-2}+\dots-r+1)=\frac{n_h}{r_h}.$$
Hence $d_g<d_h$ which is impossible.

Now if $r\equiv 1(4)$ then
$u^2-u+1=\fracNR=r^{a-1}(r^{a-1}-r^{a-2}+\dots-r+1)$ and so
$r^{a-1}-r^{a-2}<u<r^{a-1}$. This means that
\begin{eqnarray*}
r^{2a-2}-r^{2a-3}+\dots-r^a+3r^{a-1}-2r^{a-2}&<&d_g=\fracNR+2u; \\
d_g=\fracNR+2u&<&r^{2a-2}-r^{2a-3}+\dots-r^a+3r^{a-1}.
\end{eqnarray*}
Now $r^{a-1}+r^{a-2}+\dots+r+1$ divides $d_g$. But observe that
\begin{eqnarray*}
&&(r^{a-1}+r^{a-2}+\dots+r+1)(r^{a-1}-2r^{a-2}+2r^{a-3}\dots-2r+3)\\
&<&r^{2a-2}-r^{2a-3}+\dots-r^a+3r^{a-1}-2r^{a-2};
\end{eqnarray*}
\begin{eqnarray*}
&&(r^{a-1}+r^{a-2}+\dots+r+1)(r^{a-1}-2r^{a-2}+2r^{a-3}\dots-2r+4)\\
&>&r^{2a-2}-r^{2a-3}+\dots-r^a+3r^{a-1}.
\end{eqnarray*}

This gives a contradiction and all possibilities are excluded.\label{lssixb}

\subsection{$L^\dag =PSL(3,q)$}

Once again we seek to show that the hypothesis in Section \ref{section:hypothesis} leads to a
contradiction; the usual action of $PSL(3,q)$ on a Desarguesian
projective plane $PG(2,q)$ will not arise due to our restriction that
all involutions fix $u^2+u+1$ points.

Recall that, for $g$ an involution, $n_g=q^2(q^2+q+1)$ for $q$ odd and
$n_g=(q^2-1)(q^2+q+1)$ for $q$ even. We assume here that $q>2$ and we
know that $\Lalph\leq M$  where $M$ is a member of $\curlyc{1},
\curlyc{2}$ or $\curlyc{5}$. We consider the latter two possibilities
first. Observe that, in both cases, $p\equiv 1(3)$ since $p^2$
divides $|PSL(3,q):M|$. \label{lsseven}

Suppose that $M\in\curlyc{2}$. Then $v$ is divisible by $\frac{q^3(q+1)(q^2+q+1)}{6}$. Now the highest power of $q$ in $\fracNR$ is $q^2$. Since $v=\fracNR d_g$ and $(\fracNR, d_g)=1$ we must have $q^3$ dividing $d_g$ and $q^2$ dividing $r_g$. But then $u^2-u+1=\fracNR\leq q^2+q+1$. This means that $v\leq(q^2+q+1)(q^2+3q+3)$ which is a contradiction.\label{lseight}

Suppose that $M=N_{PSL(3,q)}(PSL(3,r))\in\curlyc{5}$ where $q=r^a$ and $a\geq 3$ is an odd integer. Then $|v|_p=\frac{q^3}{r^3}$. Suppose first that $|v|_p=|\fracNR|_p\leq q^2$ and  so $q\leq r^3$. Then we must have $a=3$, $r_g|(q^2+q+1)$ and $r^3$ dividing $\ordLalph$. Since $r_g|(q^2+q+1)$ we cannot have $\Lalph=PSL(3,r)$ or $PSL(3,r).3$. But since $r^3$ divides $\ordLalph$ we must have $\Lalph$ inside a parabolic subgroup $P$ of $PSL(3,r).3$. But observe that then $v$ is divisible by\label{lseightA}
$$|PSL(3,q):P|=\frac{q^3(q^3-1)(q^2-1)}{3r^3(r-1)(r^2-1)}$$\label{lsnine}
which is divisible by $9$, a contradiction. The only other possibility is that $p\not\big|\fracNR$ and $\fracNR\leq q^2+q+1$. But then $q^2\leq r_g\leq r^2(r^2+r+1)$. This is impossible.

Hence we conclude that $M\in\curlyc{1}$. Thus $\Lalph = \hatt A.B$ where $A$ is a subgroup of an elementary abelian unipotent subgroup, $U$, of order $q^2$ and $B$ is a subgroup of odd index in $GL(2,q)$. We will write $B\cap SL(2,q)=(2,q-1).B_1$ where $B_1\leq PSL(2,q)$.

 We will take $\alpha$ to be such that $\Lalph\leq P_1$ where
$$P_1= \hatt \left\{\left(\begin{array}{cc} \frac{1}{det Y} & a \ b \\ 0 & Y
\end{array}\right): Y\in GL_2(q), a,b\in GF(q)\right\}.$$

\subsubsection*{Case: $p\not\equiv 1(3)$}

In this case $|U:A|\leq 3$ and $|P:B_1\cap P|\leq 3$ for some $P\in
Syl_pPSL(2,q)$. If $B_1$ is a subgroup of $P_1^*,$ a parabolic
subgroup of $PSL(2,q)$, then $q+1$ divides the index of $B$ in
$GL(2,q)$ and $p=2$. Then $\Lalph$ is a subgroup of the Borel subgroup
of $PSL(3,q)$ and contains a normal Sylow 2-subgroup $P$. Thus
$r_g=r_g(P)=2q^2-q-1$ and so $r_g\not\big|n_g$ which is a
contradiction.\label{lsten}

If  $B_1=PSL(2,q)$ then $B\geq SL(2,q)$. In fact, in odd
characteristic, $B$ must contain all matrices of determinant $\pm1$
since $|GL(2,q):B|$ is odd. Furthermore in its action by conjugation
on the non-identity elements of $U$, $SL(2,q)$ is transitive. Hence
$A=U$. \label{lstwelve}Thus, in both odd and even characteristic,
$\Lalph$ contains all involutions of the parabolic group: $q^2(q+2)$
of them in the odd case, $(q^2-1)(q+1)$ of them in the even case. In
both cases $r_g\not\big| n_g$ which is a
contradiction.\label{lsthirteen}

For the remaining cases $p|v$ and so $p=3.$ If $B_1\leq
D_{q\pm1}$ then $q|v$ and we must have $q=3$. In this
case $n_g=3^213$ and so $u^2-u+1=\fracNR=3$ or $13$. If $\fracNR=3$
then $v=21$. This contradicts the fact that $|L:M|=13$ and this
divides $v$. So $\fracNR=13, r_g = 9, d_g=21$ and, since $B_1\leq D_{q\pm
1}$ we must have $\Lalph=[3^2]:(8.2)$. But then $\Lalph$ contains more
than 9 involutions and this case is excluded.\label{lsextraone}

If $B_1$ is a proper subgroup of $PSL(2,q)$ isomorphic to $A_4, S_4$
or $A_5$ then $q=3$ or $9$. Now $PSL(2,3)\cong A_4$ and so $q=3$ is
already excluded. If $q=9$ then 5 divides $PSL(2,q)$ and so $B_1\cong
A_5$, but $|PSL(2,9):A_5|$ is even which is impossible.\label{lsthirteenA}

If $B_1\cong PSL(2,r)$ or $B_1\cong PGL(2,r)$ for $q=r^a, a>1$ then
$\frac{q}{r}|v$. Hence $q=9$ and $r=3$. but then $5$ divides
$|PSL(2,9):B_1|$ which is a contradiction.

\subsubsection*{Case: $p\equiv 1(3)$}

In this case $3$ divides
$|PSL(3,q):M|$ and thus we assume that $B$ contains both the Sylow 2
and Sylow 3-subgroups of $GL(2,q)$. In fact $L=PSL(3,q)$ since $Z(L)$
is semiregular (see Lemma \ref{lemma:regularcentre}.) Then $B$ is a subgroup of
$GL(2,q)$ of type 4, 5, 6 or 7 in the list given earlier. Note that
$B$ contains the scalar subgroup of order 3 and so
$|GL(2,q):B|=|\hatt GL(2,q):\hatt B|$.

Observe first that there are two $P_1$-conjugacy classes
of involutions in $P_1$. Only one of these is centralized by a whole
Sylow 2-subgroup, $P$, of $P_1$. Call this conjugacy class
$\mathcal{A}$.

In the case where $\Lalph=\hatt A:B$, that is we have a split extension, we know
that $\hatt B$ contains a Sylow 2-subgroup of $P_1$ and so the involution in
the centre of $\hatt B$ must lie in $\mathcal{A}$. This implies that we can
conjugate by elements of $P_1$ (i.e. choose $\alpha$) such that this
involution $g$ is the projective image of
$$g^*=\left(\begin{array}{ccc} 1 & 0 & 0 \\ 0 & -1 & 0 \\ 0 & 0 & -1 \end{array}\right).$$
We conclude that
$$B\leq \left\{\left(\begin{array}{cc} \frac{1}{det Y} &  \\ &
Y  \end{array}\right):Y\in GL(2,q)\right\}.$$

We begin with two preliminary lemmas:

\begin{lemma}\label{lemma:aqsquared}
Let $p$ be odd and $\Lalph=\hatt A:B\leq P_1$. Suppose that $|A|=q^2$ and that $(|B|,p)=1$. Then $|B|>\frac{|GL(2,q)|}{q^2+q+1}$.
\end{lemma}
\begin{proof}
Let $h$ be an element of order $p$. Then
$$v=\frac{n_h}{r_h} d_h=\frac{(q^2-1)(q^2+q+1)}{q^2-1}d_h=(q^2+q+1)d_h.$$
We have two possibilities:
\begin{enumerate}
\item   Suppose that $h$ is quasi-central. We must have $d_h=u^2+u+1$
where $v=u^4+u^2+1$. Then $u^2-u+1=\frac{n_h}{r_h}=q^2+q+1$ and so
$d_h=q^2+3q+3$. Thus $|B|=\frac{|GL(2,q)|}{q^2+3q+3}a$ for some
integer $a$. If $a=1$ then $|B|$ is not an integer for $q>1$. If
$a\geq 2$ then $|B|>\frac{|GL(2,q)|}{q^2+q+1}$ as required. \label{lsfourteen}

\item Suppose that $h$ is not quasi-central. Then $d_h^2<v$ and so,
$$\left(\frac{v}{q^2+q+1}\right)^2<v \implies v<(q^2+q+1)^2.$$
This implies that $|B|>\frac{|GL(2,q)|}{q^2+q+1}$ as required.
\end{enumerate}
\end{proof}

\begin{lemma}\label{lemma:aq}
Let $p$ be odd and $\Lalph=\hatt A:B\leq P_1$. Suppose that $(|B|,p)=1$. Then $|A|\neq q$.
\end{lemma}
\begin{proof}
Let $h$ be an element of order $p$ and suppose that $|A|=q$. Then
$$v=\frac{n_h}{r_h} d_h=\frac{(q^2-1)(q^2+q+1)}{q-1}d_h=(q+1)(q^2+q+1)d_h.$$
But, since $v$ is odd and $q+1$ is even, this implies that $d_h$ is not an integer. This is a contradiction.
\end{proof}

We now begin our analysis of the different possibilities for
$B$. In the case where $B<GL(2,q)$ is of type 4, 6 or 7 then
Schur-Zassenhaus implies that $A.B$ is a split extension.

Suppose first that $B$ is a subgroup of type $4$ in $GL(2,q)$. Let $\alpha$ be such that $B\leq \langle D,S\rangle$ where $D$ is the subgroup of diagonal matrices and $S$ is an anti-diagonal 2-element. Note that we must have $q$ dividing $|A|$.

Now observe that, since $B$ contains a Sylow 2-subgroup of $D$, we can
choose $\alpha$ such that
\begin{eqnarray*}
\left(\begin{array}{ccc} 1 & e & f \\ 0 & 1 & 0 \\ 0 & 0 & 1 \end{array}\right)\in A
&\implies & \left(\begin{array}{ccc} -1 & e & f \\ 0 & -1 & 0 \\ 0 & 0 & 1 \end{array}\right)^2\in A \\
&\implies & \left(\begin{array}{ccc} 1 & -2e & 0 \\ 0 & 1 & 0 \\ 0 & 0
  & 1 \end{array}\right)\in A \\
&\implies & \left(\begin{array}{ccc} 1 & e & 0 \\ 0 & 1 & 0 \\ 0 & 0 & 1 \end{array}\right)\in A.
\end{eqnarray*}
We conclude that  $A=A_1\times A_2$ where
$$A_1\leq \left\{\left(\begin{array}{ccc} 1 & e & 0 \\ 0 & 1 & 0 \\ 0 & 0 & 1 \end{array}\right):e\in GF(q)\right\} \ , \
A_2\leq \left\{\left(\begin{array}{ccc} 1 & 0 & f \\ 0 & 1 & 0 \\ 0 & 0 & 1 \end{array}\right):f\in GF(q)\right\}.$$
Now consider an element, as given, of $A_1$. Then,
\begin{eqnarray*}
X=\left(\begin{array}{ccc} -1 & 0 & 0 \\ 0 & 0 & a \\ 0 & a^{-1} & 0 \end{array}\right)\in B
&\implies & \left(\begin{array}{ccc} -1 & e & 0 \\ 0 & 0 & a \\ 0 & a^{-1} & 0 \end{array}\right)^2\in A:B \\
&\implies & \left(\begin{array}{ccc} 1 & e & 0 \\ 0 & 1 & 0 \\ 0 & 0 & 1 \end{array}\right)\left(\begin{array}{ccc} 1 & -e & -ae \\ 0 & 1 & 0 \\ 0 & 0 & 1 \end{array}\right)\in A:B \\
&\implies & \left(\begin{array}{ccc} 1 & 0 & ae \\ 0 & 1 & 0 \\ 0 & 0 & 1 \end{array}\right)\in A_2.
\end{eqnarray*}
Thus, for fixed $X$, we have an injection from $A_1$ into $A_2$. There is a similar injection from $A_2$ into $A_1$ and so $|A_1|=|A_2|=\sqrt{|A|}$. Now let
$$E=B\cap \left\{\left(\begin{array}{ccc} -1 & 0 & 0 \\ 0 & 0 & a \\ 0
& a^{-1} & 0 \end{array}\right):a\in GF(q)\right\}$$
and observe that
\begin{eqnarray*}
\left(\begin{array}{ccc} 1 & e & 0 \\ 0 & 1 & 0 \\ 0 & 0 & 1 \end{array}\right)\in A_1,
\left(\begin{array}{ccc} -1 & 0 & 0 \\ 0 & 0 & a \\ 0 & a^{-1} & 0 \end{array}\right)\in E,
&\implies & \left(\begin{array}{ccc} 1 & e & 0 \\ 0 & 0 & a \\ 0 & a^{-1} & 0 \end{array}\right)^2\in A:B \\
&\implies & \left(\begin{array}{ccc} -1 & e & ae \\ 0 & 0 & a \\ 0 & a^{-1} & 0 \end{array}\right)\in A:B
\end{eqnarray*}
and this last element is an involution. We now count all the involutions in $\Lalph$ as follows:

\begin{tabular}[c]{c|c}\label{lsfifteen}
Pre-image of involution $g$ in $SL(3,q)$ & Number of such involutions in $\Lalph$ \\ \\
\hline \\
$\left(\begin{array}{ccc}  1&c&d \\ &-1& \\ && -1 \end{array}\right)$ & $|A|$\\ \\
$\left(\begin{array}{ccc}  -1&0&d \\ &-1& \\ && 1 \end{array}\right)$ & $\sqrt{|A|}$\\ \\
$\left(\begin{array}{ccc}  -1&c&0 \\ &1& \\ && -1 \end{array}\right)$ & $\sqrt{|A|}$\\ \\
$\left(\begin{array}{ccc}  -1&c&d \\ &&a \\ &a^{-1}&  \end{array}\right)$ & $|E|\sqrt{|A|}$\\ \\
\end{tabular}

Thus $r_g=\sqrt{|A|}(\sqrt{|A|}+|E|+2)$ and note that $r_g\leq q(2q+1)$ since $|E|\leq q-1$. Suppose that $(\fracNR,p)=1$. Then $r_g\geq q^2$ and we must have $|A|=q^2$. Alternatively suppose that $(\fracNR,p)\neq 1$. Then\label{lssixteen}
\begin{eqnarray*}
|\fracNR|_p=|v|_p\geq\frac{q^3}{|A|} &\implies & \frac{q^2}{\sqrt{|A|}}\geq |\fracNR|_p\geq\frac{q^3}{|A|} \\
&\implies& |A|\geq q^2.
\end{eqnarray*}
Thus, in either case, $|A|=q^2$. Then, by Lemma \ref{lemma:aqsquared}, $|B|>\frac{|GL(2,q)|}{q^2+q+1}$. \label{lsseventeen}
But $\frac{2(q-1)^2}{7}<\frac{|GL(2,q)|}{q^2+q+1}=\frac{q(q-1)^2(q+1)}{q^2+q+1}$ for $q>1$. Hence $|B|=2(q-1)^2$ and $|E|=q-1$. Then $r_g=q(2q+1)$ which makes $\fracNR$ a non-integer unless $q=1$. This is a contradiction.\label{lseighteen}

Next assume that $B$ is of type $6$ or $7$. To ensure that $B$ has odd index in $GL(2,q)$ we assume that $B\cong 2.(S_4\times C)$ or $B\cong 2.(A_4\times C).2$ where $C\leq Z(GL(2,q))/\langle -I\rangle $. \label{lsnineteen}

Then we must have $q$ dividing $|A|$ since $|v|_p\leq q^2$. We write $|A|=qp^a$ where $a\geq1$ by Lemma \ref{lemma:aq}. Since $\left(\begin{array}{ccc}  1&& \\ &-1& \\ && -1 \end{array}\right)\in B$ this means that $r_g>|A|$.

Suppose first that $q=p^a$ and $|A|=q^2$. By Lemma \ref{lemma:aqsquared},
\begin{eqnarray*}
&&\frac{|GL(2,q)|}{q^2+q+1}<|B|\leq 24(q-1) \\
&\implies& 24(q^2+q+1)>q^3-q \\ \label{lstwenty}
&\implies & q<30.
\end{eqnarray*}
Then $q=7,13$ or $19$. Note that in $GL(2,7)$ subgroups of type $6$ or
$7$ have even index and in $GL(2,19)$ subgroups of type $6$ and $7$
have index divisible by $3$. Hence we are left with $q=13$. In this
case $n_g = 3^2.13.61$ and $v$ is divisible by $|L:M|=3.7.13.61$. Now
since $u^2-u+1=\fracNR$ divides $n_g$ we must have $u=2, 4, 14$
or $23$. But in all of these case $u^2+u+1$ is not divisible by both
$7$ and $61$. Thus $v$ is not divisible by both $7$ and $61$ which is
a contradiction.\label{lstwentyone, lstwentytwo, lsextratwo}

Thus assume now that $q>p^a$ and $|A|<q^2$. Then,
\begin{eqnarray*}
\fracNR<\frac{q^2(q^2+q+1)}{|A|}&\implies & d_g<\frac{q^2(q^2+q+1)}{|A|} + 2\frac{q^2+q+1}{\sqrt{|A|}}+2 \\
&\implies& d_g<\frac{(q^2+2q+1)(q^2+q+1)}{|A|} \\
&\implies& v<\frac{(q+1)^2q^2(q^2+q+1)^2}{|A|^2}.
\end{eqnarray*}
This implies that,
\begin{eqnarray*}
&&\frac{(q^2+q+1)q^3(q-1)^2(q+1)}{|A||B|}\leq v<\frac{q^2(q^2+q+1)^2(q+1)^2}{|A|^2} \\
&\implies & |A|<\frac{(q+1)(q^2+q+1)}{q(q-1)^2}|B|
\end{eqnarray*}
which implies that $|A|<2.|B|$ for $q\geq 7$.\label{lstwentythree}

Now elements from $\hatt 2.C$ do not centralize any element of $\hatt A$. Thus let $m=\frac{(q-1)/2}{|C|}$ and observe that $\frac{q-1}{3m}=|\hatt 2.C|$ divides $|A|-1=qp^a-1$. This in turn means that $\frac{q-1}{3m}$ divides $p^a-1$. Since $q>p^a$ this means that $3m>p$. Then \label{lstwentyfour}
\begin{eqnarray*}
|B|>\frac{|A|}{2}&\implies & 48|C|>\frac{q.p^a}{2} \\
&\implies& 48\frac{q-1}{m} > q.p^a \\
&\implies& p^{a+1}< 144.
\end{eqnarray*}
Since $p\geq 7$, $a\geq 1$ we must have $p=7$, $a=1$. But when $p=7$, $2.(A_4\times C).2$ and $2.(S_4\times C)$ have even index in $GL(2,q)$ which is a contradiction.\label{lstwentyfive}

Thus we are left with the possibility that $B$ is of type $5$ in
 $GL(2,q)$. We want to show that $\Lalph=\hatt A.B$ is a split
 extension and we can choose $\alpha$ such that
$$B\leq  \left\{\left(\begin{array}{ccc} \frac{1}{det Y} &  \\ & Y
\end{array}\right):Y\in B^*\right\}\cong B^*\leq GL(2,q).$$
Observe first that each Sylow 2-subgroup of $\Lalph$ contains a
 unique element of $\mathcal{A}$. thus $\mathcal{A}\cap\Lalph$ is a
 $\Lalph$ conjugacy class. Furthermore there exist at least two
 non-conjugate maximal subgroups, $M_1$, $M_2$, of $B$ which are of
 order not divisible by $p$ and index in $B$ not divisible by $2$. Then, by
 Schur-Zassenhaus, $A:M_1$ and $A:M_2$ are subgroups of $\Lalph$. But
 $M_1,M_2$ must both have centres which are conjugate in $\Lalph$, in
 fact must lie in $\mathcal{A}$. This implies that there exist
 conjugates of $M_1$, $M_2$ which both lie in
$$\left\{\left(\begin{array}{ccc} \frac{1}{det Y} &  \\ & Y
\end{array}\right):Y\in B^*\right\}\cong B^*\leq GL(2,q).$$
These conjugates must generate a complement to $A$ as required.


Now note first that $SL(2,r)\leq GL(2,q)$ implies that
$SL(2,r)\leq SL(2,q)$. Now write $q=r^f$ and observe that, for $f=p_1\dots p_n$
where $p_i$ is prime,
$$SL(2,r)<SL(2,r^{p_1})<\dots<SL(2,r^{p_1\cdots p_{n-1}})<SL(2,q).$$
Since $B$ has odd index in $GL(2,q)$ we assume that all of these primes are odd except, possibly, for $p_1$. What is more, the chain of subgroups given here is maximal except for the first inclusion when $p_1=2$. Now there is a unique conjugacy class in $SL(2,q)$ of maximal subgroups isomorphic to $SL(2,r)$ when $q=r^a$ for $a$ an odd prime. Hence, stepping down the chain of inclusion, we assume that $SL(2,r)$ has a unique conjugacy class in $SL(2,q)$ except when $p_1=2$ in which case there are two conjugacy classes.

By examining \cite[Action Table 3.5G]{kl}) we find that, when $f$ is
even, the two conjugacy classes are fused in $GL(2,r^2)$ through
conjugation by $\left(\begin{array}{cc} \lambda & 0 \\ 0 & 1
\end{array}\right)$ where $\lambda$ generates the group
$GF(r^{2})^*$. Thus, in $GL(2,q)$ there is a unique conjugacy class of
$SL(2,r)$ and we take $\alpha$ such that $B^*$ contains the copy of
$SL(2,r)$ consisting of matrices of determinant 1 with entries in
$GF(r)$.\label{lstwentysix}

Observe that $B^*\ni\left(\begin{array}{cc} 1 & 0 \\ 0 & -1\end{array}\right)$ and so
\begin{eqnarray*}
\left(\begin{array}{ccc} 1 & e & f \\ 0 & 1 & 0 \\ 0 & 0 & 1 \end{array}\right)\in A
&\implies & \left(\begin{array}{ccc} -1 & e & f \\ 0 & -1 & 0 \\ 0 & 0 & 1 \end{array}\right)^2\in A \\
&\implies & \left(\begin{array}{ccc} 1 & e & 0 \\ 0 & 1 & 0 \\ 0 & 0 & 1 \end{array}\right)\in A
\end{eqnarray*}
Once again we conclude that  $A=A_1\times A_2$ where
$$A_1\leq \left\{\left(\begin{array}{ccc} 1 & e & 0 \\ 0 & 1 & 0 \\ 0 & 0 & 1 \end{array}\right):e\in GF(q)\right\} \ , \
A_2\leq \left\{\left(\begin{array}{ccc} 1 & 0 & f \\ 0 & 1 & 0 \\ 0 & 0 & 1 \end{array}\right):f\in GF(q)\right\}.$$
In the same way as earlier we also know that $|A_1|=|A_2|=\sqrt{|A|}$.
We count involutions in $\Lalph$:
\\
\begin{tabular}[c]{c|c}\label{lstwentyseven}
Pre-image of involution $g$ in $SL(3,q)$ & Number of such involutions in $\Lalph$ \\
\hline \\
$\left(\begin{array}{ccc}  1&c&d \\ &-1& \\ && -1 \end{array}\right)$ & $|A|$\\ \\
$\left(\begin{array}{ccc}  -1&c&d \\ &\pm1& \\ && \mp1 \end{array}\right)$ & 2$\sqrt{|A|}$\\ \\
$\left(\begin{array}{ccc}  -1&c&d \\ &\pm1&x \\ && \mp1 \end{array}\right), \ x\neq 0$ & $2(r-1)\sqrt{|A|}$\\ \\
$\left(\begin{array}{ccc}  -1&c&d \\ &v&w \\ &x&-v  \end{array}\right), \ x\neq 0$ & $r(r-1)\sqrt{|A|}$\\ \\
\end{tabular}

Thus $r_g=\sqrt{|A|}(\sqrt{|A|}+r^2+r)$. Now $SL(2,r)$ has orbits of
size $r^2-1$ in its action by conjugation on non-identity elements of
$A$. Hence either $|A|=1$ or $\sqrt{|A|}\geq r$. If $|A|=1$ then,
since $q$ divides $|\Lalph|$, we must have $r=q$ and so
$\fracNR=q^2$. This contradicts Lemma \ref{lemma:primefixed}. Hence
$\sqrt{|A|}\geq r$ and so
$|\fracNR|_p=\frac{q^2}{\sqrt{|A|}r}$. \label{lstwentyeight}

Then either $|\fracNR|_p=1$, $r=q$ and $\sqrt{|A|}=q$ or
$|\fracNR|_p=|v|_p \geq \frac{q^3}{|A|r}p^a$ where
$p^a=\frac{|G|/|L|}{|\Galph|/|\Lalph|}$. In the latter case this means
that
$$\frac{q^2}{\sqrt{|A|}r}\geq \frac{q^3}{|A|r}p^a$$
and so $|A|\geq q^2.p^{2a}.$ This implies that $|A|=q^2$ and $a=0$. In
both cases we find that $|A|=q^2$ and so $r_g=qr(\frac{q}{r}+1+r)$. In
order for this to divide $n_g$ we find that we must have
$r^4+2r^3-r+1$ divisible by $\frac{q}{r}+1+r$. For $q\geq r^6$ this is
clearly a contradiction. Examining cases individually for $q\leq r^5$
we find only contradictions.\label{lstwentynine}

Thus Proposition \ref{proposition:linearsmall} is proved.

\section{$L^\dag =U(n,q)$}

In this section we prove that, if $L^\dag =U(n,q)$, then the hypothesis in Section \ref{section:hypothesis} leads to a contradiction. This implies the following proposition:

\begin{proposition}\label{proposition:unitary}
Suppose $G$ contains a unique component $L$ such that $L^\dag $ is
isomorphic to $U(n,q)$. Then $G$ does not act transitively on a
projective plane.
\end{proposition}

We may assume that $n\geq 3$ and $(n,q)\neq (3,2)$.  We know (\cite[Proposition
2.3.2]{kl}) that our unitary geometry $(V, \kappa)$ has a hyperbolic
basis. Unless stated otherwise, we will write all matrices representing elements of $SU(n,q)$ according to this basis:
\begin{displaymath}
 \left\{
\begin{array}{ll}
\{e_1,f_1,\dots, e_m, f_m\}, &\textrm{if } n=2m; \\
\{e_1,f_1,\dots, e_m,f_m, x\}, &\textrm{if } n=2m+1.
\end{array}\right.
\end{displaymath}
where $\kappa(e_i, e_j)=\kappa(f_i,f_j)=0$, $\kappa(e_i, f_j)=\delta_{ij}$, $\kappa(e_i,x)=\kappa(f_i,x)=0$ for all $i,j$ and $\kappa(x,x)=1$.

We will also need to make use of an orthonormal basis for $(V,\kappa)$. Let
$v_i,w_i$ with $i=1,\dots,m$ be orthonormal vectors such that
$\langle v_i,w_i\rangle=\langle e_i,f_i\rangle$ for all $i=1,\dots,m.$ Our
orthonormal basis $\mathcal{B}$ will consist of these vectors $v_i,
w_i$ with $i=1,\dots,m$, as well as the vector $x$ in the case where
$n$ is odd.

Now the result of Liebeck and Saxl \cite{ls} implies that $L_\alpha$ lies inside a maximal subgroup $M$ where
\begin{itemize}
\item       for $q$ odd, $M\in\curlyc{1}, M\in\curlyc{2}$, $M^\dag =N_{U(n,q)}(U(n,q_0))$ where $q=q_0^a$ and $a$ is odd, or $M^\dag =M_{10}$ and $(n,q)=(3,5)$, or $n=4$;
\item       for $q$ even, $M\in\curlyc{1}$.
\end{itemize}

We show next that, in all cases, $M$ must lie in $\curlyc{1}$:

\begin{lemma}\label{lemma:unitarynonparabolic}
$L_\alpha$ lies inside $M$, where $M$ maximal in $L$ lies inside $\curlyc{1}$.
\end{lemma}
\begin{proof}
We may assume that $p$ is odd. Define $g$ to be the projective image of
$$g^*=\left(\begin{array}{ccccc} -1 & &&& \\  & -1 &&& \\ &&1&& \\ &&&\ddots& \\ &&&&1 \end{array}\right).$$
For $n\neq 4$, $g$ lies in the centre of a maximal subgroup $\hatt
(SU(2,q)\times SU(n-2,q)).(q+1)$. For $n=4$, $g$ lies in the centre of
a maximal subgroup $\hatt (SU(2,q)\times
SU(2,q)).(q+1).2$. Furthermore, $g$ has the same form under
our orthonormal basis $\mathcal{B}$ and, under this basis, $P\Gamma
U(n,q)=U(n,q).\langle\delta, \phi\rangle$ where $\phi$ is a field automorphism and
$\delta$ is conjugation by the projective image of
$$\left(\begin{array}{cccc} a &&& \\  & 1 && \\ &&\ddots& \\ &&&1 \end{array}\right)$$
for some $a\in GF(q^2)^*$, a primitive $(q+1)$-th root of unity. Then
$g$ is centralised by $\langle\sigma, \phi\rangle$ hence $n_g|q^{2(n-2)}b$ where $(q,b)=1$ and $b<q^{2(n-2)}$. Then, by Lemma \ref{lemma:largestprime}, $|v|_p\leq q^{2(n-2)}$.\label{uone}

Suppose that $\Lalph\leq M$ where $M\in \curlyc{2}$, or $M^\dag =N_{U(n,q)}(U(n,q_0))$ where $q=q_0^a$ and $a$ is odd, or $M^\dag =M_{10}$ and $(n,q)=(3,5)$, or $n=4$. Observe that $|U(n,q)|_p=q^{\frac{1}{2}n(n-1)}$ while, for $n\neq 4$, $|M|_p\leq q^{\frac{1}{4}n(n-1)}$. \label{utwo} Thus we must have $\frac{1}{2}n(n-1)-2(n-2)=\frac{1}{2}(n^2-5n+8)\leq\frac{1}{4}n(n-1)$. This implies that $n\leq 6$. We assume this from here on.

Note that we may also assume that $p\equiv 1(3)$ since, in all given
cases, $|U(n,q):M^\dag |$ odd implies that $p^2$ divides $|U(n,q):M^\dag |$. We may immediately rule out the possibility that $M^\dag =M_{10}$.\label{utwoa}

Consider first the case where $n\neq 4$. If $M\in\curlyc{2}$
then $|U(n,q):M^\dag |_p> q^{2(n-2)}$ for $n=3,5$ and $6$ which is a
contradiction.\label{uthreea} If $M=N_{U(n,q)}(U(n,q_0))$ then
$q=q_0^a$ where $a$ is an odd prime. Then $|M|_p\leq
q^{\frac{1}{2a}n(n-1)}$ hence we must have $\frac{1}{2}(n^2-5n+8)\leq
\frac{1}{2a}n(n-1)$ which implies that $n=3$ and
$q=q_0^3$.\label{ufour} Now, when $n=3$, $n_g=q^2(q^2-q+1)$ and
$\Lalph$ contains a Sylow $p$-subgroup of $M$. If $\Lalph\geq U(3,q_0)$
then $r_g=q_0^2(q_0^2-q_0+1)$ but then $r_g\not|n_g$ which is a
contradiction. The only other possibility is that $\Lalph\cap
U(3,q_0)\leq P_1^*$, where $P_1^*$ is a parabolic subgroup of $U(3,q_0)$. But this has even index in $U(3,q_0)$ which is a contradiction.\label{utwenty-one}

Now suppose that $n=4, p\equiv 1(3)$. Note that here $L=U(4,q)$ and
that $n_g=\frac{1}{2}q^4(q^2-q+1)(q^2+1)$. We need to consider the cases
where $M$ is a maximal subgroup of odd index not lying in $\curlyc{1}$. Furthermore we need $|U(4,q):M|_p\leq q^4$. We go through the possibilities in turn.

\begin{itemize}

\item Suppose that $M\in \curlyc{2}$. There exist two subgroups
$M\in\curlyc{2}$ such that $|U(4,q):M|_p\leq q^4$ but only one has odd
  index. We need to rule out this possibility, when $M\cong
  \hatt(SU(2,q)\times SU(2,q)).(q+1).2$ and $|U(4,q):M|_p=q^4$. Then
  $\Lalph$ must contain a Sylow $p$-subgroup of $M$. But the parabolic
  subgroup of $SU(2,q)$ has even index hence we may conclude that,
  \label{uthreeb}for some $\alpha$,
$$\Lalph> \hatt\left(\begin{array}{cc} SU(2,q)& \\  &SU(2,q) \end{array}\right).$$
Then $\Lalph$ contains $h$, the projective image of
$$\left(\begin{array}{cccc}  &1&& \\ 1 & && \\ && &1 \\ &&1& \end{array}\right).$$
Now $h$ is a $U(4,q)$-conjugate of $g$, thus $r_g\geq \frac{1}{2}(q^2-q)^2$. Hence $\fracNR<q^2(q+1)(q+2)$. If $q^4\big| \fracNR$ then we must have $\fracNR=q^4$ which is a contradiction of Lemma \ref{lemma:primefixed}. The only other possibility is that $\fracNR\leq \frac{1}{2}(q^2-q+1)(q^2+1)<\frac{1}{2}q^4$. But then $d_g<q^4$ and so $v<\frac{1}{2}q^4(q^2-q+1)(q^2+1)$ which contradicts $\Lalph\leq M$.\label{utwenty}

\item Suppose that $M\in \curlyc{6}$ or $M\in S$. The only odd index subgroup is $M=2^4.A_6$ where $q\equiv 3(8)$. But then $|U(4,q):M|_p>q^4$ which is a contradiction.\label{utwenty-two}

\item Suppose that $M\in\curlyc{5}$. If $M=N_{U(4,q)}(U(4,q_0))$ then
$q=q_0^a$ where $a$ is an odd prime. Then $|M|_P\leq q^{\frac{6}{a}}$
hence we must have $\frac{1}{2}(n^2-5n+8)=2\leq \frac{6}{a}$ which
implies that $q=q_0^3$.\label{ufourb} However this implies that $9$ divides $|U(n,q):M|$ which is a contradiction.\label{ufive}

The only other odd index subgroup in $\curlyc{5}$ is $M=PGSp(4,q)$
when $q\equiv 1(4)$. \label{utwenty-three}Now, given our original basis $\{e_1,f_1,e_2,f_2\}$ and our original hermitian form $\kappa$, define the form $\kappa_\sharp = \zeta^{-1}\kappa$ over the $GF(q)$-vector space $V_\sharp$ spanned by $\{\zeta e_1, f_1, \zeta e_2,f_2\}$. Here $\zeta$ is an element of $GF(q^2)$ such that $\zeta^q=-\zeta$. Then $\kappa_\sharp$ is a symplectic form over $V_\sharp$.\label{utwenty-fivea}

Clearly if $g^*$ is an isometry for $(\kappa_\sharp, V_\sharp)$ then $g^*$ is an isometry for $(\kappa, V)$ and we have an embedding $Sp(4,q)<SU(4,q)$.\label{utwentyfiveb} This embedding corresponds to a maximal subgroup $PSp(4,q)<U(4,q)$ when $q\not\equiv 1(4)$ and $PGSp(4,q)<U(4,q)$ when $q\equiv 1(4)$. In the latter case, there are two conjugacy classes of $PGSp(4,q)$ in $U(4,q)$; it is this case which concerns us.

Under the orthonormal basis $\{v_1,w_1,v_2,w_2\}$, the two conjugacy classes of $PGSp(4,q)$ in $U(4,q)$ are fused by $x$, the projective image of
$$\left(\begin{array}{cccc} \lambda &&& \\  & 1 && \\ &&1& \\ &&&1 \end{array}\right)$$
where $\lambda\in GF(q^2)$ is a $(q+1)$-primitive element. Thus $r_g$
is the same no matter which of the two conjugacy classes we lie
in. Assume from here on that $\Lalph\leq M=PGSp(4,q)$ preserving
$(\kappa_\sharp, V_\sharp)$.

Then $|U(4,q):M|_p=q^2$, thus $|M:\Lalph|_p\leq q^2$. The only maximal subgroup, $M_1$, of $PSp(4,q)$ such that $|PSp(4,q):M_1|$ is odd and $|PSp(4,q):M_1|_p\leq q^2$ is $(Sp(2,q)\circ Sp(2,q)).2$.\label{utwenty-four} Thus either
\begin{itemize}
\item   $\Lalph = M$ with $v$ divisible by $\frac{1}{2}q^2(q+1)(q^2-q+1)$; or
\item   $\Lalph\cap PSp(4,q)\leq B=(Sp(2,q)\circ Sp(2,q)).2$. Note
that $|(U(4,q):B|_p=q^4$. Since the parabolic subgroups of $Sp(2,q)$
are of even index we must have $\Lalph\cap PSp(4,q)=B$ and so $\Lalph
= B.2$ with $v$ divisible by $\frac{1}{4}q^4(q+1)(q^2-q+1)(q^2+1)$.
\end{itemize}

Under our original basis this implies that, for some $\alpha$,
$$\Lalph>\hatt\left(\begin{array}{cc} SU(2,q)& \\  &SU(2,q) \end{array}\right).$$\label{utwenty-six}
Now $PSp(4,q)$ is normalized in $U(4,q)$ by $h$, the projective image of
$$\left(\begin{array}{cccc}  &1&& \\ 1 & && \\ && &1 \\ &&1&
\end{array}\right).$$
Thus $h$ lies in $\Lalph$ and, as before, we know that $h$ is a
$U(n,q)$-conjugate of $g$. We may conclude that $r_g\geq
\frac{1}{2}(q^2-q)^2$ and so $\fracNR<q^2(q+1)(q+2)$. As in the case
where $M\in\curlyc{2}$ this contradicts $\Lalph=B.2$. We conclude
that $M=PGSp(4,q)$.\label{utwenty-seven}

Now observe that $C_{PSp(4,q)}(h)\cong \hatt GL(2,q).2$ thus
$r_g\geq\frac{1}{2}q^3(q+1)(q^2+1)$ and $\fracNR<q^2$. This implies that
$v<q^2(q+1)(q+2)$ which is a contradiction for $q>4$.\label{utwenty-eight}
\end{itemize}
\end{proof}

Thus $\Lalph$ lies inside a maximal subgroup $M\in\curlyc{1}$. There are two types of $M\in\curlyc{1}$ \cite[Table 3.5B]{kl}:
\begin{itemize}
\item       The parabolic subgroups, $P_m, 1\leq m\leq
\lfloor\frac{n}{2}\rfloor$. Observe that $(q+1)^m$ divides
$|L:P_m|$.\label{useven} This implies that $p=2$. If $q\equiv 1(3)$
then $(q+1)\equiv 2(3)$ and $q+1$ divides $v$. If $m>1$ and
$q\equiv2(3)$ then $9|v$. Neither of these situations are
allowed. Hence $m=1$ and we must have $q=2^a$, $a$ odd.

\item       The subgroups $B_m$ of type $GU(m,q)\perp GU(n-m,q)$ with $1\leq
m<n/2$. In this case $q^{m(n-m)}$ divides $|L:B_m|$ and we must have
$p\equiv 1(3)$. Observe that $q^{m(n-m)}>q^{2(n-2)}$ for
$\frac{n}{2}>m>2$. But we know, by the argument in the previous lemma,
that $|v|_p\leq q^{2(n-2)}$ hence $m\leq 2$.\label{ueight}
\end{itemize}

We now examine these two situations in turn and seek a contradiction.

\subsection{Case: $p=2$, $q=2^a$, $a$ odd, $\Lalph\leq P_1$}

Set $n_e$ to be the even element of $\{n,n-1\}$ while $n_o$ is the odd element. Then $i:=|U(n,q):P_1|=\frac{(q^{n_e}-1)(q^{n_o}+1)}{q^2-1}$.\label{unine} We know that $3|(q+1)\big|i$. In addition, $\polysumsq{n_e-2}|i$ and so for all $r|\frac{n_e}{2}, \polysumsq{2r-2}|i$ which means that for all $r|\frac{n_e}{2}, r\equiv 1(3)$. A similar argument allows us to conclude from the fact that $(q^{n_o-1}-\dots+q^2-q+1)|i$ that for all $r|n_o, r\equiv 1(3)$. We may conclude from this that $n$ is even and $n\equiv 2(12)$. Thus $n\geq 14$.\label{uten}

Now $\Lalph = [q^{2n-3}]:B\leq P_1$ where $B\leq\hatt((q^2-1)\times
SU(n-2,q))$ . We consider the two possibilities given by Lemma \ref{lemma:directproduct}:

\begin{itemize}
\item   $B\leq \hatt((q^2-1)\times B_1)$ for some $B_1<SU(n-2,q)$. We know that $B_1$ must lie in a parabolic subgroup of $SU(n-2,q)$ by Liebeck, Saxl \cite{ls}. However any parabolic subgroup of $SU(n-2,q)$ has index divisible by $q+1$ which would result in $9|v$ which is a contradiction.

\item $B=\hatt(A_1\times SU(n-2,q))$ for some $A_1\leq (q^2-1)$. For some $\alpha$
$$\Lalph\geq \hatt \left(\begin{array}{ccc} SU(n-2,q)&& \\ & 1 &  \\  &&1  \end{array}\right).$$
Now consider transvections in $SU(n,q)$. All transvections are
conjugate to
$$g^*:V\to V, v\mapsto v+s\kappa(v,e_1)e_1$$
for some $s\in GF(q^2), s+s^q=0$\cite[p119]{taylor}. For $W=\langle e_1\rangle$,
define $X_{W,W^\perp}$ to be the subgroup of $SU(n,q)$ consisting of
all transvections of this form. Now suppose that $h\in SU(n,q)$ preserves $W$. Then, for $v\in V$,
\begin{eqnarray*}
v(h^{-1}g^*h) &=& (vh^{-1}+s\kappa(vh^{-1},e_1)e_1)h \\
&=& v+s\kappa(vh^{-1},e_1hh^{-1})e_1h \\
&=& v+s\kappa(v,e_1h)e_1h \\
&=& v+stt^q\kappa(v,e_1)e_1
\end{eqnarray*}
where $t\in GF(q)^*$ is defined via $e_1h=te_1$. Then
$(stt^q)^q+stt^q=tt^q(s+s^q)=0$. Thus $X_{W,W^\perp}$ is normal in the
parabolic subgroup of $SU(n,q)$ stabilizing $W$. Since
$|X_{W,W^\perp}|=q$\cite[p114]{taylor}, we may conclude that, for $g$
the projective image of $g^*$, $\frac{|P_1|}{q-1}$ divides $C_L(g)$. Then, since the only maximal subgroup of $U(n,q)$ whose order is
divisible by $\frac{|P_1|}{q-1}$ is $P_1$, we find that $n_g\leq
\frac{|U(n,q)|(q-1)(n,q+1)2\log_2q}{|P_1|}$. \label{ueleven, utwenty-nine}

Furthermore, $g\in\Lalph$ and, by the same argument, $r_g\geq \frac{|SU(n-2,q)|}{|P_1^*|}$ where $P_1^*$ is a parabolic subgroup of $SU(n-2,q)$. Thus,
$$\fracNR\leq \frac{|U(n,q)|(q-1)(n,q+1)2\log_2q}{|P_1|}\frac{|P_1^*|}{|SU(n-2,q)|}<q^8.$$\label{utwelve}
Then $v<q^{17}$ which is a contradiction.

\end{itemize}

\subsection{Case: $p\equiv1(3), \Lalph \leq B_m, m\leq 2$}

Observe that $|L:B_m|=q^{m(n-m)}\frac{(q^n-(-1)^n)\dots(q^{n-m+1}-(-1)^{n-m+1})}{(q+1)\dots (q^m-(-1)^m)}.$ \label{uthirteen}Consider two situations:
\begin{itemize}
\item       Suppose $n$ is odd. Then $L$ contains the projective image, $g$, of
$$g^*=\left(\begin{array}{cccc} -1 & && \\  & \ddots && \\ &&-1& \\ &&&1 \end{array}\right).$$
Then $g$ is centralized in $U(n,q)$ by $\hatt GU(n-1,q)$. Furthermore, as in Lemma
\ref{lemma:unitarynonparabolic}, $g$ has
the same form, under the basis $\mathcal{B}$, as above and so is centralised
by $\langle\sigma, \phi\rangle$. Hence $n_g|(q^{n-1})(q^{n-1}-\dots-q+1)$. Thus
$|v|_p\leq q^{n-1}$. Suppose that $m\geq 2$, in which case
$|L:B_m|$ is divisible by $q^{2(n-2)}$. Thus we need $2(n-2)\leq
n-1$ which gives $n\leq 3$. For $n=3$ we know that $m=1$. Thus, in
general, $\Lalph\leq B_1=\hatt GU(n-1,q)$. Furthermore $\Lalph$
contains a Sylow $p$-subgroup of $\hatt GU(n-1,q)$.\label{ufourteen}

Thus either $\Lalph\geq\hatt SU(n-1,q)$ or $\Lalph$ lies in a parabolic subgroup of $\hatt GU(n-1,q)$. But $(q+1)$ divides $|\hatt GU(n-1,q):P|$ for $P$ a parabolic subgroup of $\hatt GU(n-1,q)$ which is impossible. Thus $\Lalph\geq\hatt SU(n-1,q)$ and $\Lalph$ contains all the involutions of $\hatt GU(n-1,q)$.

Now, for $n>3$, consider a different involution $g$ as in Lemma
\ref{lemma:unitarynonparabolic}. Then
$n_g=q^{2(n-2)}\frac{(q^n+1)(q^{n-1}-1)}{(q+1)(q^2-1)}$ and $r_g\geq
r_g(\hatt
GU(n-1,q))=q^{2(n-3)}\frac{(q^{n-1}-1)(q^{n-2}+1)}{(q+1)(q^2-1)}$.
This implies that $\fracNR\leq q^4$ and so $\fracNR\leq q^4-q^2+1$ and
$v<q^8+q^4+1$. But $|L:B_1|=q^{n-1}(q^{n-1}-\dots-q+1)$ which is
greater than $q^8+q^4+1$ for $n\geq 7$. For $n=5$,
$2|U(5,q):B_1|>q^8+q^4+1$ and so have $L=U(5,q)$, $\Lalph=B_1$ and
$v=q^4(q^4-q^3+q^2-q+1)$. But, since $q^4>\sqrt{v}$, this implies that
$d_g=q^4$ which contradicts Lemma
\ref{lemma:primefixed}.\label{ufifteen}

For $n=3$ there is a unique conjugacy class of involutions of size $q^2(q^2-q+1)$. Since $\hatt SU(2,q)\leq \Lalph\leq \hatt GU(2,q)$, $\Lalph$ must contain precisely the involutions lying in $\hatt GU(2,q)$ of which there are $q^2-q+1$. Then $\fracNR=q^2$ which contradicts Lemma \ref{lemma:primefixed}.\label{unineteen}

\item   Suppose $n$ is even and let $g$ be as in the proof of Lemma
\ref{lemma:unitarynonparabolic}. Now $|U(n,q):B_1|$ is even and thus
$\Lalph<B_2\cong\hatt (SU(n-2,q)\times SU(2,q)).(q+1)$ and, since
$|v|_p\leq q^{2(n-2)}$, $\Lalph$ contains a Sylow $p$-subgroup of
$\hatt (SU(n-2,q)\times SU(2,q))$. Note that, since $B_2$ is
non-maximal in $L=U(4,q)$, we may assume that $n\geq
6$. \label{useventeen}

Now the index of the parabolic subgroups of $SU(n-2,q)$ in $SU(n-2,q)$
is even. Hence we must have $\Lalph> \hatt SU(n-2,q)$. For some $\alpha$, we may assume that
$$\Lalph\geq \hatt\left(\begin{array}{ccc} SU(n-2,q)&& \\ & 1 &  \\  &&1  \end{array}\right).$$
Now $g$ is centralized in $L$ by some conjugate of $B_2$. This implies
that
$$n_g=q^{2(n-2)}\frac{(q^n-1)(q^{n-1}+1)}{(q+1)(q^2-1)} \ \ {\rm and}
\ \ r_g\geq q^{2(n-4)}\frac{(q^{n-2}-1)(q^{n-3}+1)}{(q+1)(q^2-1)}.$$
Thus $\fracNR\leq q^6(q^2+1)$ and $v\leq q^{16}+q^{15}$ and, for $n\geq 8$, this contradicts $\Lalph\leq B_2$.\label{ueighteen}

We are left with the possibility that $n=6$. But $2|U(6,q):B_2|>q^{16}+q^{15}$, thus $\Lalph=B_2$ and $v=q^8(q^4+q^2+1)(q^4-q^3+q^2-q+1)$. But then $q^8\geq \sqrt{v}$ and so $d_g=q^8$ which contradicts Lemma \ref{lemma:primefixed}.
\end{itemize}

Thus Proposition \ref{proposition:unitary} is proven.

\section{$L=PSp(n,q)$}

In this section we prove that, if $L=PSp(n,q)$, then the hypothesis in Section \ref{section:hypothesis} leads to a contradiction. This implies the following proposition:

\begin{proposition}\label{proposition:symplectic}
Suppose $G$ contains a minimal normal subgroup isomorphic to
$PSp(n,q)$ with $n\geq 4$. Then $G$ does not act transitively on a
projective plane.
\end{proposition}

We know \cite[Proposition 2.4.1]{kl} that our symplectic geometry
$(V,\kappa)$ has a symplectic basis. Unless stated otherwise, we will
write all matrix representations of $Sp(n,q)$ according to this basis,
$\{e_1, f_1,\dots, e_m,f_m\}$, where $n=2m$. Here $\kappa (e_i,
e_j)=\kappa (f_i, f_j)=0$ and $\kappa (e_i, f_j)=\delta_{ij}$.

By Liebeck and Saxl \cite{ls}, we know that $L_\alpha$ lies inside a maximal subgroup $M$ where
\begin{itemize}
\item       for $q$ odd, $M\in\curlyc{1},\curlyc{2}$ or
$M=N_{PSp(n,q)}(PSp(n,q_0))$ or $n=4$;
\item       for $q$ even, $M\in\curlyc{1}$.
\end{itemize}
Note that when $n=4$ we can assume that $q>3$ since $PSp(4,3)\cong
U(4,2)$ which has already been covered.

\begin{lemma}
$\Lalph$ lies inside a maximal subgroup from family $\curlyc{1}$.
\end{lemma}
\begin{proof}
Assume that $q$ is odd and that $\Lalph\leq M$ where $M$ is
a maximal subgroup of $PSp(n,q)$ that does not lie in
$\curlyc{1}$. Observe that in $PSp(n,q)$ there exists a subgroup
$B\cong Sp(2,q)\circ Sp(n-2,q)$.

For $n\neq 4$, by
\cite[Lemma 3.2.1 and Table 3.5.c]{kl}, $B$ is normal in a $P\Gamma
Sp(n,q)$-maximal subgroup $B_\Gamma$ such that $|P\Gamma
Sp(n,q):B_\Gamma|=|L:B|$. Thus, for $n\neq 4$, the involution $g\in
Z(B)$ has $n_g=|L:B|=q^{n-2}(\polysumsq{n-2})$.

When $n=4$ the same argument applies to $B\cong (Sp(2,q)\circ
Sp(2,q)).2$ and the involution $g\in Z(B)$ has
$n_g=\frac{1}{2}q^2(q^2+1)$.

Therefore the highest power of $p$ in $v$ is at most $q^{n-2}$. The
lowest index of $p$ among maximal subgroups $M\in \curlyc{2}$ or
$M=N_{PSp(n,q)}(PSp(n,q_0))$ is $q^{\frac{1}{8}n^2}$. This implies
that $n-2\geq \frac{1}{8}n^2$ which is a contradiction for $n>4$.\label{spone}

Now suppose that $M$ is maximal in $PSp(4,q)$, $M\not\in\curlyc{1}$,
$|PSp(4,q):M|$ is odd and $|PSp(4,q):M|_p\leq q^2$. We must have
$M=(Sp(2,q)\circ Sp(2,q)).2$. Then $\Lalph\leq M$ and $\Lalph\geq P$ for
some $P$ a Sylow $p$-subgroup of $M$. Since the parabolic subgroups of
$Sp(2,q)$ have even index in $Sp(2,q)$ we must have
$\Lalph=(Sp(2,q)\circ Sp(2,q)).2$.

Now we can choose $\alpha$ such that
$$\Lalph= \hatt\left\langle\left(\begin{array}{cc} Sp(2,q)& \\ & Sp(2,q)
\end{array}\right), h^*:=\left(\begin{array}{cc} & I_{2\times 2} \\ I_{2\times 2} &
\end{array}\right)\right\rangle.$$


Observe that $h$ is conjugate to $g$ in PSp(4,q).
Now $h$ has at least
$\frac{1}{2}q^2(q^2-1)$ $\Lalph$-conjugates in $\Lalph$, thus
$\fracNR\leq \frac{\frac{1}{2}q^2(q^2+1)}{\frac{1}{2}q(q^2-1)}\leq
2q$. Then $v\leq 8q^2$. But $v>|L:\Lalph|=\frac{1}{2}q^2(q^2+1)$ which
is a contradiction for $q>3$.

Hence in all cases $M\in\curlyc{1}$.\label{sptwo, spextraone}
\end{proof}

In $\curlyc{1}$ we have subgroups of two types:
\begin{itemize}
\item       Parabolic subgroups, $P_m\cong
[q^a].(\frac{q-1}{(q-1,2)}).(PGL(m,q)\times PSp(n-2m,q))$ where $1\leq m\leq \frac{n}{2},
a=\frac{m}{2}-\frac{3m^2}{2}+mn.$ If $\Lalph\leq P_m$ then
$(q+1)\big||PSp(n,q):P_m|$ divides $v$. Hence we must have $p=2$.\label{spthree}

\item       Subgroups, $B_m$, of type $Sp_{m}\perp Sp_{n-m}$
isomorphic to $Sp(m,q)\circ Sp(n-m,q)$ where $2\leq m<\frac{n}{2}$ and
$m$ is even. In this case $q^2$ divides $|PSp(n,q):B_m|$ which in
turn divides $v.$ Hence we must have $p\equiv1(3)$.\label{spfour}
\end{itemize}

\subsection{Case: $p=2, \Lalph\leq P_m$}

The index of $P_m$ in $Sp(n,q)$ is divisible by $q^2+1$ for all $m>1$ which is
impossible and so $m=1$. Then $P_1\cong[q^{n-1}]:((q-1)\times Sp(n-2,q))$ and
$|Sp(n,q):P_1|=(q+1)(\polysumsq{n-2}).$ We conclude that $q\equiv
2(3)$ and that every prime dividing $\frac{n}{2}$ is congruent
to $1(3)$. Hence $n\geq 14$ and $n\equiv 2(4)$. This implies that $n-2\equiv 0(4)$ and every parabolic subgroup of $Sp(n-2,q)$ has index divisible by $q^2+1$. Thus $\Lalph= [q^{n-1}]:(A\times Sp(n-2,q))$ for some $A\leq q-1$.\label{spseven}

Now consider $Sp(n,q)$ acting on a vector space $V$ preserving a
symplectic form $\kappa$. For $u\in V, a\in GF(q)$ we have transvections in $Sp(n,q)$ defined by,
$$g_{a,u}:V\to V, v\mapsto v+a\kappa(v,u)u.$$
Set $W=\langle u\rangle$ and let $X_{W,W^\perp}=\{g_{a,u}:a\in
GF(q)\}$. Then $X_{W,W^\perp}<Sp(n,q)$ is of size $q$. The parabolic
subgroup of $Sp(n,q)$ which preserves $W$ normalizes
$X_{W,W^\perp}$.\label{speight}

Now let $g=g_{1,u}$. Then, since the
only maximal subgroup whose order is divisible by $\frac{|P_1|}{q-1}$
is $P_1$, we have
$$n_g\leq \frac{|Sp(n,q)|}{|P_1|}(q-1)\log_2q.$$
Similarly $r_g\geq \frac{|Sp(n-2,q)|}{|P_1^*|}$ where $P_1^*$ is a
parabolic subgroup of $Sp(n-2,q)$. Then
$$\fracNR\leq
\frac{|Sp(n,q)||P_1^*|(q-1)\log_2q}{|Sp(n-2,q)||P_1|}\leq q^4.$$
Thus $v\leq q^9$ which contradicts $n\geq 14$ and this case is excluded.\label{spnine}

\subsection{Case: $p\equiv 1(3), \Lalph<B_m$}

We know that the maximum power of $p$ in $v$ is at most $q^{n-2}$. Now $|PSp(n,q):B_m|_p=\frac{q^{\frac{1}{4}{n^2}}}{q^{\frac{1}{4}{m^2}}q^{\frac{1}{4}(n-m)^2}}$. Thus we need,
$$n-2\geq\frac{1}{4}(n^2-m^2-(n-m)^2)=\frac{1}{2}m(n-m).\label{spfive}$$
This implies that $m=2$ and so $\Lalph\leq Sp(2,q)\circ
Sp(n-2,q)$. If $n=4$ then $B_2$ is not maximal and so we assume that
$n>4$. Furthermore we know that $\Lalph$
must contain a Sylow $p$-subgroup of $Sp(2,q)\circ Sp(n-2,q)$. But the
indices of a parabolic subgroup of $Sp(2,q)$ in $Sp(2,q)$ and of a
parabolic subgroup of $Sp(n-2,q)$ in $Sp(n-2,q)$ are both divisible by
$q+1$, hence are even. Thus we conclude that $\Lalph=Sp(2,q)\circ
Sp(n-2,q)$.

Now $r_g\geq \frac{1}{2}q^{n-4}(q^{n-4}+\dots q^2+1)$ and so
$\fracNR\leq 2q^2(q^2+1)$ and $v\leq 8q^4(q^2+1)^2$. But
$v>|L:\Lalph|=q^{n-2}(q^{n-2}+\dots q^2+1)$ which is a contradiction
for $n>6$.

Thus we must assume that $n=6$ and $|L:\Lalph|=q^4(q^4+q^2+1)$ and
$\fracNR\leq 2q^2(q^2+1)$. If $|\fracNR|_p=|v|_p\geq q^4$ then
$\fracNR=q^4$ which contradicts Lemma \ref{lemma:primefixed}. Thus
$|\fracNR|_p=1$ and so $\fracNR\big|q^4+q^2+1$. If $\fracNR=q^4+q^2+1$
then $d_g$ is not divisible by $q^4$ which contradicts the fact that
$|L:\Lalph|$ divides $v$. If $\fracNR<\frac{1}{2}(q^4+q^2+1)$
then $v<|L:\Lalph|$ which is also a contradiction. \label{spsix, spextratwo}

\section{$L=\Omega(n,q)$, $nq$ odd}\label{section:orthogonalodd}

Throughout the next two sections, Greek letters such as $\epsilon, \eta$ and
$\zeta$ will stand for either $+,-$ or $\circ$. We will write polynomials
such as $x-\epsilon$ to mean $x-\epsilon1$. We write
$\Omega^\circ(n,q)$ to mean $\Omega(n,q)$ when $n$ is odd.


In this section we assume that $n\geq 7$ and $q$ is odd and we prove that, if $L=\Omega(n,q)$, then the hypothesis in Section \ref{section:hypothesis} leads to a contradiction. This implies the following proposition:

\begin{proposition}\label{proposition:orthogonalodd}
Suppose that $n$ is odd, $n\geq 7$ and $G$ has a minimal normal
subgroup isomorphic to $\Omega (n,q)$. Then $G$ does not act transitively on a projective plane.
\end{proposition}

Observe that $L$ contains $\Omega^\epsilon(n-1,q).2$ for $\epsilon=-$
and $\epsilon=+$. One of these groups contains a central involution and hence
$L$ contains an involution $g$ such that
$r_g(L)=\frac{1}{2}q^{\frac{n-1}{2}}(q^{\frac{n-1}{2}}+\epsilon)$.
Examining \cite[Table 3.5.D]{kl} for fusion of conjugacy classes, we
see that $n_g=r_g(L)$ and thus $|v|_p\leq
q^{\frac{n-1}{2}}$.\label{ooNINE}

We begin by proving that $\Lalph$ must lie in a maximal subgroup $M\in\curlyc{1}$:

\begin{lemma}\label{lemma:orthogonalodd}
$\Lalph$ does not lie inside a subgroup $M\in\curlyc{i}, i>1$.
\end{lemma}
\begin{proof}
We examine the list of odd index maximal subgroups in $G$ as given by
Liebeck and Saxl\cite{ls}. The following possibilities are available for a maximal subgroup $M$ of odd index. We exclude them in turn.
\begin{itemize}

\item       $L=\Omega(7,q)$ and $M=\Omega(7,2)$. We know that $|v|_p\leq
q^3$ and so $|\Lalph|$ must be divisible by $q^6$. This is impossible for $\Lalph\leq M$. \label{ooONE}

\item       $M\in\curlyc{2}$ or
$M=N_{\Omega(n,q)}(\Omega(n,q_0))$ where $q=q_0^c$
for $c$ an odd prime. In both cases $|M|_p\leq
\sqrt{|\Omega^\epsilon(n,q)|_p}$.\label{ooFOUR} Now $|\Omega^\epsilon(n,q)|_p=q^{\frac{1}{4}(n-1)^2}$ and so we must have,
$$\frac{1}{8}(n-1)^2+\frac{1}{2}(n-1)\geq\frac{1}{4}(n-1)^2.$$
This is impossible for $n\geq 7$.\label{ooFIVE}
\end{itemize}
\end{proof}

Thus $\Lalph$ lies inside a parabolic subgroup or a subgroup $B_m$ of
type $O(m,q)\perp O^\eta(n-m,q)$ for some odd $m<n$. In fact parabolic
subgroups have even index in $P\Omega(n,q)$ hence we may assume that $\Lalph\leq B_m$ for some $m$. \label{ooSIX}

Since $|v|_p\leq q^{\frac{n-1}{2}}$ we know that $\Lalph\leq
B_1=\Omega^\eta(n-1,q).2$ and that $\Lalph$ contains a Sylow
$p$-subgroup of $\Omega^\eta(n-1,q)$.\label{ooSEVEN} Now the parabolic
subgroups of $\Omega^\eta(n-1,q)$ have even index. Hence we must have
$\Lalph=\Omega^\eta(n-1,q)$ and $v$ is divisible by
$|\Omega(n,q):\Omega^\eta(n-1,q).2|=\frac{1}{2}q^{\frac{n-1}{2}}(q^{\frac{n-1}{2}}+\eta)$.

Now consider the involution $h$ centralized in $L$ by
$(\Omega^\zeta(2,q)\times\Omega(n-2,q)).[4].$ Then
$n_h=\frac{q^{n-2}(q^{n-1}-1)}{2(q-\zeta)}$. Now $\Omega^\eta(n-1,q)$
contains a conjugate of $h$ centralized by, at most,
$(\Omega^\zeta(2,q)\times \Omega^{\zeta\eta}(n-3,q)).[4]$. then
$r_h\geq
\frac{q^{n-3}(q^{\frac{n-3}{2}}+\eta\zeta)(q^{\frac{n-1}{2}}-\eta)}{2(q-\zeta)}$.
This implies that $\frac{n_h}{r_h}\leq q(q+1)$ and so $v\leq
2q^2(q+1)^2$. But then $v<|L:\Lalph|$ which is a
contradiction.\label{ooTHIRTEEN, ooCHECKONE}

Hence we have proved Proposition \ref{proposition:orthogonalodd}.


\section{$L=P\Omega^\epsilon(n,q)$, $n$ even}

In this section we assume that $n\geq 8$ and we prove that, if $L=P\Omega^\epsilon(n,q)$, then the hypothesis in Section \ref{section:hypothesis} leads to a contradiction. This implies the following proposition:

\begin{proposition}\label{proposition:orthogonaleven}
Suppose that $n$ is even, $n\geq8$ and $G$ has a minimal normal
subgroup isomorphic to $P\Omega^\epsilon(n,q)$. Then $G$ does not act transitively on a projective plane.
\end{proposition}

First we examine what happens when $p=2$:
\begin{lemma}
Suppose $n\geq 8$ is even and $G$ has a minimal normal subgroup isomorphic
to $P\Omega^\epsilon(n,2^a)$. Then $G$ does not act transitively on a set of size $x^2+x+1$.
\end{lemma}
\begin{proof}
Write $q=2^a$. We know that $\Lalph\leq P_m$ for some integer $m$. If
$m>1$ then $q^b+1$ divides $|P\Omega^\epsilon(n,q):P_m|$ where $b$ is
some even integer.\label{oeTEN} Since
$q^b+1\equiv 2(3)$ this is impossible. Thus $\Lalph$ lies inside some
parabolic subgroup $P_1$. Now
$$|P\Omega^\epsilon(n,q):P_1|=\frac{(q^{\frac{n}{2}}-\epsilon)(q^{\frac{n-2}{2}}+\epsilon)}{q-1}.$$

If $q\equiv2(3)$ then $q^\frac{n-2}{2}+1\equiv q^\frac{n}{2}+1\equiv2(3)$. Since one of these divides
 $|P\Omega^\epsilon(n,q):P_m|$, this is impossible. Hence $q\equiv
 1(3)$. Now let $n_e$ be the even one of $\frac{n}{2}$ and $\frac{n-2}{2}$, $n_o$
 the odd one. Then one of the following holds:
\begin{itemize}
\item   $|\Omega^\epsilon(n,q):P_1|=\frac{q^{n_e}-1}{q-1}(q^{n_0}+1)$
and $9$ divides $|\Omega^\epsilon(n,q):P_1|;$ or
\item   $|\Omega^\epsilon(n,q):P_1|=\frac{q^{n_o}-1}{q-1}(q^{n_e}+1)$ and $q^{n_e}+1\equiv 2(3).$
\end{itemize}
Both of these cases are impossible.\label{oeELEVEN}
\end{proof}

Throughout the rest of the section $p$ is odd. Now $L$ contains
maximal subgroups in $\curlyc{1}$ of type $O^\zeta(2,q)\perp
O^\eta(n-2,q)$ for
$\zeta\eta=\epsilon$. One of these groups contains a central
involution and hence $L$ contains an involution $g$ such that
$|L:C_L(g)|=\frac{q^{n-2}(q^{\frac{n-2}{2}}+\eta)(q^{\frac{n}{2}}-\epsilon)}{2(q-\zeta)}$.
\label{oeEIGHT} Examining for fusion of conjugacy classes in
\cite[Tables 3.5.E and 3.5.F]{kl} we see that, except when
$(n,\epsilon)=(8,+)$, $n_g=|L:C_L(g)|$. When $(n,\epsilon)=(8,+)$, we
know that $n_g\leq 3|L:C_L(g)|$ and so, in all cases, $|v|_p\leq q^{n-2}$.

We begin by proving that $\Lalph$ must lie in a maximal subgroup
$M\in\curlyc{1}$:

\begin{lemma}
$\Lalph$ does not lie inside a subgroup $M\in\curlyc{i}, i>1$.
\end{lemma}
\begin{proof}
We examine the list of odd index maximal subgroups in $G$ as given by
Liebeck and Saxl\cite{ls}. The following possibilities are available
for a maximal subgroup of odd index $M\not\in\curlyc{1}$. We exclude them in turn.
\begin{itemize}

\item       $L=P\Omega^+(8,q)$ and either $M=\Omega^+(8,2)$ or $M=2^3.2^6.PSL(3,2)$. We know that
$|v|_p\leq q^6$ and so $|\Lalph|_p\geq q^6$. This is impossible for
  $\Lalph\leq M$ in both cases.\label{oeTWO,oeTHREE}

\item       $M\in\curlyc{2}$ or
$M=N_{P\Omega^\epsilon(n,q)}(P\Omega^\epsilon(n,q_0))$ where $q=q_0^c$
for $c$ an odd prime. In both cases $|M|_p\leq
\sqrt{|P\Omega^\epsilon(n,q)|_p}$.\label{oeFOUR} Now $|P\Omega^\epsilon(n,q)|_p=q^{\frac{1}{4}n(n-2)}$ and so we must have
$$\frac{1}{8}n(n-2)+ n-2\geq\frac{1}{4}n(n-2).$$
This is impossible for $n>8$. When $n=8$, no subgroup $M$ of odd index
has $|M|_p\geq 6$ so the result stands.\label{oeFIVE}
\end{itemize}
\end{proof}

Thus $\Lalph$ lies inside a parabolic subgroup $P_m$ or a subgroup $B_m$ of
type $O(m,q)^{\zeta_1}\perp O^{\eta_1}(n-m,q)$ for some $m<\frac{n}{2}$. In fact parabolic
subgroups have even index in $P\Omega^\epsilon(n,q)$ for $p$ odd.\label{oeSIX}
Hence we assume that $\Lalph\leq B_m$ for some
integer $m$. We know that $|v|_p\leq q^{n-2}$ and so
$|P\Omega^\epsilon(n,q):B_m|_p\leq q^{n-2}$. This implies that $m=1$
or $m=2$. Note also that $p\equiv 1(3)$.\label{oeTWELVE}

Suppose first that $\Lalph\leq B_2$ where $B_2$ is of type
$O^{\zeta_1}(2,q)\perp O^{\eta_1}(n-2,q)$ for $\zeta_1\eta_1=\epsilon$. Then
$|P\Omega^\epsilon(n,q):B_2|=\frac{q^{n-2}(q^{\frac{n-2}{2}}+\eta_1)(q^{\frac{n}{2}}-\epsilon)}{2(q-\zeta_1)}$
and so $\Lalph$ must contain a Sylow $p$-subgroup of $B_2$. Since the
parabolic subgroups of $P\Omega^{\eta_1}(n-2,q)$ have even index we must
have $\Lalph>\Omega^{\eta_1}(n-2,q)$.\label{oeFOURTEEN}

In the case where $\Lalph\leq B_1$ then $\Lalph\leq \Omega(n-1,q).c_1$
where $c_1\in\{1,2\}$. Now
$|P\Omega^\epsilon(n,q):B_1|_p=q^{\frac{n-2}{2}}$ hence
$|B_1:\Lalph|_p\leq q^{\frac{n-2}{2}}$. Examining the proof of
Lemma \ref{lemma:orthogonalodd} this means that
$\Lalph\cap\Omega(n-1,q)$ lies inside a maximal subgroup of
$\Omega(n-1,q)$ in family $\curlyc{1}$. \label{oeFIFTEEN}

Since the parabolic subgroups of $\Omega(n-1,q)$ have even index in
$\Omega(n-1,q)$ this means that $\Lalph\cap\Omega(n-1,q)\leq
B_{m_1}^*$; here $B_{m_1}^*$ is a maximal subgroup of $\Omega(n-1,q)$ of type $O_{m_1}(q)\perp
O^\gamma(n-1-m_1,q)$ for some odd $m_1<n-1$. In fact $|B_1:\Lalph|_p\leq
q^{\frac{n-2}{2}}$ implies that $m_1=1$ and that $\Lalph$ contains a
Sylow $p$-subgroup of $B_1^*=\Omega^{\eta_1}(n-2,q).c_2$ where
$c_2\in\{1,2\}$. Once again, since the
parabolic subgroups of $\Omega^{\eta_1}(n-2,q)$ have even index we must
have $\Lalph>\Omega^{\eta_1}(n-2,q)$.\label{oeSIXTEEN}

Thus in both cases, when $m=1$ and when $m=2$, we see that
$\Lalph>\Omega^{\eta_1}(n-2,q)$ is a
subgroup of $P\Omega^\epsilon(n,q)$ which preserves a decomposition of
the associated vector space $V$ into subspaces, $V=W_2\perp
W_{n-2}$, where $\dim W_i=i$ and the $W_i$ are non-degenerate
subspaces of $V$.\label{oeSEVENTEEN}

Then $H=\Omega^{\eta_1}(n-2,q)$ contains $h$ a conjugate of $g$, and
$C_H(h)$ is isomorphic to either $(\Omega^{\gamma_1}(2,q)\times
\Omega^{\gamma_2}(n-4,q)).2$ or $2.(P\Omega^{\gamma_1}(2,q)\times P
\Omega^{\gamma_2}(n-4,q)).[4]$ (see \cite[Proposition 4.1.6]{kl}). In either case
$r_g\geq\frac{q^{n-4}(q^{\frac{n-4}{2}}+\gamma_2)(q^{\frac{n-2}{2}}-\eta_1)}{2(q-\gamma_1)}$.

If $n>8$ this means that $\fracNR\leq \frac{q^2(q+1)^3}{(q-1)^2}$ and
so $v\leq 2q^4(q+1)^4$. Since $|L:\Lalph|<v$ we must have $n=10, q=7$
and $\Lalph=B_1$. But then $|L:B_1|$ is divisible by
$\frac{1}{2}7^4(7^5\pm1)$. This is impossible since then $|L:B_1|$ is
divisible by a prime $s\equiv 2(3)$.\label{oechecktwo}

If $n=8$ then $\fracNR<4q^2(q+1)^2$.\label{oeEIGHTEEN} Then
$v<28q^4(q+1)^4$ which is less than $|L:B_2|$. Thus $\Lalph=B_1$. But
then $|L:\Lalph|$ is even which is a contradiction.\label{oecheckthree}

Proposition \ref{proposition:orthogonaleven} is now proven.\label{oeNINETEEN}

\section{$L$ is an exceptional group of Lie type in odd characteristic}\label{section:exceptionalodd}

In this section we prove that, if $L$ is an exceptional group of Lie type in odd characteristic, then the hypothesis in Section \ref{section:hypothesis} leads to a contradiction. This implies the following proposition:

\begin{proposition}\label{proposition:exceptionalodd}
Suppose that $G$ has a minimal normal subgroup $L$ where $L$ is an
exceptional group of Lie type in odd characteristic or that $G$ has a
unique component $L$ such that $L^\dag $ is isomorphic to a simple group
$E_6(q)$ or ${^2E_6(q)}$ where $q$ is odd. Then $G$ does not act
transitively on a projective plane.
\end{proposition}

We introduce some extra notation for this section and the following
one. We will write $E_6^-$ for ${^2E_6}$,
$E_6^+$ for $E_6$. Similarly $SL^-$ will stand for $SU$, $SL^+$ for
$SL$. We will use $\epsilon$ to denote either $\pm1$ or $\pm$
depending on the context. Generally our notation refers to the adjoint
version of the exceptional group, any variation on this will be
specified. For a group $G$, we will write $\frac{1}{2}G$ to
mean a subgroup in $G$ of index $2$. We define
$P(G):=\min\{{|G:H|}:H<G\}$. Finally, for a group $H$ we write
$O^{p'}H$ to mean the unique smallest normal subgroup $N$ of $H$ such
that $|H/N|_p=1$.

We have eight possibilities for $L$ which we will examine in
turn. As usual we will examine odd-index maximal subgroups of $L$,
treating these as candidates to contain a stabilizer $\Lalph$, and
seek to show a contradiction.

We immediately exclude the case where $L={^2G_2(q)}$, $q>3,$ by
examining the list of maximal subgroups  of ${^2G_2(q)}$ given in
\cite[Theorem C]{kleidman2} (see also \cite{ward}). We see that any maximal subgroup of odd
index must have index divisible by $9$ and hence cannot contain a
point-stabilizer. Hence this case is excluded.\label{eoONE} Note that
the list given by Kleidman \cite{kleidman2} contains a maximal
subgroup of odd index (with structure $(2^2\times
D_{\frac{1}{2}(q+1)}):3$) which has been omitted by Liebeck and
Saxl\cite{ls} and by Kantor\cite{kantor}.

For the remaining cases we will refer to the results of Liebeck and
Saxl giving the maximal subgroups $M^\dag $ of odd index in $L^\dag $.\cite{ls} These
maximal subgroups $M^\dag $ take one of two forms: Either
$M^\dag =N_{L^\dag }(L^\dag (q_0))$, where $q=q_0^a$ for $a$ an odd prime and the
subgroup $L^\dag (q_0)$ of $L^\dag (q)$ corresponds to the centralizer of a field
automorphism of $L^\dag (q)$ (see \cite[Theorem C]{kantor}), or $M^\dag $ is
enumerated in \cite[Table 1]{ls}.

Note that, by \cite[Table 5.1.B]{kl}, ${\rm Out} L$, the outer
automorphism group of $L$, has order strictly less than $q$ provided
$L\neq {^3D_4(3)}, {^2E_6(5)}$. We also use the following lemma:\label{eoTWO}

\begin{lemma}\label{lemma:inndiag}
Let $\phi$ be a field automorphism of $L(q)$ of prime order $a$. Let
$L(q_0)=O^{p'}C_{L(q)}(\phi)$ where $q=q_0^a$. Then
$N_{L(q)}(L(q_0))\lesssim {\rm Inndiag}(L(q_0))$ and, furthermore, ${\rm Inndiag}(L(q_0))=L(q_0).d$ where
\begin{displaymath}
d=\left\{ \begin{array}{ll}
(3,q_0-\epsilon) & L=E^\epsilon_6 \\
(2,q_0-1) & L=E_7 \\
1 & \rm{otherwise}
\end{array}\right.
\end{displaymath}
\end{lemma}
\begin{proof}
Our notation is consistent with that in \cite{gorenstein}. Write
$L(q)=O^{p'}C_{\overline{L}}(\sigma)$ where $\overline{L}$ is a
simple adjoint $\overline{\mathbb{F}_p}$-algebraic group,
$\overline{\mathbb{F}_p}$ is the algebraic closure of $GF(q)$ and $\sigma$ is
a Steinberg automorphism \cite[Definition 2.2.1]{gorenstein}.

By \cite[Proposition 2.5.17]{gorenstein}, there exists a Steinberg
automorphism $\tau$ of $\overline{L}$ such that $\tau^a=\sigma$ and
$\tau$ induces $\phi$ on $L$. Then
$L(q_o)=O^{p'}C_{\overline{L}}(\tau)$ and, by \cite[Proposition
2.5.9]{gorenstein}, $N_{\overline{L}}(L(q_0))= C_{\overline{L}}(\tau).$

Thus $N_{L(q)}(L(q_0))=C_{L(q)}(\tau)\leq C_{L(q)}(\phi)\lesssim {\rm
Inndiag}(L(q_0))$ by \cite[Proposition 4.9.1]{gorenstein}. The
structure of ${\rm Inndiag}(L(q_0))$ is given in \cite[Theorem
2.5.12]{gorenstein}. \label{eoTHREE}
\end{proof}

\subsection{Case: $L=E_8(q)$}

Referring to \cite[Table 4.5.1]{gorenstein}, we see that $E_8(q)$
contains an involution $g$ such that $C_L(g)\geq 2.(PSL(2,q)\times
E_7(q))$. There is one such $E_8(q)$ conjugacy class of involutions in $L$ and
so $n_g$ divides $$\frac{2q^{56}(q^{10}+1)(q^{12}+1)(q^6+1)(q^{30}-1)}{q^2-1}.$$\label{eoFOUR}

Using Lemma \ref{lemma:largestprime} this implies that $|v|_p\leq q^{56}$ and hence that
$|\Lalph|_p\geq q^{64}$. The list in \cite[Table 1]{ls} contains no
maximal subgroups $M$ such that $|M|_p\geq q^{64}$. Similarly Lemma
\ref{lemma:inndiag} implies that
$|N_L(E_8(q_0))|_p=|E_8(q_0)|_p=q_0^{120}$. Since $q=q_0^a$ where $a$
is an odd prime, $q_0^{120}\leq q^{40}$ and so this possibility is excluded.\label{eoFIVE}

\subsection{Case: $L=E_7(q)$}

Referring to \cite[Table 4.5.1]{gorenstein}, we see that  $E_7(q)$
contains an involution $g$ such that $C_L(g)$ contains
$SL^\epsilon(8,q)/(4,q-\epsilon)$ for $\epsilon$ either $+$ or
$-$. There is one such ${\rm Inndiag}(E_7(q))$ conjugacy class
of involutions in $L$ and so $n_g$ divides $$(4,q-1)q^{35}(q^7+\epsilon)(q^5+\epsilon)(q^3+\epsilon)(q^8+q^4+1)(q^{12}+q^6+1).$$\label{eoSIX}

This implies that $|v|_p\leq q^{35}$ and hence that
$|\Lalph|_p\geq q^{28}$. The list in \cite[Table 1]{ls} contains one
maximal subgroup such that $|M|_p\geq q^{28}$, namely $M=N_L(2.(PSL(2,q)\times P\Omega^+(12,q))$. Then $|L:M|_p=q^{32}$
and so $p\equiv 1(3)$. But this implies that $9$ divides $|L:M|$
and so it is not possible that $\Lalph\leq M$.\label{eoSEVEN}

Similarly Lemma \ref{lemma:inndiag} implies that
$|N_L(E_7(q_0))|_p\leq|E_7(q_0).2|_p=q_0^{63}$. Since $q=q_0^a$ where $a$
is an odd prime, $q_0^{63}\leq q^{21}$ and so this possibility is excluded.\label{eoEIGHT}

\subsection{Case: $L^\dag =E_6^\epsilon(q)$}

Referring to \cite[Table 4.5.1]{gorenstein}, we see that $L$
contains an involution $g$ such that $C_L(g)$ contains
$Spin_{10}^\epsilon(q)$. Here $Spin^\epsilon_{10}(q)\cong
(4,q-\epsilon).P\Omega^\epsilon(10,q)$. There is only one such ${\rm
Inndiag}(E_6^\epsilon(q))$ conjugacy class of
involutions in $L$ and so,
$$n_g=q^{16}(q^6+\epsilon q^3+1)(q^2+\epsilon q+1)(q^8+q^4+1).$$\label{eoNINE}

This implies that $|v|_p\leq q^{16}$ and hence that
$|\Lalph|_p\geq q^{20}$. Then Lemma \ref{lemma:inndiag} implies that
$|N_{L^\dag }(L^\dag (q_0))|_p\leq|L^\dag (q_0).(3,q-\epsilon)|_p$ which
divides $3q_0^{36}$. Since $q=q_0^a$ where $a$
is an odd prime, $q_0^{36}\leq q^{12}$ and so this possibility is excluded.\label{eoTEN}

\subsubsection{Subcase: $\epsilon=+$}

In this case the list in \cite[Table 1]{ls} contains two maximal
subgroups $M^\dag $ such that $|M^\dag |_p\geq q^{20}$:
$M^\dag =N_{L^\dag }((4,q-1).P\Omega^+(10,q))$ or $M^\dag $ is parabolic of type
$D_5$. If $p\equiv 1(3)$ in either case then $9$ divides $|L:M|$ which
is a contradiction. Hence $p\not\equiv 1(3)$, the universal and
adjoint versions coincide and $L$ is simple.

In the non-parabolic case,  $|L:M|_p>p^2$ which is impossible for
 $p\not\equiv 1(3)$. Hence $M$ is a parabolic subgroup of $E_6^+(q)$ of
 type $D_5$ and $|L:M|=(q^6+q^3+1)(q^2+q+1)(q^8+q^4+1)$.\label{eoELEVEN}

Now  $M\cong [q^{16}]: (Spin^+_{10}(q)H)$ where $H$ is a Cartan
subgroup of $E_6(q)$ and $H$ normalizes $Spin^+_{10}(q)$. Here
$Spin^+_{10}(q)\cong (4,q-1).P\Omega^+(10,q)$ and $P\Omega^+(10,q)$
has parabolic subgroups of even index. This implies that $\Lalph\geq
[q^{16}]:(Spin^+_{10}(q).2)$ for $p\neq 3$.

Furthermore, for $p=3$, every non-parabolic subgroup of
$P\Omega^+(10,q)$ has index divisible by $9$\cite{kl}. This
means that $\Lalph\geq [\frac{q^{16}}{3}].(Spin^+_{10}(q).2)$. Now $E$, the
commutator subgroup of the Levi complement in $M,$ is isomorphic to
$Spin^+_{10}(q)$ and $|E:\Lalph\cap E|$ is at most $\frac{3}{2}(q-1)$. But
$P(Spin^+_{10}(q))>\frac{3}{2}(q-1)$ \cite[Table 5.2.A]{kl}. Thus
$\Lalph>E$.

Now if $q=3^a$ then $|E|$ is divisible by $3^{8a}-1;$ in particular, $|E|$ is
divisible by the primitive prime divisors of $3^{8a}-1.$ This implies that if
$\phi:E\to GL(m,3)$ is a non-trivial representation of $E$ over $GF(3)$ then
$m\geq 8a$. Now consider the action of $E$ on the unipotent radical of
the full parabolic group, $[q^{16}],$ considered as a module over
$GF(3)$. We know that $E$ does not act trivially on any submodule of the
unipotent radical (otherwise $Z(E)$ would have too large a
centralizer; see \cite[Table 4.5.1]{gorenstein}). Thus the action must be
either irreducible or split into two modules both of size $q^8$. In
either case we must have $\Lalph\geq
[q^{16}]:(Spin^+_{10}(q).2)$.\label{eoTWELVE, eoTWELVEB}

We return to the general case where $p\not\equiv 1(3)$ and assume that $M$ contains $C_L(g)=Spin^+_{10}(q)H$. Furthermore we
know that $L$ acts on the cosets of $M$ as a rank 3 permutation group with
subdegrees $1, q(q^3+1)(q^8-1)/(q-1)$ and
$q^8(q^4+1)(q^5-1)/(q-1)$(\cite{kantor}). Then we have two possibilities:\label{eoTHIRTEEN}

\begin{itemize}
\item   Suppose $C_{M}(h)\geq Spin^+_{10}(q)$ for all $h$ in $\Lalph$ where
$h$ is $L$-conjugate to $g$. Now if $M=[q^{16}]:C_L(g)$ then $M$ contains
$q^{16}$ $M$-conjugates of $C_L(g)$ each containing a unique copy of
$Spin_{10}^+(q)$. Any other $L$-conjugate of $C_L(g)$ lies inside a
non-trivial conjugate of $M$. But these intersect $M$ with non-trivial
indices as above. These intersections cannot contain
$Spin_{10}^+(q)$. Hence $M$ contains only $M$-conjugates of $g$ and,
in fact, all these must lie in $\Lalph$. Thus $r_g=q^{16}$ and
$\fracNR=(q^8+q^4+1)(q^6+q^3+1)(q^2+q+1)$. Set
$$u=q^8+\frac{1}{2}q^7+\frac{3}{8}q^6+\frac{5}{16}q^5\frac{99}{128}q^4+\frac{127}{256}q^3+\frac{423}{1024}q^2+\frac{749}{2048}q+\frac{39587}{32768}.$$

Then $u^2-u+1>\fracNR$ for $q\geq 47$. If we set
$u_1=u-\frac{1}{32768}$ then $u_1^2-u_1+1<\fracNR$ for $q>1$. Thus we
need to check $q<47$ but no such $q$ satisfies $u^2-u+1=\fracNR$ for
integer $u$.\label{e0CHECKONE}


\item   Suppose there exists $h$ in $\Lalph$ which is $L$-conjugate to
$g$ and $C_{M}(h)$ does not contain a copy of $Spin^+_{10}(q)$. Then
$C_L(h)$ lies inside a non-trivial conjugate of $M$. Hence
$|M:C_M(h)|$ is a multiple of $q(q^3+1)(q^8-1)/(q-1)$ or
$q^8(q^4+1)(q^5-1)/(q-1)$. Furthermore we know that $q^{16}$ divides
$|M:C_M(h)|$ since $|M|_p=q^{16}|C_L(g)|_p$. Hence $|M:C_M(h)|\geq q^{16}(q^4+1)(q^5-1)/(q-1)$.


Now, if $\Lalph\geq [q^{16}]: (Spin^+_{10}(q).2)$ then $r_g=r_g(M)$
since $\Lalph\unlhd M$ and $|M:\Lalph|$ is odd. Thus
$r_g\geq q^{16}(q^4+1)(q^5-1)/(q-1)$ and $\fracNR< q^8+q^4+1$. Then
$d_g\leq q^8+q^4+1<(q^6+q^3+1)(q^2+q+1)$. Thus
$v<|L:M|$ which is a contradiction.\label{eoFIFTEEN}

\end{itemize}

\subsubsection{Subcase: $\epsilon=-$}

In this case the list in \cite[Table 1]{ls} contains one maximal
subgroup $M^\dag $ in $L^\dag $ such that $|M^\dag |_p\geq q^{20}$, namely
$M^\dag =N_{L^\dag }((4,q+1).P\Omega^-(10,q))$. In fact $|M|_p=q^{20}$ and so
$p \equiv 1(3)$ and the universal and adjoint versions of $E_6^-$ coincide and
$L$ is simple. Then
$M=N_L(Spin^-_{10}(q))\cong Spin^-_{10}(q).(q+1)$ (\cite[Table
4.5.2]{gorenstein}). Furthermore $\Lalph$ must contain a Sylow $p$-subgroup of
$M$. But the parabolic subgroups of $P\Omega^-_{10}(q)$ have even
index, hence $Spin^-_{10}(q).2\leq \Lalph\leq
Spin^-_{10}(q).(q+1)$.\label{eoSIXTEEN}

Now, using \cite[Table 4.5.2]{gorenstein}, we see that $E_6^-(q)$
contains two conjugacy classes of involutions: those conjugate to $g$,
centralized by $Spin^-_{10}(q)$, and those conjugate to $g_1$ say,
centralized by $SL(2,q)\circ SU(6,q)$. Then
$n_g=q^{16}(q^2-q+1)(q^6-q^3+1)(q^8+q^4+1)$ and
$N_{g_1}=q^{20}(q^4+1)(q^2+1)(q^6-q^3+1)(q^8+q^4+1)$.\label{eoSEVENTEEN}

We examine the involutions lying in $Spin^-_{10}(q)$ using \cite[Table
4.5.2]{gorenstein}. Apart from the central involution,
$Spin^-_{10}(q)$ contains two conjugacy classes of involutions. Let $h$ be an
involution in $Spin^-_{10}(q)$ centralized by
$Spin^+_4(q)\circ Spin^-_6(q)$. Then $\Lalph$ contains at least
$\frac{1}{4}q^{12}(q^4+q^3+q^2+q+1)(q^2-q+1)(q^4+1)(q^2+1)$ conjugates
of $h$. If $h$ is $L$-conjugate to $g$, then $\fracNR<4q^8$ which is a
contradiction. Thus assume that $h$ is $L$-conjugate to $g_1$.\label{eoEIGHTEEN}

In this case $\fracNR\leq 4q^{16}+4q^{12}+4q^8$. Then
$$d_g<\fracNR+2\sqrt{\fracNR}+2<4q^{16}+4q^{12}+6q^8+2q^4+2.$$
This implies that $v<19|L:M|$ for $q\geq7$.

Now suppose that $q^{16}$ does not divide $\fracNR$. Then
 $\fracNR$ divides $(q^2-q+1)(q^6-q^3+1)(q^8+q^4+1)$ and so
 $d_g<3q^{16}$ and $v=|L:M|$. This contradicts Lemma
 \ref{lemma:kantorinequality}. Thus $v=7|L:M|$ or $v=13|L:M|$ and
 $q^{16}\big| \fracNR.$

If $\fracNR \geq 7q^{16}$ then $v>49q^{32}>13|L:M|$ which is a
 contradiction. Thus, by Lemma \ref{lemma:primefixed},
 $\fracNR=3q^{16}$. This implies that
 $3q^{16}<d_g<3q^{16}+2\sqrt{3}q^8+2$ and so
 $9q^{32}<v<9q^{32}+12q^{24}+6q^{16}$. But then $7|L:M|<v<13|L:M|$
 which is a contradiction.\label{eoNINETEEN}

\subsection{Case: $L={^3D_4(q)}$}

We know that ${^3D_4(q)}$ has a single conjugacy class of
involutions\cite{gorenstein} which is centralized by a maximal
subgroup isomorphic to $(SL(2,q^3)\circ SL(2,q)).2$
\cite{kleidman3}. Hence, for $g$ an involution in $L$,
$n_g=q^8(q^8+q^4+1)$ and so $|v|_p\leq q^8$ and $|\Lalph|_p\geq
q^4$.\label{eoTWENTY}

If $\Lalph<M=N_L({^3D_4(q_0)}))$ then this condition implies that
$q=q_0^3$. No such subfield subgroup exists.\label{eoTWENTY-ONE}

There are two other odd index maximal subgroups $M$ such that
$|M|_p\geq q^4$.\cite{ls} The first possibility is that $M=G_2(q)$ and
$|L:M|_p=q^6$. But then odd index subgroups of $G_2(q)$ have $p$-index
strictly greater than $q^2$.\cite{ls} Thus $\Lalph=G_2(q)$. Now
$r_g(G_2(q))=q^4(q^4+q^2+1)$ and so $\fracNR=q^4(q^4-q^2+1)$. But this
implies that $|v|_p\leq q^4$ which is impossible.\label{eoTWENTY-TWO}

The second possibility is that $\Lalph\leq M=2.(PSL(2,q)\times PSL(2,q^3)).2$. Then
$|L:M|=q^8(q^8+q^4+1)$ and so $p\equiv 1(3)$ and $\Lalph$ contains a
Sylow $p$-subgroup of $M$. But the parabolic subgroups of $PSL(2,q)$
have even index, hence we conclude that $\Lalph=M$.\label{eoTWENTY-THREE}

Now $r_g(2.(PSL(2,q)\times
PSL(2,q^3)))\geq1+\frac{1}{2}q^3(q^3-1)\frac{1}{2}q(q-1)$. This
implies that $\fracNR<7q^8$. Suppose that
$|\fracNR|_p=1$ and hence $\fracNR\leq q^8+q^4+1$. Then $d_g<3q^8$ and
so $d_g=q^8$. This contradicts Lemma \ref{lemma:kantorinequality}.\label{eoTWENTY-FOUR}

Thus $|\fracNR|_p>1$ and so we must have either $\fracNR=q^8$ (contradicting
Lemma \ref{lemma:primefixed}) or $\fracNR=3q^8$. If $\fracNR=3q^8$
then $d_g<\frac{13}{3}(q^8+q^4+1)$ which is the smallest possibility for
$d_g$ that is larger than $\fracNR$. Thus we have a contradiction.\label{eoTWENTY-FIVE}

\subsection{Case: $L=G_2(q)$}

Referring to \cite[Table 4.5.1]{gorenstein}, we see that  $G_2(q)$
contains an involution $g$ such that $C_L(g)$ contains
$SL(2,q)\circ SL(2,q)$. There is one such conjugacy class of involutions in $L$
and, examining \cite{kleidman2}, we see that $C_L(g)\cong (SL(2,q)\circ
SL(2,q)).2.$ Hence $n_g=q^4(q^4+q^2+1).$ Using Lemma
\ref{lemma:largestprime}, we may conclude that $|v|_p\leq q^4$ and
hence that $|\Lalph|_p>q^{2}$.\label{eoTWENTY-SIX}

Examining the odd-index maximal subgroups \cite{kl}, we find that all
have $p$-index divisible by $p^2$ and so $p\equiv 1(3)$. We have a number of
possibilities for $M$ an odd-index maximal subgroup, $|M|_p\geq q^2$, $M$ containing $\Lalph$:\label{eoTWENTY-SEVEN}

\begin{itemize}

\item Suppose $M=N_L(G_2(q_0))$. Then using Lemma \ref{lemma:inndiag} we find
that $q=q_0^3$. But this means that $9$ divides $|L:M|$ which is
impossible.\label{eoTWENTY-NINE}

\item Suppose $M=(SL(2,q)\circ SL(2,q)).2.$ Then $\Lalph\geq2.P.2$ where
$P$ is a Sylow $p$-subgroup of $PSL(2,q)\times PSL(2,q)$. Since the
parabolic subgroup of $PSL(2,q)$ have even index we must have
$\Lalph=M$ and $v=q^4(q^4+q^2+1)a$ for some integer $a$. Then
Lemma \ref{lemma:kantorinequality} implies that $a\not=1$ and so
$a\geq7$.\label{eoTHIRTY}

Now $PSL(2,q)\times PSL(2,q)$ has at least $\frac{1}{4}q^2(q\pm 1)^2$
involutions and thus so does $SL(2,q)\circ SL(2,q)$. Then
$$\fracNR<4q^2\frac{q^4+q^2+1}{q^2-2q+1}<7q^4$$
for $q\geq 7$. Thus either $\fracNR=q^4$ (contradicting Lemma
\ref{lemma:primefixed}) or $\fracNR=3q^4$ or $\fracNR$ divides $q^4+q^2+1$.\label{eoTHIRTY-ONE}

If $u^2-u+1=\fracNR=3q^4$ then
$u^2+u+1=d_g<3q^4+2\sqrt{3q^4}+2<4q^4+4q^2+4$. This implies that
$v<12q^4(q^4+q^2+1)$ and so $a=7$. But then
$d_g=\frac{7}{3}(q^4+q^2+1)$ which is less than $\fracNR$ for $q\geq
7$. This is impossible.\label{eoTHIRTY-TWO}

If $u^2-u+1=\fracNR=q^4+q^2+1$ then $u=q^2+1$ and
$d_g=q^4+3q^2+3$. But then $(v,p)=1$ which is impossible. If
$\fracNR<q^4+q^2+1$ then $u\leq q^2$ which implies that $\fracNR\leq
q^4-q^2+1$ and $d_g\leq q^4+q^2+1$. Then $\fracNR d_g<|L:M|$ which is a contradiction.\label{eoTHIRTY-THREE}

\item Suppose $M=SL^\epsilon(3,q).2$ and so $p\equiv 1(3)$. Consider
first the situation where $\Lalph=M$. When $\epsilon=+$,
$M=\langle SL(3,q),\phi\rangle$ where $\phi$ is a graph automorphism
\cite[(2.6)]{chang}. When $\epsilon=-,$ $M\leq P\Gamma U(3,q)$
\cite[Proposition 2.2]{kleidman2}. In both cases $M$ is equal to a
universal version of $A_2^\epsilon(q)$ extended by a graph
automorphism \cite[Definition 2.5.13]{gorenstein}.

Examining \cite[Table 4.5.2]{gorenstein} we see that $M$ has 2
conjugacy classes of involutions. These have size $q^2(q^2+\epsilon
q+1)$ and $q^2(q^2+\epsilon q+1)(q-\epsilon).$ When $\epsilon=+$ this
gives $r_g=q^3(q^2+q+1)$ and $\fracNR=q(q^2-q+1)$. This is impossible
since either $|\fracNR|_p=1$ or $|\fracNR|_p\geq q^3$. When $\epsilon = -$ we have
$r_g=q^2(q^2-q+1)(q+2)$ and $\fracNR=\frac{q^2(q^2+q+1)}{q+2}$. This is
not an integer for $q>1$ hence can be excluded.\label{eoTHIRTY-FIVE}

Thus we must have $\Lalph<M$ and we know that  $|M:\Lalph|_p\leq
q$. Examining the subgroups of $SL^\epsilon(3,q)$ we find that
$\Lalph\cap SL^\epsilon(3,q)\leq P_1,$ a parabolic subgroup of
$SL^\epsilon(3,q)$.

When $\epsilon=-$, $|SL^\epsilon(3,q):P_1|$ is even hence this
possibility can be excluded.\label{eoTHIRTY-FOUR, eoTHIRTY-SIX}

When $\epsilon=+$, $M=\langle SL(3,q),m\rangle$ where $m$ is a graph automorphism
of $SL(3,q)$. Since $\Lalph$ has odd index in $G_2(q)$, $\Lalph$ must
contain a graph automorphism. Examining \cite[Table 3.5.A]{kl} we find that
$\Lalph\cap SL(3,q)$ lies inside a subgroup $M_1$ of $SL(3,q)$ of type
$GL(2,q)\oplus GL(1,q)$ or of type $P_{1,2}$. In the former case we find that
$|v|_p\geq q^5$. Since $|n_g|_p=q^4$ we must have $|\fracNR|_p=1$
which implies that $\fracNR\leq q^4+q^2+1$ and $|d_g|_p\geq q^5$ which
contradicts Lemma \ref{lemma:counting}. In the latter case, we find
that $|SL(3,q):M_1|$ is even and this case can be excluded.\label{eoTHIRTY-EIGHT}
\end{itemize}

We have covered all possible odd-index maximal subgroups in $G_2(q)$.

\subsection{Case: $L=F_4(q)$}

Referring to \cite[Table 4.5.1]{gorenstein}, we see that  $F_4(q)$
contains an involution $g$ such that $C_L(g)$ contains
$Spin(9,q)$. There is one such conjugacy class of involutions in $L$
and so $n_g=q^8(q^8+q^4+1).$\label{eoTHIRTY-NINE}

This implies that $|v|_p\leq q^{8}$ and hence that
$|\Lalph|_p\geq q^{16}$. Then Lemma \ref{lemma:inndiag} implies that
$|N_L(F_4(q_0))|_p\leq|F_4(q_0)|_p=q_0^{24}$. Since $q=q_0^a$ where $a$
is an odd prime, $q_0^{24}\leq q^{8}$ and so $\Lalph$ does not lie in
$|N_L(F_4(q_0))$.\label{eoFORTY}

The list in \cite[Table 1]{ls} contains one maximal subgroup $M$ such
that $|M|_p\geq q^{16}$. Then $M\cong 2.\Omega(9,q)$, $\Lalph$ must
contain a Sylow $p$-subgroup of $M$ since
$|L:M|_p=q^{16}$. Furthermore, $p\equiv 1(3)$. Now the
parabolic subgroups of $\Omega(9,q)$ have even index, hence we may
conclude that $\Lalph=M$ and $v=q^8(q^8+q^4+1)a$ for some integer
$a$. Lemma \ref{lemma:kantorinequality} implies that $a\not=1$ and
hence $a\geq 7$.\label{eoFORTY-ONE}

Now suppose $r_g\geq\frac{1}{2}q^4(q^4-1)$. Then $\fracNR\leq
2q^4(q^4+3)<\frac{7}{3}q^8$. Then $d_g<\frac{14}{3}q^8$ and
$v<7q^{16}$ which is a contradiction. Also $r_g$ is clearly greater
than $1$. Thus there is an involution $g\in 2.\Omega(9,q)$ such that
$$1<|2.\Omega(9,q):C_{2.\Omega(9,q)}(g)|<\frac{1}{2}q^4(q^4-1).$$\label{eoFORTY-TWO}

Now let $B$ be the central subgroup of $\Lalph$ of order $2$, so that
$\Lalph/B\cong P\Omega(9,q)$. Let $h=Bg$ an involution in $P\Omega(9,q)$. Then we must have
$$|\Omega(9,q):C_{\Omega(9,q)}(h)|<\frac{1}{2}q^4(q^4-1).$$
Examining \cite[Table 4.5.1]{gorenstein} we see that all involution
centralizers in $\Omega(9,q)$ have index at least
$\frac{1}{2}q^4(q^4-1)$. Hence we have a contradiction.\label{eoFORTY-THREE}

Proposition \ref{proposition:exceptionalodd} is now proven.

\section{$L$ is an exceptional group of Lie type in characteristic $2$}

In this section we prove that, if $L$ is an exceptional group of Lie type in characteristic $2$, then the hypothesis in Section \ref{section:hypothesis} leads to a contradiction. This implies the following proposition:

\begin{proposition}\label{proposition:exceptionaleven}
Suppose $G$ has a minimal normal subgroup $L$ where $L$ is an
exceptional group of Lie type in characteristic $2$ or that $G$ has a
unique component $L$ such that $L^\dag $ is isomorphic to $E_6(q)$ or
${^2E_6(q)}$ where $q=2^a$. Then $G$ does not act transitively on a
projective plane.
\end{proposition}

We have nine possibilities for $L$ and, by Tits' Lemma
\cite[1.6]{seitz}, we know that $\Lalph$ must lie in a parabolic
subgroup $M$ of $L$. We demonstrate that this is impossible, generally
by showing a contradiction with Lemma \ref{lemma:projsyl}.\label{eONE}

\subsection{Case: $L={^3D_4(q)}$; ${G_2(q)}, q>2$ }
In each case, for any parabolic subgroup $M$, $|L:M|$ is
divisible by $(q^4+q^2+1)(q+1)$. If $q\equiv 1(3)$ then $|L:M|$ is
divisible by $q+1\equiv 2(3)$, while if $q\equiv 2(3)$ then $9$
divides $|L:M|$. Thus $M$ cannot contain $\Lalph$ (Lemma
\ref{lemma:projsyl}) and we are done.\label{eTHREE}

\subsection{Case: $L={^2B_2(q)}, q>2; {^2F_4(q)}', F_4(q), E_7(q), E_8(q)$}\label{subsection:exceptionalevenb}

Examining the indices of the parabolic subgroups $M$ in $L$ in these
cases, we find that they are nearly always divisible by $q^m+1$ for
some even integer $m$. Since $q^m+1\equiv 2(3)$ these cases are
excluded. We deal with the exceptions which are as follows:\label{eFOUR}
\begin{enumerate}

\item $L=E_7(q)$ and $M$ is of type $E_6$. Then $|L:M|$ is divisible by
$(q^5+1)(q^9+1)$. If $q\equiv 1(3)$ then $q^5+1\equiv
2(3)$ and if $q\equiv 2(3)$ then $9$ divides $|L:M|$. Both of these
are impossible hence $M$ cannot contain $\Lalph$.\label{eSIX}

\item $L=E_7(q)$ and $M$ is of type $D_6$. Then $|L:M|$ is divisible by
$(q^8+q^4+1)(q^{12}+q^6+1)$ which is in turn divisible by $9$. Hence
$M$ cannot contain $\Lalph$.\label{eSEVEN}

\item $L=E_7(q)$ and $M$ is of type $D_5\times A_1$. Then $|L:M|$ is
divisible by $(q^5+1)(q^8+q^4+1)$. If $q\equiv 1(3)$ then $q^5+1\equiv
2(3)$ and if $q\equiv 2(3)$ then $9$ divides $|L:M|$. Both of these
are impossible hence $M$ cannot contain $\Lalph$.\label{eEIGHT}
\end{enumerate}

Note that Kantor's argument to exclude the last two cases ($L=
E_7(q)$ and $M$ of type $D_6$ or $D_5\times A_1$) when the action is
primitive is incorrect\cite{kantor}.\label{eNINE}

\subsection{Case: $L^\dag =E_6^\epsilon(q)$}

We proceed as in Subsection
\ref{subsection:exceptionalevenb}; we need only examine the parabolic
subgroups $M$ in $L$ which are not divisible by $q^m+1$ for some even
integer $m$. There are two such possibilities:

\begin{enumerate}
\item   $L^\dag =E_6^+(q)$ and $M$ is of type $D_5$. Then $|L:M|=(q^6+q^3+1)(q^8+q^4+1)(q^2+q+1)$. For $q\equiv 1(3)$, $|L:M|$ is
divisible by $9$ hence $M$ cannot contain $\Lalph$. Thus we assume
that $q\equiv 2(3)$ and so $L$ is simple.\label{eFIVE}

Now we know that $M':=[q^{16}].\Omega_{10}^+(q)\leq
\Lalph\leq M\cong [q^{16}]:(\Omega_{10}^+(q) H)$ where $H$ is the Cartan
subgroup of $L$. This is because all parabolic subgroups of
$\Omega_{10}^+(q)$ have index divisible by $q^4+1\equiv 2(3)$.

By \cite[(15.1),(15.5)]{aschseitz}, $L$ contains an involution $g$
such that $C_L(g)=[q^{21}]:SL(6,q)$ and so
$n_g=(q^6+q^3+1)(q^8+q^4+1)(q^8-1).$ Now if $r_g\geq
(q^6+q^3+1)(q^8-1)$ then $\fracNR\leq (q^4+1)^2-(q^4+1)+1$ and so
$d_g\leq (q^4+1)^2+(q^4+1)+1$. But then $\fracNR d_g< |L:M|$ which is
a contradiction. Thus, for all $h\in \Lalph$ conjugate in $G$ to $g$,
$|K:C_K(h)|<(q^6+q^3+1)(q^8-1)$. \label{eELEVEN}

Now $\Omega_{10}^+(q)\not\leq C_L(g)$. Furthermore the only maximal
subgroups of $\Omega_{10}^+(q)$ with index less than
$(q^6+q^3+1)(q^8-1)$ are the parabolic subgroups and $Sp_8(q)$. All
but one type of parabolic subgroups have index
divisible by $q^3+1$. Since $q^3+1$ does not divide $n_g$, there
must be $h\in \Lalph$ conjugate in $G$ to $g$ such that
$C_{K}(h)$ lies in either
$[q^{16}].([q^8]:\frac{1}{2}((q-1)\times SO_8^+(q)))$ or
$[q^{16}].Sp_8(q)$.\label{eTWELVE}

Consider the first possibility. Now $SO_8^+(q)\not\leq C_L(g)$ and so
$$r_g\geq P(SO_8^+(q))\frac{|\Omega_{10}^+(q)|}{|[q^8]:\frac{1}{2}((q-1)\times SO_8^+(q))|}.$$
Using the value for $P(SO_8^+(q))$ given in \cite[Table 5.2.A]{kl} we
conclude that $r_g>(q^6+q^3+1)(q^8-1)$ which is impossible.\label{eTHIRTEEN}

Similarly $Sp_8^+(q)\not\leq C_L(g)$ and so
$$r_g\geq P(Sp_8^+(q))\frac{|\Omega_{10}^+(q)|}{|Sp_8^+(q))|}.$$
Once again we find that $r_g>(q^6+q^3+1)(q^8-1)$ which is
impossible.\label{eFOURTEEN}

\item $L^\dag =E_6^-(q)$ and $M$ is of type ${^2D_4}(q)$. Then $|L:M|$ is
divisible by $(q^5+1)(q^9+1)$; we exclude this possibility in the
same way as in Subsection \ref{subsection:exceptionalevenb}, when $L=E_7(q)$ and $M$ is of type $E_6$.\label{eSIXb}
\end{enumerate}

Theorem A is now proven. 

\bibliographystyle{amsalpha}
\bibliography{paper}

\end{document}